\documentstyle[12pt]{article}
\begin{document}
\title{Induced expansion for quadratic polynomials}
\author{Jacek Graczyk  \\
Institute of Mathematics, Warsaw University, \\
ul.Banacha 2, 02-097 Warszawa, Poland.
\and
Grzegorz \'{S}wi\c{a}tek
\thanks{Partially supported by NSF grant \#431-3604A and the Sloan 
Foundation} \\
Institute for Mathematical Sciences,\\
SUNY at Stony Brook,\\
Stony Brook, NY 11794, USA.}
\date{August 17, 1993}

\newtheorem{lem}{Lemma}[section]
\newtheorem{conjec}{Conjecture}
\newtheorem{theo}{Theorem}
\newtheorem{coro}{Corollary}
\newtheorem{com}{Comment}

\newtheorem{prop}{Proposition}
\newtheorem{prob}{Problem}
\newtheorem{con}{Construction}[section]
\newtheorem{defi}{Definition}[section]
\newcommand{\hf}{\hat{f}}
\newcommand{\A}{{\cal A}}
\newcommand{\B}{{\cal B}}
\newcommand{\de}{{\bf \delta}}
\newcommand{\sr}{\mbox{{{\em I\hskip -2pt  R}}}}
\newcommand{\sn}{\mbox{{{\em  N\hskip -2pt I}}}}
\newtheorem{fact}{Fact}[section]
\newenvironment{proof}
{{\bf Proof:}\newline}{\begin{flushright}$\Box$\end{flushright}}
\newenvironment{sol}
{{\bf Solution:}\newline}{\begin{flushright}
$\Box$\end{flushright}}
\newcommand{\ex}{{\cal EX}}
\newcommand{\Cr}{{\bf Cr}}
\newcommand{\Bo}{\Box^{n}_{i}}
\newcommand{\D}{{\cal D}} 
\newcommand{\R}{{\bf R}}
\newcommand{\DPo}{{\bf DPoin}}
\newcommand{\td}{{\underline \tau}}
\newcommand{\tg}{{\overline \tau}}
\newcommand{\gd}{{\underline \gamma}}
\newcommand{\gu}{{\overline \gamma}}
\newcommand{\al}{{\overline \alpha}}
\newcommand{\dist}{\mbox{dist}}
\newcommand{\mod}{\mbox{mod }}
\newenvironment{double}{\renewcommand{\baselinestretch}{2}\protect\Large
\protect\normalsize}{}

\maketitle
\begin{abstract}
We prove that non-hyperbolic non-renormalizable quadratic polynomials
are expansion inducing. For renormalizable polynomials a counterpart of this statement is 
that in the case of unbounded combinatorics renormalized mappings
become almost quadratic.   
Technically, this follows from the decay of the box geometry. 
Specific estimates of the rate of this decay are shown which are sharp
in a class of S-unimodal mappings combinatorially related to rotations
of bounded type. We use real methods based on cross-ratios and Schwarzian
derivative complemented by complex-analytic estimates in terms of
conformal moduli.  
\end{abstract}

\section{Introduction}
\subsection{Overview}
In recent years, a rather dramatic progress occurred in the study of
real quadratic polynomials, or, more broadly, S-unimodal mappings with
a singularity of quadratic type. Examples of progress include better
understanding of measure-theoretical attractors and quasisymmetric
classification of quadratic polynomials which lead to a proof of the
monotonicity conjecture in the real quadratic family. 

This progress was partly based on ingenuous new estimates. 
The main breakthrough, however, was in achieving a better
understanding of the rich dynamics of unimodal mappings in conjunction
with their  geometry, and thus being able to apply appropriate tools
in different cases. Most of the progress in this direction seems to be
due to the application of the idea of {\em inducing}. The first
striking application of inducing to the study of unimodal maps was in
the work~\cite{invmes}. In that work useful geometrical and analytic 
estimates were obtained only in judiciously chosen special cases.
Another notable step was the work of~\cite{gujo}. An attempt was made
there to handle all cases, though some patterns emerged as
analytically unmanageable. Then, independently, two approaches appeared. 
One of them was the inducing construction of~\cite{yours}, which is
also the underlying approach of the present paper. Here a clear
and complete topological model of unimodal dynamics was obtained, with
satisfactory estimates in most cases, except in what was called an {\em
infinite box case}. The infinite box case was subsequently solved for
mappings with a quadratic singularity based on the phenomenon of {\em
decaying box geometry}. The phenomenon was first noticed in~\cite{yours}, 
however proved there only in some cases. An estimate called the {\em starting
condition} was provided which, if satisfied,
allowed one to prove the decaying geometry in general for S-unimodal
maps. It should be noted that about the same time a similar case of
decaying
geometry was observed independently by~\cite{tanver} and~\cite{gracz}
for circle mappings with a flat piece. It is not known whether there
is more than an analogy between both cases.   
 
A breakthrough work was~\cite{yoc}, see~\cite{miln} for a
description. This was done for non-renormalizable quadratic
polynomials, complex as well as real. 
Inducing was not directly mentioned, but implicitly
present in the construction of a Markov partition. It is believed that
the approaches of~\cite{yours} and~\cite{yoc} give equivalent
sequences of partitions for real polynomials. A great achievement
of~\cite{yoc} was in being able to get estimates in the ``hard case'' which
again emerged (and was called {\em persistently recurrent}). This was
done by watching how pieces of the partition nest in one another with
certain moduli, and using a computation quite similar to one done
in~\cite{brahu} (curiously, this last work was about cubic polynomials.) 

In the real situation, analytic methods were refined in the study of
one particular topological class of the so-called {\em Fibonacci
polynomial}. This map was proposed in~\cite{hokel} as an interesting
example to study. From the point of view of ~\cite{yours}, the Fibonacci polynomial shows the simplest example of the infinite box case. 
In the work of~\cite{limo} a complex-analytic idea was
applied to obtain estimates in the Fibonacci case. Another approach was used
in~\cite{nowike} where the same results as in~\cite{limo} were proved  
by purely real methods based on negative Schwarzian. Together with 
arguments of~\cite{yours}  or~\cite{ma} based on inducing, the non-existence
of non-Feigenbaum type Cantor attractors (proved in~\cite{opo}) implies
{\em induced expansion}. 

Finally, there was a work of~\cite{hyde} which addressed the geometry
of renormalizable quadratic polynomials. Again, a hard case emerged
when the trajectory of the restrictive interval was allowed to be
arbitrarily long (unbounded case), and in terms of the construction
of~\cite{yours} needed to be tracked through a long sequence of box
returns. In~\cite{hyde}, the idea borrowed from~\cite{limo} was used
which consists of
introducing an artificial map for which the starting condition holds.
Then, a conjugacy with good quasiconformal properties is constructed
between this mapping and the given quadratic polynomial which forces
the starting condition for the polynomial.

\subsection{Statement of results}
\paragraph{Introduction}
In this paper we rely on the inducing approach of~\cite{yours}, and
our goal is to verify the starting condition in the box case. Our
method is a mix of real and complex arguments and works just as well
for renormalizable polynomials of unbounded type as for
non-renormalizable ones. The method  is direct, that is independent of
quasiconformal conjugacies with artificial maps, the
tableau computation of~\cite{brahu} or moduli estimates of ~\cite{yoc}. On the technical level, our main theorems 
follow from Theorems A and B stated further in the text. The derivation is  by
arguments of~\cite{hyde} (fully presented here.) 

\paragraph{Main theorems.}
\subparagraph{Definition of our class of mappings.} 
\begin{defi}\label{defi:9ma,1}
We define class $\cal F_{\eta}$ to comprise all unimodal mappings of the
interval $[0,1]$ into itself normalized so that $0$ is a fixed point
which satisfy these conditions:
\begin{itemize}
\item
Any $f\in {\cal F}$ can be written as $h((x-\frac{1}{2})^{2})$ where
$h$ is a polynomial defined on a set containing $[0,1/4]$ with range 
$(-\eta, 1+\eta)$. 
\item
The map $h$ has no critical values except on the real line.
\item
The Schwarzian derivative of $h$ is non-positive.
\item
The mapping $f$ has no attracting or indifferent periodic cycles.
\end{itemize}
We also define 
\[ {\cal F} := \bigcup_{\eta>0} {\cal F}_{\eta}\; .\]
\end{defi}

For any mapping $f\in {\cal F}$ we define the {\em fundamental
inducing domain} as follows. From the non-existence of attracting or
indifferent periodic points it follows that there is a repelling fixed
point $q>1/2$. The fundamental inducing domain is the interval
$(1-q,q)$. Almost every orbit passes infinitely many times through the
fundamental inducing domain.    

\subparagraph{Theorem about non-renormalizable mappings.}
\begin{theo}\label{t:1}
Let $f\in {\cal F}_{\eta}$ be non-renormalizable. Then on an open, dense and
having full measure subset of the fundamental inducing domain one can define a
continuous function $t(x)$ with values in positive integers so that
$f^{t(x)}$ is an expanding  Markov mapping. That is, restricted to a
maximal interval on which $t(x)$ is defined an constant, $f^{t(x)}$ is
a diffeomorphism onto $(1-q,q)$, expanding, and with distortion
(measured as the variation of the logarithm of the jacobian) bounded
by depending on $\eta$ only. 
\end{theo}

Theorem~\ref{t:1} has a number of consequences (see~\cite{yours}). It
gives an alternative proof of the non-existence of ``exotic''
attractors in class $\cal F$ (already known from~\cite{opo}). 
It also gives an approach to constructing invariant measures. 

\subparagraph{Theorem in the renormalizable case.}
\begin{defi}\label{defi:14ma,1}
Let $f\in {\cal F}$. A point $x$ in the domain of $f$ is called {\em
almost parabolic with period $m$ and depth $k$} provided that:
\begin{itemize}
\item
the derivative of $f^{m}$ at $x$ is one,
\item
$f^{m}$ is monotone between $x$ and the critical point, 
\item
$k$ consecutive images $f^{m}(1/2), \ldots, f^{km}(1/2)$ are between
$x$ and $1/2$. 
\end{itemize} 
\end{defi}

\begin{theo}\label{th:9ma,2}
Let $f\in {\cal F}_{\eta}$ be renormalizable, and let $n$ be the return time
of the maximal restrictive interval into itself. Denote by $k(n)$ the
maximum of depths of almost parabolic points with periods less than $n$.
Specify a number $D>0$. For every given $k$, a number 
$N(\eta,D,k)$ exists independent of $f$ so that if $n>N(\eta,D,k)$ and
$k(n)\leq k$, then $f^{n}$ on
a neighborhood of $1/2$ is affinely conjugate to a mapping from ${\cal
F}_{D}$.   
\end{theo}

Theorem~\ref{th:9ma,2} ``almost complements'' the
theory of renormalizable mappings developed in~\cite{miszczu}. In
fact, it says that such a theory at least in some aspects is much
simpler for renormalizable mappings of unbounded type. The exclusion
of the unbounded case with almost parabolic returns is the only gap.
Theorem~\ref{th:9ma,2} is a critical step in the proof of monotonicity
in the real quadratic family, see~\cite{hyde} (where, by the way, the
theorem is stated wrongly without excluding the almost parabolic case.)   
         
\subparagraph{Technical theorems.}
The strongest results of our paper are contained in technical theorems
A,B, C, D. They imply theorems 1 and 2. In addition, Theorem C gives
exact bounds on the exponential rate of decay of box geometry for
S-unimodal rotation-like maps. Theorem D concerns the decay of {\em
complex} box geometry and can serve as an important step in proving
the density of hyperbolicity in the real quadratic family,
see~\cite{hyde}. The technical theorems are stated for objects that we
call {\em box mappings} and their statement must be postponed until
those are defined.

\paragraph{Plan of the work.}
All three theorems follow from Theorem B which we will formulate later
which essentially says that the starting condition holds for any map
from $\cal F$ after a bounded number of box inducing steps. To prove
this we will use two complementary methods. The first one is based purely
on real-variable considerations. The advantage of a real-variable
approach is in robust estimates. The disadvantage is sensitivity to the
combinatorial complexity of the problem. Thus, we conduct our
estimates only in the case we call ``rotation-like'' in which
estimates are subtle and the combinatorics not harder than for circle
rotations. In this way we prove Theorem A for real mappings with
rotation-like behavior.  
Our second method is complex-analytic and based on watching
annuli which are mapped by the dynamics and nest inside one another.     
This method will allow us to cover the full realm of combinatorial 
possibilities, but estimates are weaker. In particular, they would be
too weak in the rotation-like case. More precisely, to prove the
starting condition we will need to show that certain ratios are small
after a number of box steps. Depending on the dynamics, by each box
step the geometrical ratios decrease either by a multiplicative
constant uniformly less than $1$, or by being raised to a power
uniformly greater than $1$. The complex approach is capable of
accounting for the second type of phenomenon, but, at least at
present, misses the first one. It turns out that it is exactly the
rotation-like case that exhibits the first, slower, type of decay.    

After proving the technical theorems, we will derive Theorems 1 and 2. 
  
\paragraph{Acknowledgements.} Jacek Graczyk gratefully acknowledges the
hospitality of the Institute for Mathematical Sciences in Stony Brook
where, amidst beautiful Long Island spring,  most of this work was
done.  Both authors thank M. Jakobson for dicussions regarding 
the rotation-like case. We are also grateful to S. Sutherland for
sharing his very enlightening computer-generated pictures, and to J.
Milnor for remarks concerning references and the historical background
of our problem.    

\subsection{Box mappings.}
\paragraph{Real box mappings.}
The method of inducing was applied to the study of unimodal maps first
in~\cite{invmes}, then in~\cite{gujo}.  
In~\cite{kus} and~\cite{yours} an elaborate approach was developed to
study induced maps, that is, transformations defined to be iterations
of the original unimodal map restricted to pieces of the domain. We
define a more general and abstract notion in this work, namely:

\begin{defi}\label{defi:704a,1}
Consider a transformation $\phi$ defined on an open dense subset of
an interval $I$ into $I$. Call restrictions of $\phi$ to connected
components of the domain {\em branches}. If each branch is at least
three times differentiable and the Schwarzian derivative is
non-positive wherever defined, we will call $\phi$ a {\em piecewise map.}
\end{defi}

We will deal with two types of branches, namely monotone and folding.
A monotone branch is a diffeomorphism onto its image, while a folding
branch arises as a quadratic polynomial pre- and post-composed with 
diffeomorphisms. Examples include induced maps studied in works cited above.

Next, we want to impose more specific conditions on the images of
the branches of a generalized induced map. 
\begin{defi}\label{defi:1,1}
A box map on $I$ is a piecewise map on $I$ more combinatorial
structure. Namely, there
is a finite nesting sequence of intervals, called boxes or real boxes,
\[ I = B_{0} \supset \cdots \supset B_{n} \]
indexed by a growing sequence of integers and halves of odd integers. 
All boxes are assumed symmetric with respect to the critical point. 
It is further assumed that there is only one
central folding branch whose domain is $B_{n}$ and whose image is
contained in some $B_{k}$ with
the boundary of $B_{n}$ going into the boundary of $B_{k}$. It is
then said that the {\em rank} of the central branch is $k$. 
All other branches are
monotone and each maps onto some $B_{l}$ which determines its
rank.  

Furthermore, no branch of positive rank can share a
common endpoint with another branch of positive rank or folding, and
if the domain of some branch has a non-empty intersection with a box,
then it is contained in this box. 
\end{defi}

This is an extension of the notion of induced map used in~\cite{kus}.
Namely,
\begin{defi}\label{defi:1,2}
A box map is called full if $B_{1}$ is the domain of the central
branch which is of rank $0$, and the highest rank is $1$.
\end{defi}

\paragraph{Complex box mappings.}
We will now give a precise meaning to complex extension of real box
mappings. Given a real box $B$, a corresponding {\em complex box} must be
an open disk symmetrical with respect to the real axis which intersects the
real line exactly along $B$. Given a real box structure, a
compatible {\em complex box structure} will be a sequence of nesting
complex boxes corresponding to the real boxes. 

\subparagraph{Hole structures.}
Now suppose that a real box mapping $\phi$ is given whose all branches
are real-analytic and pick a compatible complex box structure. 
A {\em hole structure} means that for the domain of every branch of
$\phi$ save monotone branches of rank $0$ an open disk (hole) is
chosen symmetrical with respect to the real line and intersecting the
line along the domain. It is further assumed that the branch has an
analytic continuation to its hole of the same topological type as the
real branch, i.e. either univalent for monotone branches or degree $2$
for folding branches. We will not hesitate to talk of monotone or
folding branches for complex box mappings.
This analytic extension will be called the {\em
complex branch}. Lastly, the image of the complex branch is exactly
the complex box corresponding to the rank of the real branch.        

\begin{defi}\label{defi:5mp,1}
Given a real-analytic real box mapping $\phi$, the choice of a complex
box structure and a hole structure, defines a {\em complex box
mapping.} Thus, by a complex box
mapping we will formally understand the conglomerate of all three
compatible structures: a real box map, complex box structure, and hole
structure. With a slight abuse of language, we will also extend this
name to the union of complex branches of the complex box map. 
\end{defi}

\paragraph{Special type box mappings.}
\begin{defi}\label{defi:9mp,1}
A type I box mapping of rank $n$ is a first of all a box
mapping whose central domain is $B_{n}$.  The central branch has rank
$n-1 \leq n' < n$ corresponding to the element of the box structure
directly preceding $B_{n}$.   All monotone branches of
positive rank have rank $n$. This definition makes equal sense for
real and complex box mappings. 
\end{defi} 
 Note
that the full map is type I of rank $1$, with $n' = 0$. 
Most of our work will be done in terms of type I box maps.  

\begin{defi}\label{defi:9mp,2}
A type II box mapping of rank $n$ is a box mapping whose central
domain is $B_{n}$ and all branches of positive rank map onto $B_{n'}$
which the element of the box structure directly preceding $B_{n}$.
\end{defi}

\paragraph{Separation symbols for complex box mappings.}
\subparagraph{Definition of the symbols.}
Now, let $\varphi$ be any type I complex box mapping. Fix the notations so
that the central hole of $\varphi$ is labeled $B_{n}$ with $n>1$. 
Pick a monotone branch $B$ of rank $n$. Also, pick another hole $C$
contained in $B_{n-1}$. We may assign to $B$ its 
{\em separation symbol} from $C$ which is simply an ordered quadruple
of real non-negative numbers:
 \[ s(B) := (s_{1}(B),\cdots,s_{4}(B))\; .\]
A valid separation symbol by definition implies the existence of
certain annuli with moduli estimated from below in terms of the
components of the symbol.  
We first assume that there are annuli $A_{1}$ and $A_{2}$,
both selected for the given $B$ even though the dependence on $B$ is
not emphasized in the notation. Both annuli are contained in
$B_{n'}$. The annulus $A_{2}$ surrounds the hole
$C$ separating it from the boundary of $B_{n-1}$ as well as the
hole $B$. In addition, $A_{2}$ is not allowed to intersect any holes
which meet the real on the side of $B$ opposite to $C$. This, in
effect means that $A_{2}$ also separates $C$ from all holes strung on
the real line ``behind'' $B$.
Then $A_{1}$ separates $A_{2}$ from the boundary of
$B_{n-1}$. The number $s_{2}(B)$ is a lower bound on the modulus of $A_{2}$
and $s_{1}(B)$ is a lower bound on the sum of moduli of $A_{1}$ and
$A_{2}$. 

We then proceed to select annuli around $B$ which will give the
meaning of the two remaining components of the symbol. First, the
annulus $A'$ is chosen to separate $B$ from the domain of the extension of
the branch defined on $B$. This extension is 
of rank $n'$. Then, the
existence of $A_{3}$ is postulated which surrounds $A'$ separating it
from $C$ and the boundary of $B_{n-1}$. Also, $A_{3}$ is not allowed
to intersect any holes that meet the real line on the side of $C$
opposite to $B$. Finally, $A_{4}$ separates
$A_{3}$ from the boundary of
$B_{n-1}$. The component $s_{3}$ is then assumed to be a lower bound
for the sum of moduli $A'$ and $A_{3}$, while $s_{4}$ is a lower bound
of the sum of all three moduli: $A'$, $A_{3}$ and $A_{4}$. 

In this paper we will always assume $C$
to be the central hole. We will then call the separation symbol of $B$
from $C$ the {\em critical symbol} of $B$. 

This fully lists the geometrical properties implied by a valid
separation symbol. 

\paragraph{Normalized symbols.}
We will now arbitrarily impose certain algebraic relations
among various components of a separation symbol.
Choose a number $\beta$, and $\alpha:=\beta/2$, together
with $\lambda_{1}$ and $\lambda_{2}$. Assume
$\alpha\geq\lambda_{1},\lambda_{2}$, $\lambda_{1}\geq -\frac{\alpha}{2}$,
$\lambda_{2}\geq -\frac{\alpha}{2}$ and $\lambda_{1} + \lambda_{2}\geq 0$.   
If these quantities are connected with a separation symbol $s(B)$ as
follows
\[ s_{1}(B) = \alpha + \lambda_{1})\; ,\]
\[ s_{2}(B) = \alpha - \lambda_{2}\; ,\]
\[ s_{3}(B) = \beta - \lambda_{1}\; ,\]
\[ s_{4}(B) = \beta + \lambda_{2}\; .\]

we will say that $s(B)$ is normalized with norm $\beta$ and
corrections $\lambda_{1}$ and $\lambda_{2}$.
\subparagraph{Separation norm of a box mapping.}
For a type I complex box mapping $\phi$ its {\em separation norm} is
defined as the supremum of values of $\beta$ for which valid normalized
critical symbols with norm $\beta$ exist for all univalent branches.   

\subsection{Box inducing process}
\paragraph{Description of the process.}
We will now describe a {\em standard inducing step} which is a
procedure that takes a box mapping called $\phi$ and returns a type
I box mapping called $\Phi$. $\Phi$ is induced in terms of $\phi$, i.e.
all branches of $\Phi$ are compositions of branches of $\phi$. 

\subparagraph{Close and non-close returns.}
If the escape time of the critical point is $1$, we say that the
mapping shows a non-close return. Otherwise, it is a close return. 
The box case occurs when the image of the critical point is
found in the domain of a monotone branch of positive rank. The
remaining two cases when the push-forward image is in an rank $0$
domain, or beyond the domain of the induced map (the Misiurewicz
case), are well understood by the methods of~\cite{yours}, so we
always assume that a box case occurs. The procedure we describe applies
equally well to complex and real box mappings. So suppose that a box
mapping $\phi$. Let the central domain be $B_{n}$.

\subparagraph{The first filling.}
Construct $\phi'$ by replacing $\phi$ on the central domain (hole) with the
identity. This gives a valid box mapping. 
Next, define $\phi_{1}$ by replacing the central branch $\psi$ of
$\phi$ with $\phi'\circ\psi$. The new central domain is the
preimage by $\psi$ of the domain containing the critical value. The
new central domain is adjoined to the hole structure as $B_{n+1}$.  

\paragraph{Filling-in.} In order to obtain a type I mapping, we "fill in" all monotone branches obtained by the first filling.  This
process was introduced in~\cite{yours}. That is, each point in the domain
of a monotone branch is mapped by monotone of branches of positive
rank until it leaves
the union of their domains. Its first image outside of those domains,
by definition, is the image of the point under a new map. This
definition works except for a Cantor set of points.  The mapping
obtained in this way has all monotone (univalent) branches of rank
$n+1$.  

This construction can be carried out regardless of whether the return
was close or not. In the case of a non-close return this
completes a {\em standard inducing step} and the map so obtained is
$\Phi$. For close returns, it is
convenient to bundle together several steps. 

\paragraph{Close returns.}
In the case of a close return, we proceed as follows. 
We construct $\phi'$ in the usual way be replacing the central branch
with the identity, and then define $\phi_{1}$ by substituting the
central branch $\psi$ with $\phi'\circ\psi$. In other words, the first
filling occurs like for a non-close return. The central domain is
$B_{n+1}$, but we will more often use the notation $B^{1}$. 
We then look whether the
critical value of $\phi_{1}$ (which is still the same as the critical
value of $\psi$) is in $B^{1}$. If so, we repeat the same step with
$\phi_{1}$ instead of $\phi$ to obtain $\phi_{2}$ and continue in this
way until the critical value is no longer in the central domain
$B^{k}$. This must happen for some $k$, or the intersection of all
central domains would be a restrictive interval and $\phi$ would be
suitable. The mapping has a bunch of
monotone branches all of which have rank $n$ and are compositions of
the central branch applied a number of times with branches of $\phi$.
We then make it a type I mapping $\phi$ by filling-in of those
monotone branches so that they become of rank $n+k$, i.e. all map onto
$B_{k}$. The central branch remains unchanged of rank $n+k-1$. Call
this new type I mapping $\tilde{\phi}$.   
Then the last step is to
execute the standard inducing step on $\tilde{\phi}$ which, by
definition, will give $\Phi$. Finally, we rearrange the box structure
of $\Phi$ defining $B_{n+1}$ as the central domain of $\Phi$ and
$B_{n+1/2}$ as $B^{k}$.   

The reader may check that this is the same map that we would get
simply following the standard inducing step described for non-close
returns $k+1$ times. However, our description gives a more direct
insight into the origin of branches of $\Phi$.      

\paragraph{Immediate preimages.}
Regardless of how $\Phi$ was constructed, its {\em immediate} or {\em
primary branches} are those univalent branches which are restrictions of the
central branch $\psi$. Note that there may be none or two
immediate preimages, depending on whether the real range of the
central branch covers the critical point. In the case of a close return the immediate preimages are formed at the last stage of inducing on $\tilde{\phi}$. The well-known "Fibonacci case" is an example when all monotone branches occurring in the constru

\paragraph{Filling-in and hierarchy of branches.}
Consider a abstract setting in which one has a bunch of univalent branches with common range $B'$ and fills them in to get branches mapping onto some $B\subset B'$. 
The original branches mapping onto $B'$ will be called {\em parent branches} of the filling-in process. Clearly, every branch after the filling-in has a dynamical extension with range $B'$. For two branches, the domains of these respective extensions m
\begin{itemize} 
\item
one branch is mapped onto the central branch, or
\item
the domains of their dynamical extensions are disjoint.
\end{itemize}
In the second case we say that the original branches were {\em independent}. 
Otherwise, the one mapped onto a monotone domain is called {\em subordinate} to
the other one which was transported onto the central domain.

We then distinguish the set of "maximal" branches subordinate to none. They 
are mapped by their dynamical extensions directly onto the
central domain. Therefore, the domains of extensions of maximal branches are
disjoint. They also cover domains of all branches. The extensions of maximal
branches exactly the {\em parent branches} of the filling-in process.

For example, in the non-close return the first filling gives a set of parent branches, two of which may be immediate, which later get filled in. In the close return filling-in is done twice, so we will be more careful in speaking about parent branches i
\section{Statement of main results}
\paragraph{Rotation-like returns.}
Let $\phi$ be a type I map obtained in a standard inducing step. Then
its immediate branches are defined as in the previous paragraph. 

\begin{defi}\label{defi:7mp,1}
We say that $\phi$ exhibits a rotation-like return in the following
situation. When the return is not close, it is rotation-like if and
only if the critical value lands in an immediate branch of $\phi$.
When the return is close map the critical value by the central branch
until the first exit from the central domain. The return is
rotation-like if upon the first exit from the central domain the
critical value gets into an immediate branch of $\phi$ and is mapped
by this branch into $B^{k}$. 
\end{defi}

The above definition is slightly technical however its advantage
is that the existence of rotation-like
maps follows immediately from the definition.
\begin{prop}\label{p:1}
Rotation like maps do exist.
\end{prop}

\begin{proof}
Take a one parameter full family of S-unimodal maps. For each of them
do inducing procedure. Then the position of the critical 
value of the central domain during each step of induction 
depends continuously on the parameter value. Hence in particular
we obtain all rotation-like maps.
\end{proof}

Now we are in position to give a simple characterization of
rotation-like maps which also will justify the name of this
class of unimodal maps. 

\subparagraph{Rotation-like sequences.}
We will say that a sequence of type I box mappings is rotation-like if
the each arises from the previous one by a standard inducing step, and
each return, with the exception of the first one for which immediate
preimages are not defined, is rotation-like.

\begin{fact}\label{rot}
For a rotation-like sequence, there is an inductive formula relating 
consecutive
central branches of the inducing procedure
\[ f_{n+1}=f_{n-1}\circ f_{n}^{a_{n}},\]
where $f_{j}$ denotes the central branch of $\phi_{j}$, and $a_{n}$ is
a the smallest $i$ such that the $i$-th iterate of
the central branch of the $n$-th mapping maps the critical point outside of the
central domain. E.g., $a_{n}=1$ is equivalent to saying that the
$n$-th map shows a non-close return. 
We define $a_{0}$ by a requirement that $f_{2} = f^{a_{0}}\circ f_{1}^{a_{1}}$.
\end{fact}

Following the analogy with the circle homeomorphisms
we will introduce a concept of a rotation number
for our class of maps. 

\begin{defi}
The rotation number $\rho(f)$ of $\phi$ is equal to  
\[\rho_{n}=\frac{1}{a_{1} +\frac{1}{a_{2}+\cdots}}\]
which can be written shortly as $[a_{0},a_{1},\cdots  ]$ 
using the formalism of continued fractions.
\end{defi}

\subsection{Technical theorems}
\paragraph{The starting condition.}
\begin{defi}\label{defi:10ma,1}
We say that a type I or type II box mapping of rank $n$ satisfies the {\em
starting condition} with norm $\delta$ provided that  
Let $|B_{n}|/|B_{n'}| < \delta$ and if $D$ is a
monotone domain of $\phi$, then also $|D|/\dist(D,\partial B_{n'})< \delta$. 
\end{defi} 

\begin{fact}\label{fa:8mp,1}
Let $\phi$ be a real box mapping of type I of rank $n$.  
Pick $\tau > 1$ and assume that  the central branch is $\tau$-extendable.
Then, there is a number $\delta(\tau)$ bounded away from $0$ for
$\tau$ in any closed subset of $(1,\infty)$ with the following
property.     
Suppose that $\phi_{i}$
with $\phi_{0} = \phi$ be a sequence of type I box mappings such that
$\phi_{j+1}$ arises from $\phi_{j}$ by a standard inducing step and
let $\phi$ satisfy the starting condition with norm $\delta(\tau)$. 
Then, in the box construction of~\cite{yours}, the ratios
$\delta_{n+i}/\delta_{n+i-1}$ tends to $0$ at least exponentially fast
with $i$ with rate given by an absolute constant.  
\end{fact}
\begin{proof}
This Fact was proved in~\cite{yours}. 
\end{proof}

The box construction of~\cite{yours} is slightly different from the
inducing construction we use. We do not need the details now. 

\paragraph{Theorem A about real box mappings.}
We state the theorem as follows.\newline
{\bf Theorem A}

\begin{em}
Let $\phi=\phi_{0}$ be a type I real box mapping. Let $B_{n}$ be its
central domain. Suppose that the ratio of lengths $|B_{n}|/|B_{n'}|$
is $1-\epsilon$. Next. let $\phi_{i}$ be a rotation-like sequence
derived from $\phi$. Specify a $delta>0$. Then, there is a function 
$K(\epsilon,\delta)$ of
$\epsilon$ and $\delta$ only and independent of $\phi$,
bounded on any compact compact subinterval of $(0,1]^{2}$, with the
property 
that
for $i\geq K(\epsilon,\delta)$ the mapping $\phi_{i}$ satisfies the starting
condition with norm $\delta$.   
\end{em}

The proof of Theorem A uses purely real methods. It generalizes the
result of~\cite{nowike}. A natural question is whether an analogous
result can be demonstrated by real methods for an arbitrary box
sequence of induced maps. In principle, that should be possible, but
technical difficulties are daunting. 

\paragraph{Theorem B about complex box mappings.}
The technical complexity of the general box case becomes tractable
when one works with complex box mappings. Hence the theorem:\newline
{\bf Theorem B}

\begin{em}
Let $\phi=\phi_{0}$ be a type I complex box mapping, and $\phi_{i}$ be
a sequence of complex box mappings such that $\phi_{i+1}$ is derived
from $\phi_{i}$ in a standard inducing step. Suppose that the
separation norm of $\phi$ is $\beta$. Also, specify a $\delta>0$.
Then, there is a function
$K(\beta, \delta)$ depending solely on $\beta$ and $\delta$ 
and bounded on any closed set
contained in $(0,\infty)^{2}$ with the property that if $i\geq
K(\beta, delta)$,
then $\phi_{i}$ as a real box mapping satisfies the starting
condition with norm $\delta$.
\end{em}

\paragraph{Sharp estimates for rotation-like mappings.}
\noindent {\bf Theorem C}\begin{em}
For any S-unimodal rotation-like map 
there exist positive constants $K_{1}$, $K_{2}$ 
 and $\kappa_{1}, \kappa_{2} <1$ depending only on the initial
geometry, i.e. the number $\epsilon$ of Theorem A, so that 
\[ K_{1}\kappa_{1}^{a_{1}\cdot ...\cdot a_{n-1}} \leq \tau_{n} \leq 
K_{2}\kappa_{2}^{a_{1}\cdot ... \cdot a_{n-1}}.\]
\end{em}

This is an improvement of Theorem A which also gives the lower bound
on the rate of decay of box geometry. The proof is by purely real
methods. 

\paragraph{Growth of conformal moduli.}
\noindent{\bf Theorem D}\begin{em}
Let $\phi=\phi_{0}$ be a type I complex box mapping of rank $n$, and $\phi_{i}$ form a sequence of complex box mappings derived from $\phi$ by the box inducing process.  Suppose that box ratios on the real line decrease at least exponentially fast, i.e.
\[ \frac{|B_{n+j}|}{|B_{(n+j)'}|} < C^{j}\]
with $C<1$. Let $\beta_{j}$ denote the separation norm of $\phi_{j}$.
Then, there is a number
$\overline{C}(C,\beta_{0})$ so that if the separation norm of
$\phi_{0}$ is at least $\beta_{0}$, then 
$\beta_{j} \geq \overline{C}\cdot j$. The constant $\overline{C}$ only
depends on its specified parameters. 
\end{em}

Theorem D claims a decay of the conformal geometry in a sequence of complex box mappings derived by inducing. This phenomenon seems to be the basis of many recently 
obtained results, see \cite{hyde}.  In view of Theorem B and Fact~\ref{fa:8mp,1} the assumption of Theorem D regarding the decrease of ratios is automatically satisfied for any complex box mapping with constant $C$ depending only the initial separation 
\section{Real induction}
In this section we prove Theorem A. Suppose that in the situation of
Theorem A a rotation-like sequence $\phi_{i}$ is given,
$i=0,1,\cdots$. 
\subsection{Rotation-like sequences}
\paragraph{Notations.}
For the real induction, we fix the critical point at $0$. The central
branch of $\phi_{i}$ will be called $f_{i}$.   
Denote the endpoints of  the box $B_{n}$ by
$z_{n}^{-}$ and $z_{n}^{+}$. We adopt the following convention
of ascribing signs to points:  $`+'$ written as a superscript 
indicates the endpoint of $B_{n}$ which lies closer to the 
critical value of $f_{[n+1]}$ (where $[\cdot]$ denotes the integral
part.) In other words $f_{[n+1]}(0) \in (0,z_{n}^{+})$.
Each central branch can be represented as a composition of 
a diffeomorphism $h_{n}$ and a quadratic map $g$. We know
that the image of $B_{n}$ is contained in $B_{n-1}$.
The following Lemma describes extendability of diffeomorphisms
$h_{n}$.

Generally, we will use $(x,y)$ to denote the interval from $x$ to $y$,
regardless of the ordering of $x$ and $y$.
\begin{lem}\label{ext}
If $n\geq 3$, the diffeomorphism  $h_{n}$ extends on some neigh\-bor\-hood
of $g(B_{n})$ so that the image of the extension coincides with
 $(f_{n-2}(0),z_{(n-2)'}^{-})$.
\end{lem}

\begin{proof} 
By construction, $h_{n}$ extends to a neighborhood
of $g(B_{n})$ which is mapped on $B_{n-1}$.
 Take arbitrary $n$. We know that $f_{n}= f_{n-2}\circ f_{n-1}^{a_{n-1}}$.
Consider two cases.

\begin{itemize}
\item If $a_{n}=1$ then the range of a monotone extensions of $h_{n}$
is the same as that of the composition $f_{n-2}\circ h_{n-1}$
 which is clearly $(f_{n-2}(0),z_{(n-2)'}^{-})$.
and next pull it back by the extension of $f_{n}$ whose
image by the inductive hypothesis covers $B_{n-2}$.
 The resulting interval is the
\item For $a_{n}>1$  choose an extension of $h_{n-1}$ so that
 $f_{n-1}^{a_{n}-1}\circ h_{n-1}$ maps it diffeomorphically 
onto $(0,z_{n-2}^{+})$. The image of monotone branches of 
$f_{n-2}$ gives the desired range of the extension.
\end{itemize}

This proves Lemma~\ref{ext}.
\end{proof}

In the ``real part'' of this paper we will extend diffeomorphisms
$h_{n}$ every time only over a one side of their domains. We
distinguish between two directions
of one-sided extendability of $h_{n}$. The key observation is that the points
$\;\{f_{n-2}(0),f_{n}(z_{n}^{+}),f_{n}(0),z_{n-3}^{-}\}\;$
are always arranged according either to the natural or reversed order
of the real line. Observe, that $f_{n}$ extends further ``through the
head'' meaning in the direction of the critical value, where $h_{n}$
can be extended up to $z^{-}_{(n-2)'}$, than ``through the legs''
meaning in the direction of $f_{n}(\partial B_{n})$ where the
extension is only up to $f_{n-2}(0)$ which is closer to $B_{n'}$ than
$z^{+}_{(n-2)'}$.  
The extendability of central branches plays a crucial
role in estimates of the distortion. By real K\"{o}be's Lemma the
distortion depends
only on the relative scale of the images of domains
with respect to the images of their extensions. As it happens (and
will be proved), these scales
will improve during inducing procedure finally forcing {\em the starting
condition}. We will need however some initial extension to start with.

\subparagraph{Estimates a priori.} Denote by $\tau_{n}$, $n$ integer, the
ratio $|B_{n}|/|B_{n'}|$. Let $\tau$ be the supremum of $\tau_{n}$
with respect to $n$. 
We have the following Lemma:
\begin{lem}\label{lem:8mp,1}
Under the assumptions of Theorem A, there is a function
$K_{3}(\epsilon)$ bounded on any closed subset of $(0,1]$ with the
property that $\tau_{n} < \tau=0.37$ for $n\geq K_{3}(\epsilon)$.
\end{lem}

The proof of this lemma is technical and apart from main ideas of this
paper, so we put it in a separate section at the end of the real estimates.

\paragraph{Ratios and cross-ratios.}
Suppose we have three points $a,b,c$ arranged so that
$a\notin [b,c]$. 
Let us define a few relative scales  of the interval
$(b,c)$ with respect to $(a,c)$. 
\begin{defi}
The exclusive ratio of the interval $(b,c)$ with respect to $(a,c)$
is given by 

\[\R_{e}(b,c;a)=\frac{|b-c|}{\dist((b,c),a)},\]
whereas their inclusive ratio by
\[\R_{i}(b,c;a)=\frac{|b-c|}{\max(|b-a|,|c-a|)}.\]
Set $\R(b,c;a)$ to be equal to the geometric mean
of the inclusive and exclusive ratios.
\[\R(b,c;a) = \sqrt{\R_{i}(b,c;a)\R_{e}(b,c;a)}\]
\end{defi}
Together with ratios we will often use 
cross-ratios. The types of cross-ratios we use
are expanded by homeomorphisms with negative Schwarzian derivative.

\begin{defi}
Suppose we have a quadruple $a,b,c,d$ ordered so that $a<b<c<d$ or
reversely.  
Define  their inclusive cross-ratio as
\[\Cr_{i}(a,b,c,d)= \frac{|b-c||d-a|}{|c-a||d-b|},\]
and their exclusive cross-ratio as
\[\Cr_{e}(a,b,c,d)= \frac{|b-c||d-a|}{|b-a||d-c|}.\]
Finally set 
\[\Cr(a,b,c,d)= \sqrt{\Cr_{i}(a,b,c,d)\Cr_{e}(a,b,c,d)}\]
\end{defi}

\paragraph{Distortion.}
Suppose that we have an expression $A$ which is defined in terms
of distances between points (like ratios and cross-ratios.)
Then we consider $f_{*}(A)$ obtained by replacing given points with
their images by $f$. 
We will measure the distortion of this transformation
by the ratio $f_{*}(A)/A$. For example, if we set $A=\Cr(a,b,c,d)$
then the distortion by $f$ is equal to
\[\frac{\Cr(f(a),f(b),f(c),f(d))}{\Cr(f(a),f(b),f(c),f(d))}.\]

\subsection{Induction parameters}
In this subsection we will introduce quantities which
will describe geometry of partitions given by our inducing
procedure. Next we will compose a quasi-invariant which
after a finite number of inducing steps  will decrease
at least exponentially fast. The real induction parameters formulated
here will directly  correspond to these in the complex
part. The same concerns induction formulae. This 
suggests that estimates from the complex part of our work can
somehow be translated into
the corresponding ones in the real line. This would enable one to give a proof
by purely real methods. However, the combinatorial
complexity of such an approach seems formidable.

Denote by  $(x_{n}^{-},x_{n}^{+})$ the domain of the 
{\em primary}
branch of rank $n$  which contains the critical value of
$f_{n}^{a_{n}}$. {\bf Set} $v_{n}=f_{n}^{a_{n}}(0)$.

\begin{defi}
The center of $(x_{n}^{-},x_{n}^{+})$, denoted by $c_{n}$, is 
defined by the condition $\phi_{n}(c_{n})=0$.
\end{defi}
Following the convention of ascribing superscripts ${+}$ and ${-}$
we will define $x_{n}^{+}$ by the condition
$v_{n} \in (x_{n}^{+},c_{n})$.

\begin{lem}
For rotation-like maps  points $x_{n}^{-}$ lie closer to zero
than $x_{n}^{+}$.
\end{lem}

\begin{proof}
For rotation-like maps the ranges of central branches always contain 
the critical point. In particular, it means that the image of $B_{n+1}$ 
by $f_{n}^{a_{n}}$ covers both $(0)$ and the interval $(c_{n},x_{n}^{-})$.
\end{proof}
Parameters of the induction measure
sizes of domains of branches as well as and their separation
from the critical point and the boundary of the relevant box.
The distortion of these quantities 
will be controlled by bounds on $\tau$ and extendability of branches. 
Here, we provide a full list of parameters.
\begin{itemize}
\item $\alpha_{n} = \R(z_{n}^{+},z_{n}^{-};\; z_{n'}^{+}) $,
\item $\gamma_{n} = \R(x_{n}^{+},x_{n}^{-};\;0)$,
\item $\beta_{n} = \Cr(z_{n'}^{?},x_{n}^{+},x_{n}^{-},0)$
\end{itemize}
where $?$ is chosen as $+$ or $-$ so that the points have allowable ordering. 

We will examine how these quantities change after a standard inducing
step.  Generally, none of these quantities is  decreases monotonely 
in the inducing procedure. Nevertheless,  we can choose products of
these parameters that show monotone decay.  Consider
the products  $\alpha_{n}\gamma_{n}$
and $\alpha_{n}\beta_{n}$ (they correspond to the sums 
$s_{1}+s_{3}$ and $s_{2}+s_{4}$ in the complex induction).
We will see soon that  that {\em primary} (immediate) domains 
$(x_{n}^{-},x_{n}^{+})$ stay always at the definite
distance from the boundary of $B_{n'}$. This implies
that in the case of rotation-like maps these products are equivalent and 
it is enough to consider only one of them.

\begin{prop}\label{p:2}
Consider the quantity $\alpha_{n}\gamma_{n}$ for a rotation-like
sequence 
\[\phi_{0},\ldots, \phi_{n}, \ldots\; .\] 
If $\tau_{n} \leq
\tau=0.37$ for every $n$, then there is an absolute constant
$\Lambda<1$ and a fixed integer $N$ with the property that
\[ \alpha_{n}\gamma_{n} \leq \Lambda^{n} \alpha_{0}\gamma_{0}\]
for all $n>N$. 
\end{prop} 
\paragraph{Auxiliary quantities.}
Before we pass to the proof of Proposition~\ref{p:2} which will
occupy the next three subsections we will introduce three auxiliary induction
quantities $\gu_{n}$, $\omega_{n}$ and $\Omega_{n}$.
\[\gu_{n}=\R(v_{n},x_{n}^{-};\;0)\]
\[\omega_{n}=\frac{|f_{n}(0)|}{|z^{+}_{n'}|}\;\;\;\; \mbox{ and} \;\;\;\;
\Omega_{n}=\frac{|v_{n}|}{|z^{+}_{n'}|}\]
(we recall that  $v_{n}=f_{n}^{a_{n}}(0)$.)
\ \\
In addition, we have already defined $\tau_{n}$ and $\tau$. 
The estimates in the next two subsections
will be quite complicated.  It may help the reader to think
of central branches as quadratic polynomials, and of monotone branches
as affine. In this model estimates
are easier and actually give the right idea of the real situation. 
Then the distortion
 might be treated as a correction to formulae obtained in 
the ``linear-quadratic'' model.

\subsection{A non-close return} Throughout this subsection we assume
that $\phi_{n}$ makes a non-close return. In particular it means that
$(n+1)'=n$.

The distance of a point $z$ to zero is denoted by $|z|$. Observe that
 $|z_{n}^{+}|=|z_{n}^{-}|$ and $g(z_{n}^{+})=g(z_{n}^{-})$. Hence,
we will often drop superscripts $``+''$ and $``-''$ from 
the notation of distances if only no confusion can arise. 
We will start with the following simple observation.
\[ \alpha_{n+1}^{2}= 4 \cdot \R_{e}(g(z_{n+1}),g(0); g(z_{n})).\]
The image of $\R_{e}(g(z_{n+1}),g(0); g(z_{n}))$ by $(h_{n})_{*}$ is equal to 
\[ \gu_{n}\frac{\sqrt{|x_{n}^{-}||f_{n}(0)|}}
{|x_{n}^{-}|+|z_{n'}^{-}|}.\]
To find the distortion of $(h_{n})_{*}$ on this ratio, we will
complete the ratio to the cross-ratio
$\Cr_{e}(g(z_{n}),g(z_{n+1}),g(0),h_{n}^{-1}(z_{(n-2)'}^{-}))$.
Since the cross-ratio is expanded, we get

\begin{equation}\label{eq2}
\alpha_{n+1}^{2} \leq 4\cdot \gu_{n}\;
\frac{\sqrt{|x_{n}||f_{n}(0)|}}{|x_{n}|+|z_{n'}|}\; 
\frac{|z_{n'}|+|z_{n-2'}|}{|z_{n-2'}|-|f_{n}(0)|}.
\end{equation}

\begin{fact}\label{f:2} 
The distortion  of $\;\gu_{n}$ and $\gamma_{n}$ by a quadratic map is at 
least $2$.
\end{fact}
\begin{proof}
This follows directly from the definition of $\gu_{n}$.
\end{proof}
We pass to estimating $\gu_{n+1}$. 
Take the image of $\;\gu_{n+1}$ by the  quadratic map $g$. 
Fact \ref{f:2} implies that
\[\gu_{n} \leq \frac{1}{2} \cdot \R(g(v_{n+1}(0)),g(x_{n});\;g(0)).\]
Complete $g_{*}(\gu_{n+1})$ to the cross-ratio
\[\Cr\{h_{n}^{-1}(f_{n-2}(0)),g(f_{n+1}(0)),g(x_(n)),g(0)\}\]
and then push it forward by $h_{n}$. By the property of expanding cross-ratios
we have that
\begin{equation}\label{eq3}
\gu_{n+1} \leq
\frac{1}{2}\cdot \frac{|z_{n+1}|+|f_{n+2}(0)|}
{\sqrt{|f_{n}(0)|^{2}-|z_{n+1}|^{2}}}
\frac{|f_{n-2}(0)|+|f_{n}(0)|}{\sqrt{|f_{n-2}(0)|^{2}-|z_{n+1}|^{2}}}.
\end{equation}

\begin{com}\label{c1:1}
Note that the estimate~(\ref{eq3}) remains true if
we replace $f_{n+2}(0)$ by $z_{n+1}^{+}$.
\end{com}

Our next task is to combine estimates on $\gu_{n+1}$ and $\alpha_{n+1}$
and get the best possible upper bound of their product in terms
of $\gu_{n}$ and $\alpha_{n}$.
To this end we  prove

\begin{lem}\label{oszac}
For arbitrary $n$ the following inequality holds.
\[\frac{\sqrt{|x_{n}||f_{n}(0)|}}{|z_{n'}|+|x_{n}|}\;
\frac{|f_{n+2}(0)|+|z_{n+1}|}{\sqrt{|f_{n}(0)|^{2}-|z_{n+1}|^{2}}}\leq
 \frac{1}{4}\cdot \alpha_{n}\alpha_{n+1}(1+\omega_{n+2})\;
\frac{|z_{n-1}|}{|f_{n}(0)|+|z_{n-1}|}.\] 
\end{lem}

\begin{proof}
By the definition of $\alpha_{n}$  we have that
\begin{equation}\label{eq3+}
\frac{|z_{n+1}|}{\sqrt{|f_{n}(0)|^{2}-|z_{n+1}|^{2}}}
\leq \frac{1}{4}\cdot \alpha_{n+1}\alpha_{n}
\frac{\sqrt{|z_{n-1}|^{2}-|z_{n}|^{2}}}{|z_{n}|}
\sqrt{\frac{|z_{n}|^{2}-|z_{n+1}|^{2}}{|f_{n}(0)|^{2}-|z_{n+1}|^{2}}}.
\end{equation}
The last factor in the inequality $(\ref{eq3+})$  
is decreasing with respect to $|z_{n+1}|$. Thus, the right-hand
side of $(\ref{eq3+})$  is bounded by
\[\leq \frac{1}{4}\cdot \alpha_{n+1}\alpha_{n}\;\frac{|z_{n-1}|}{|f_{n}(0)|}.\]
To complete the reasoning we will need the following elementary fact:\newline 
\begin{em}
For any three positive numbers $0<x<y<z$ the inequality 
\[\frac{\sqrt{xy}}{z+x} < \frac{y}{z+y}\]
holds.
\end{em}\newline
which can be readily proved by calculus. 
>From there, 
\[\frac{\sqrt{|x_{n}||f_{n}(0)|}}{|z_{n'}|+|x_{n}|} \leq
\frac{|f_{n}(0)|}{|f_{n}(0)|+|z_{n'}|}\]
which completes the proof.
\end{proof}

\begin{com}\label{c1:2}
 Replace $|f_{n+2}(0)|$ by $|z_{n+1}|$ in the estimate of 
Lemma \ref{oszac}. By the same reasoning we obtain
\begin{equation}\label{oszac2}
\frac{\sqrt{|x_{n}||f_{n}(0)|}}{|z_{n'}|+|x_{n}|}
\frac{|z_{n+1}|}{\sqrt{|f_{n}(0)|^{2}-|z_{n+1}|^{2}}}\leq
\frac{1}{2}\cdot \alpha_{n}\alpha_{n-1}\;
\frac{|z_{n-1}|}{|z_{n-1}|+|f_{n}(0)|}.
\end{equation}
\end{com}
Multiply the inequalities $(\ref{eq2})$ and $(\ref{eq3})$ and then combine
with the inequality in  Lemma \ref{oszac}. As a result 
 we get the following  recursive formula
\[\alpha_{n+1}\gu_{n+1} \leq \lambda_{n}\alpha_{n}\gu_{n},\]
where $\lambda_{n}$ is less than
\begin{equation}\label{eq5}
\frac{1}{2}\cdot(1+\omega_{n+2})\;\frac{|z_{n'}|}{|f_{n}(0)|+|z_{n'}|}\;
\frac{|f_{n-2}(0)|+|f_{n}(0)|}{\sqrt{|f_{n-2}(0)|^{2}-|z_{n+1}|^{2}}}\;
\frac{|z_{n'}|+|z_{(n-2)'}|}{|z_{(n-2)'}|-|f_{n}(0)|}.
\end{equation}
We will bound from above $\lambda_{n}$ by maximizing $(\ref{eq5})$
with respect to a location of $f_{n}(0)$. To this end consider
\[\epsilon_{n} = \frac{|f_{n-2}(0)|+|f_{n}(0)|}
{|z_{(n-2)'}|-|f_{n}(0)||z_{n'}|+|f_{n}(0)|}\]
as a function of  $|f_{n}(0)|$ on the interval $(0,|z_{n'}|)$.

\begin{lem}\label{anal1}
The function $\epsilon_{n}$ achieves a global maximum in $0$.
\end{lem} 

\begin{proof} 
The sign of the derivative of $\epsilon_{n}$ with respect to $|f_{n}(0)|$
is the same as the sign  of
\[-(|f_{n-2}(0)|-|z_{n'}|)|z_{(n-2)'}|
+2|f_{n}(0)||f_{n-2}(0)|+|f_{n}(0)|^{2}.\]
The smaller root of the above quadratic polynomial is always
less than zero. Thus, the function $\epsilon_{n}$ can have only a local
minimum in the interval $(0,|z_{n-1}|$. Direct computation
shows that if $\tau^{2} \leq 1/3$, then 
$\epsilon_{n}(0)\geq \epsilon_{n}(|z_{n-1}|)$  
\end{proof}

%wyliczenie pochodnej
%\[\frac{d\epsilon_{n}}{d|f_{n}(0)|}=
%\frac{|z_{n-3}|-|f_{n}(0)||z_{n-1}|+|f_{n}(0)|-
%(|f_{n-2}|(0)+|f_{n}(0)|)(-2|f_{n}(0)|-|z_{n-1}|+|z_{n-3}|)}
%{(|z_{n-3}|-|f_{n}(0)||z_{n-1}|+|f_{n}(0)|)^{2}}.\]
Finally, by Lemma \ref{anal1} and the definition of $\tau$,
$\lambda_{n}$ is less than
\begin{equation}\label{iv}
\frac{1+\omega_{n+2}}{2}\; \frac{1+\tau^{2}}{\sqrt{1 - \tau^{6}}}.
\end{equation}

\begin{com}\label{c1:3}
The same computation based on Comments~\ref{c1:1} and~\ref{c1:2} yields 
\begin{equation}\label{com}
\gamma_{n}\alpha_{n} \leq 
\frac{1+\tau^{2}}{\sqrt{1 - \tau^{6}}}\gu_{n-1}\alpha_{n-1}.
\end{equation}
\end{com}
\subsection{A close return.} Throughout this subsection we assume
that $n$-th return is close. In particular it means that
$(n+1)'=n+1/2$. The scheme of the proof is much the same
as in the previous case. The only difference is that the 
reasoning is  a bit more way around and
requires several repetitions of the estimates similar 
to these found in the last subsection.

Let $\zeta_{n}$ be the ratio of $B_{n+1/2}$ to $B_{n}$.
Put $\overline{\alpha_{n+1}}= \R(z_{n+1}^{+},z_{n+1}^{-},z_{n}^{+})$.
Denote the primary preimage
of $B_{n+1/2}$ contained in $(x_{n}^{-},x_{n}^{+})$
by $(x_{n+1/2}^{-},x^{+}_{n+1/2})$.

\paragraph{Recursion.} We will write a  recursion for the sequence 
$\alpha_{n}\gamma_{n}$. By definition,
\begin{equation}\label{eq-1}
\alpha_{n+1} \leq \frac{\al_{n+1}}{\zeta_{n}\sqrt{1-\tau^{2}}}.
\end{equation}  
For $i$ ranging from $1$ to $a_{n}-1$ let $x_{n+1/2,i}^{-}$ stand for 
$f_{n}^{-i}(x_{n+1/2}^{-})$ and $x_{n+1/2,i}^{+}$ for 
$f_{n}^{-i}(x_{n+1/2}^{+})$. To bound $f_{n*}(\gamma_{n+1}\al_{n+1})$
from above
we will use similar arguments as in the previous section.
Push forward $\gamma_{n+1}$ by the quadratic map $g$. Then 
\[\gamma_{n+1}\leq \frac{1}{2} R(g(x_{n+1}^{+}),g(x_{n+1}^{-});\; g(0)).\]
Complete $g_{*}(\gamma_{n+1})$ to the cross-ratio
\[\Cr\{h_{n}^{-1}(f_{n-2}(0)),g(x_{n+1}^{+}),g(x_{n+1}^{-}),g(0)\}\]
and then push it forward by $h_{n}$. By the property of expanding cross-ratios
we have that
\begin{equation}\label{equ:15mp,3}
\gamma_{n+1} \leq
 \frac{|z_{n+1}|}{\sqrt{|f_{n}(0)|^{2}-|z_{n+1}|^{2}}}
\frac{|f_{n-2}(0)|+|f_{n}(0)|}{\sqrt{|f_{n-2}(0)|^{2}-|z_{n+1}|^{2}}}
\; . 
\end{equation}
In the same way as in the proof of Lemma~\ref{oszac} (see the inequality $(\ref{eq3+})$)
we obtain 
\[\frac{|z_{n+1}|}{\sqrt{|f_{n}(0)|^{2}-|z_{n+1}|^{2}}} \leq 
 \frac{1}{4}\cdot \al_{n+1}\alpha_{n}\;\frac{|z_{n'}|}{|f_{n}(0)|}.\] 
By definition of $\al_{n+1}$,
\[\al_{n+1} = 4\cdot \R_{e}(g(z_{n+1}),g(0);\;g(z_{n})) \leq
4 \Cr(g(z_{n}),g(z_{n+1}),g(0),h_{n}^{-1}(z^{-}_{(n-2)'})).\]
Again using the expanding property of cross-ratios we get
\[\al_{n+1}^{2} \leq 4\cdot 
\frac{|x_{n+1/2,a_{n}-1}^{-}-f_{n}(0)|}{|x^{-}_{n+1/2}|+|z_{n'}|}\; 
\frac{|z_{n'}|+|z_{n-2'}|}{|z_{n-2'}|-|f_{n}(0)|}.\]
The estimates for $\gamma_{n+1}$ and $\al_{n+1}^{2}$ lead to the
following formula

\[\al_{n+1}\gamma_{n+1}  \leq
\alpha_{n}\eta_{n}\; 
\frac{|f_{n}(0)|-|x_{n+1/2,a_{n}-1}^{-}|}
{\sqrt{|f_{n}(0)|^{2}-|z_{n+1}|^{2}}|},\]
where $\eta_{n}$ is equal to
\[\frac{|z_{n'}|}{|z_{n'}|+|x_{n+1/2,a_{n}-1}^{-}|}
\frac{|z_{n'}|+|z_{(n-2)'}|}{|z_{(n-2)'}|-|f_{n}(0)|}
\frac{|f_{n-2}(0)|+|f_{n}(0)|}
{\sqrt{|f_{n-2}(0)|^{2}-|z_{n+1}|^{2}}}.\]
Let $0\leq i < a_{n}$. We shall write $\gamma_{n,i}$ 
for the ratio 
\begin{equation}\label{equ:15mp,2}
\frac{|f_{n}^{a_{n}-i}(0)|-|x_{n+1/2,i}^{-}|}
{\sqrt{|f_{n}^{a_{n}-i}(0)||x_{n+1/2,i}|}}\; .
\end{equation}
Now, compute
\[\frac{|f_{n}(0)- x_{n+1/2,a_{n}-1}^{-}|}
{\sqrt{|f_{n}(0)|^{2}-|z_{n+1}|^{2}}} \div 
\R(g\circ f_{n}(0),g(x_{n+1/2,a_{n}-1}^{-});\; g(z_{n+1})) \leq \]
\[\leq\frac{|x_{n+1/2,a_{n}-1}^{-}|}{|x_{n+1/2,a_{n}-1}^{-}|+|f_{n}(0)|}\; .\]
Increase this factor to $1/2$ and decrease 
$\eta_{n}$ replacing\newline
 $|x_{n+1/2,a_{n}-1}^{-1}|$ with $|f_{n}(0)|$. As a result we obtain
an upper bound of $\al_{n+1}\gamma_{n+1}$ which can be written in the form
\[\frac{1}{2}\alpha_{n} \eta_{n}'\R(g\circ f_{n}(0),g(x_{n+1/2,a_{n}-1}^{-});\; g(z_{n+1})),\]
where
\[\eta'_{n} = \frac{|z_{n'}|}{|f_{n}(0)|+|z_{n'}|}\;
\frac{|f_{n-2}(0)|+|f_{n}(0)|}{\sqrt{|f_{n-2}(0)|^{2}-|z_{n+1}|^{2}}}\;
\frac{|z_{n'}|+|z_{(n-2)'}|}{|z_{(n-2)'}|-|f_{n}(0)|}.\]
By Lemma \ref{anal1} we can  only worsen estimates letting $|f_{n}(0)|=0$
in the expression for $\eta_{n}'$. Therefore,
\begin{equation}\label{eq5-}
\al_{n+1}\gamma_{n+1} \leq \frac{1}{2}\alpha_{n} 
\frac{1+\tau^{2}}{\sqrt{1-\tau^{6}}}
\R(g\circ f_{n}(0),g(x_{n+1/2,a_{n}-1}^{-}),g(z_{n+1})).
\end{equation}
Complete the last ratio to an appropriate cross-ratio by adjoining a
fourth point at $h_{n}^{-1}(f_{(n-2)}(0))$ and
then push it forward by $h_{n}$. The points $f_{n}(z_{n+1})$ and
$x_{n+1/2,a_{n}-2}^{-})$ lie on the opposite sides
of zero. The resulting cross-ratio can be only increased if
we move the point $f_{n}(z_{n+1})$ in the direction of zero. Hence,
for $a_{n}>2$
\begin{equation}\label{eq5'}
\R(g\circ f_{n}(0),g(x_{n+1/2,a_{n}-1}^{-}),g(z_{n+1})) \leq
\frac{\gamma_{n+1/2,a_{n}-2}}{{1-\tau^{2}}}
\end{equation}
If $a_{n}=2$ then put $1-\tau\Omega_{n}$ in the place of the denominator
of $(\ref{eq5'})$ in order to have the correct estimate.

We can use the sequence of estimates starting from~(\ref{equ:15mp,2})
again to prove
\begin{equation}\label{eq6} 
\gamma_{n,i} \leq \frac{1}{2}\cdot \frac{\gamma_{n+1/2,i-1}}{1-\tau^{2}}
\end{equation}
provided $i>1$ and
\begin{equation}\label{eq7} 
\gamma_{n,1} \leq \frac{1}{2}\cdot \frac{\gamma_{n+1/2,i-1}}{1-\tau\Omega_{n}}
\end{equation}
when $i=1$.
>From inequalities~(\ref{eq5-}), (\ref{eq5'}), (\ref{eq6}) and
(\ref{eq7}) we obtain 
\begin{equation}\label{eq8}
\al_{n+1}\gamma_{n+1} \leq {\overline \Lambda}_{n} \alpha_{n}\gamma_{n+1/2,0}
\end{equation}
where ${\overline\Lambda}_{n}$ is bounded from above by
\begin{equation}\label{eq9}
 \frac{1}{2^{a_{n}-1}}\frac{1+\tau^{2}}
{\sqrt{(1-\tau^{6})}(1-\tau^{2})^{a_{n}-2}(1-\tau\Omega_{n})}. 
\end{equation}
In the last step of our reasoning we exploit the fact that
$\gamma_{n+1/2,0\;}$ is substantially less than $\gamma_{n}$. We claim that
\begin{lem}\label{gniazdo}
\[\gamma_{n+1/2,0} \leq \zeta_{n}\gamma_{n}.\]
\end{lem}
\begin{proof}
We will actually prove
\[ \R(x_{n+1/2}^{+},x_{n+1/2}^{-};0)  \leq \zeta_{n}\gamma_{n}\]
and use 
\[ \gamma_{n+1/2,0} < \R(x_{n+1/2}^{+},x_{n+1/2}^{-};0)\]
which follows directly from the definition of $\gamma(n+1/2,0)$.  
First observe that $|x_{n+1/2}^{-}-x_{n+1/2}^{+}| \leq
 \zeta_{n}|x^{-}_{n}-x^{+}_{n}|$.
Indeed, the centers of $B_{n+1/2}$ and $B_{n}$ coincide and
the hyperbolic length of $B_{n+1/2}$ with respect
to $B_{n}$ is not increased by the pullback by a diffeomorphism
with a non-positive Schwarzian. Since  the element of 
the hyperbolic length is the smallest in the middle of an interval 
we conclude that pullbacks are nested with the ratio at most $\zeta_{n+1}$. 
Denote by $s_{1}$ and $s_{2}$ the centers of
$(x_{n+1/2}^{-},x_{n+1/2}^{+})$ and $(x_{n}^{-},x_{n}^{+})$.
A straightforward calculation
shows that if these  centers coincide then the Lemma follows.
Suppose that $s_{1}$ is less than $s_{2}$ since if otherwise
then we are done. Push $x_{n}^{+}$ toward $s_{1}$ so far that the centers
coincide again. This operation can only increase the ratio
of $\gamma_{n}$ to $\R(x_{n+1/2}^{+},x_{n+1/2}^{-};0)$. The ratio of
the lengths of
the resulting, concentrical intervals, is again at most $\zeta_{n}$, which
completes the proof.
\end{proof}  
Finally, Lemma $\ref{gniazdo}$ and inequalities $(\ref{eq-1})$
and $(\ref{eq8})$ imply that
\begin{equation}\label{indu}
\alpha_{n+1}\gamma_{n+1} \leq \Lambda_{n}\alpha_{n}\gamma_{n},
\end{equation}
where $\Lambda_{n}$ is less than
\[ \frac{1}{2^{a_{n}-1}}
\frac{1+\tau^{2}}
{\sqrt{(1-\tau^{6})(1-\tau^{2})}(1-\tau^{2})^{a_{n}-2}(1-\tau\Omega_{n})}.\] 
Clearly, $\Lambda_{n}$ is the largest for $a_{n}=2$ since
$\tau^{2}<\frac{1}{2}$.  
For $a_{n} = 2$ we get  
\begin{equation}\label{blad1}
\frac{1}{2}\frac{(1+\tau^{2})}
{(1-\tau\Omega_{n})\sqrt{(1-\tau^{6})(1-\tau^{2})}}
\end{equation}
as an  upper bound of $\Lambda_{n}$ in the case of a close return.  

\begin{com}
 In this subsection in the unlike for non-close returns
 we worked with the quantity $\gamma_{n}$ instead of $\;\gu_{n}$.
The estimates become stronger
if we decrease the left-hand side of $\;(\ref{indu})$ substituting
$\gamma_{n}$ by $\gu_{n}$.
\begin{equation}\label{indu1}
\alpha_{n+1}\gu_{n+1} \leq \lambda_{n} \alpha_{n}\gamma_{n},
\end{equation}
 where $\lambda_{n}$ is less than
\begin{equation}\label{blad2}
\frac{1+\omega_{n+2}}{4}\;\frac{(1+\tau^{2})}
{(1-\tau\Omega_{n})\sqrt{(1-\tau^{6})(1-\tau^{2})}}.
\end{equation}  
\end{com}

The proof is the same, except that the estimate~(\ref{equ:15mp,3})
ought to be replaced by~(\ref{eq3}). 

\subsection{Conclusion of the proof}
The sequence $\omega_{n}$ plays a crucial role in the inductive formulae 
derived in the last two subsections. The decay
of geometry depends directly on the separation of this sequence
from $1$. We will begin by writing a recursion for the sequence
$\Omega_{n}$. Clearly, $\omega_{n}\leq \Omega_{n}$.
We will consider two cases. 
\paragraph{A non-close return.} In this case $\Omega_{n}=\omega_{n}$.
\begin{lem}\label{l:3}  
Assume that $\phi_n$ makes a non-close return. Then
\[\Omega_{n+1}^{2} \leq (1+\tau^{2})\frac{\Omega_{n}+\Omega_{n+2}\tau^{2}}
{1+\Omega_{n}}.\]
\end{lem}

\begin{proof}
Take the image of $\Omega_{n}$ by the quadratic map $g$
\[\Omega_{n+1}^{2}= \R_{i}(g(v_{n+1}),g(0); g(z_{n}))\]
and next push it forward by $h_{n}$. As a result we obtain
\begin{equation}\label{a} 
f_{n*}(\Omega_{n+1})=\frac{|f_{n}(0)|+|f_{n+2}(0)|}{|f_{n}(0)|+|z_{n'}|}.
\end{equation}
To compute the distortion brought in by $h_{n}$ complete $g_{*}(\Omega_{n})$ 
to the cross-ratio
\[\Cr_{i}(g(v_{n+1}),g(x_{n+1}),g(0),h_{n}^{-1}(z_{(n-2)'}^{-})).\]
that fact that the cross-ratio is expanded. We obtain
\[\frac{|z_{n'}|+|z_{(n-2)'}|}{|z_{(n-2)'}|+|f_{n+2}(0)|}.\]
as a correction to $(\ref{a})$.
Observe that $|f_{n+2}(0)|/|z_{n'}| \leq \Omega_{n+2}\tau^{2}$ which
establishes the claim of Lemma \ref{l:3}.
\end{proof}

\paragraph{A close return.} Assume that $f_{n}$ shows a close return. 
Then $\Omega_{n+1}=|v_{n+1}|/|z_{n+1/2}|$. By definition
\[\Omega_{n+1}^{2}= \R_{i}(g(v_{n+1}),g(0); g(z_{n+1/2})).\]
Complete the last ratio to an appropriate cross-ratio by adjoining a
fourth point at $h_{n}^{-1}(z_{(n-2)'}^{-}(0))$ and
then push it forward by $h_{n}$. By the property of 
expanding cross-ratios we get
\[\Omega_{n+1}^{2} \leq \frac{|f_{n}(0)|+|f_{n+2}(0)|}
{|f_{n}(0)|+|f_{n}(z_{n+1/2})|}\;
\frac{|f_{n}(z_{n+1/2})|+|z_{(n-2)'}|}{|z_{(n-2)'}|+|f_{n+2}(0)|}.\]
Let us denote the ratio $|f_{n}(z_{n+1/2})|/|z_{n}|$ by $\sigma$.
Replace $|f_{n}(0)|$ by $|f_{n}(z_{n+1/2})|$ in the inequality
above. We obtain a new bound of $\Omega_{n+1}$ equal to
\begin{equation}\label{wyj}
\frac{1}{2}\cdot(1+\tau_{n+1}\frac{|z_{n+1/2}|}{|z_{n}|\sigma}\Omega_{n+2})
(1+\sigma\tau^{3})
\end{equation}
 Clearly, $\frac{|z_{n+1/2}|}{|z_{n}|\sigma} < 1$. 
We will estimate $\sigma$ from above under the assumption that
the critical value of $f_{n}$ remains in the central domain
for at least two iterates. 
The point $f_{n}(z_{n+1/2})$ is $f_{n}^{-a_{n}+1}(z_{n}^{-})$. Thus,
$\sigma$ is
the greatest when $f_{n}(z_{n+1/2})$ coincides with a boundary
point of $B_{n}$. The next inequality is obtained by
completing the ratio $g_{*}(\sigma)$ to an ``inclusive'' cross-ratio
with a fourth point at $h_{n}^{-1}(z_{(n-2)'})$ and
then pushing it forward by $h_{n}$. 
\[\sigma^{2} \leq \frac{|f_{n}(0)|+|z_{n}|}{|f_{n}(0)|+|z_{n'}|}\;
\frac{|z_{(n-2)'}|+|z_{n'}|}{|z_{(n-2)'}|+|z_{n}|},\]
and finally,
\begin{equation}\label{ceq1}
\sigma \leq \sqrt{\frac{2\tau(1+\tau^{2})}{1+\tau}}
\end{equation}
So, we write (\ref{wyj}) as
\begin{equation}\label{ceq0}
\Omega_{n+1}^{2} \leq \frac{1}{2}(1+\tau\Omega_{n+2})
(1+\sigma\tau^{3}),
\end{equation}
where $\sigma$ bounded by $(\ref{ceq1})$.

\paragraph{A bound.} 
We want to find an upper bound of $\Omega_{n}$. To this end, observe that
\begin{itemize}
\item  $\Omega_{n+1}$ is
 an increasing function of $\Omega_{n}$ and $\Omega_{n+2}$.
\item As long as the value of $\Omega_{n}$ is greater
than $0.7$, the estimate of Lemma \ref{l:3}  gives a lower value of
$\Omega_{n+1}$ than inequality (\ref{ceq0}).
\end{itemize}
The last statement can be easily justified by direct computation. Indeed,
the right-hand side of the estimate of Lemma \ref{l:3} is smaller than
$(1/2)(1+\tau^{2})(1+\tau^{2}\Omega_{n+2})$ while that of (\ref{ceq0}) is
larger than $(1/2)(1+\tau\Omega_{n+2})$. We will be done once we
show that
\[(1+\tau^{2})(1+0.7\cdot \tau^{2}) < 1+0.7\cdot \tau,\]
which clearly holds for $\tau < 0.4$.

So we consider the recursion given by assuming equality in (\ref{ceq0}).  
The function 
\[ y\rightarrow \sqrt{1+\frac{\tau y}{2}}\] 
has exactly one attracting 
fixed point for $y \leq 0.823562$ in the positive domain. It follows
that if $\Omega_{n+1}$ in~(\ref{ceq0}) is greater than this fixed
point, then $\Omega_{n+2}$ has to be less than $\Omega_{n+1}$. 
Thus, we set 
$\Omega =0.823562$ as an bound of $\Omega_{n}$ and $\omega_{n}$.
We note that this bound is attained in for all values of $n$
sufficiently large depending only on $\epsilon$ stipulated by Theorem
A. Indeed, by Lemma~\ref{lem:8mp,1} $\tau_{n}$ gets smaller than
$\tau=0.37$ for $n$ sufficiently large in terms of $n$, and then it is
clear that $\Omega_{n}$ decreases at least by a uniform amount for
each step of the recursion given by~(\ref{ceq0}) as long as it is
greater than $\Omega$.  

\paragraph{Final estimates.}
After these preparations we will prove Proposition \ref{p:2}.
We will conduct estimates by  splitting the sequence of
box mappings into blocks. Each block save perhaps the first
will have a mapping with a close return immediately followed by 
a maximal sequence of consecutive non-close returns. E.g., a single
mapping with a close return is a block. 
The only exception from the above rule of constructing  blocks 
occurs when $\phi_{0}$ makes a 
non-close return. Then the first block consists of a maximal sequence of
box mappings with non-close returns.

Below we list the rules which will give recursive estimates within 
a given block of box mappings. 
\begin{enumerate}
\item  Suppose $\phi_{n}$ exhibits no close return
and is the last such mapping in a given block. Then we use formula 
$(\ref{com})$ to estimate   
\[\alpha_{n+1}\gamma_{n+1} \leq 1.13837 \cdot \alpha_{n}\gu_{n},\]
If $\phi_{n}$ is not last in its block, we use the inequality $(\ref{iv})$.
\[\alpha_{n+1}\gu_{n+1} \leq 0.56919 \cdot(1+\omega_{n+2})\; \alpha_{n}\gu_{n}.\]
\item Let $n$-th  box mapping exhibits close return.\newline
If $\phi_{n}$ is not a block in its own right, then we apply 
formula~(\ref{blad1}):
\[\alpha_{n+1}\gu_{n}\leq  0.80402 \cdot  \alpha_{n}\gamma_{n}.\]
If $\phi_{n}$ is a block by itself, then~(\ref{blad2}) implies that

\[\alpha_{n+1}\gamma_{n+1} \leq 0.881181 \cdot \alpha_{n}\gamma_{n}.\]
\end{enumerate}

We will consider two cases:

\subparagraph{Blocks with at least two box mappings with non-close returns.}
Suppose that a series of at least two box mappings with non-close
returns begins at the moment $n$. We will show that the separation
of the critical value $f_{n}(0)$ from the boundary of the box $B_{n'}$
improves with $n$ growing.
Indeed, by Lemma \ref{l:3}
\[\Omega_{n+1} \leq (1+\tau^{2})\sqrt{\frac{\Omega}{1+\Omega}}
\leq 0.76403\]
and next
\[\Omega_{n+2}\leq \sqrt{(1+\tau^{2})\frac{0.76403+\Omega\tau^{2}}
{1.76403}} \leq 0.751714.\]
By monotonicity of this formula with respect to $\Omega_{n}$, all  
$\Omega_{k} \leq 0.751714$
for $k \geq n+2 $. Consequently, for $k \geq n$ we obtain  that
\[\alpha_{k+1}\gu_{k+1} \leq 0.9971 \cdot \alpha_{k}\gu_{k}.\]
It happens that $0.80402\cdot 1.13837 <1$. 
Thus, if a block of length $k\geq 3$ starting with $\phi_{n}$ is taken
as whole, then
\[ \alpha_{n+k}\gamma_{n+k} \leq (0.9971)^{k/3} \alpha_{n}\gamma_{n}\;
.\]

\subparagraph{Shorter blocks.}
For a block of a single map with a close return, or one close and one
non-close, it follows immediately from our rules that the product
$\alpha_{n}\gamma_{n}$ decreases after passing through a block by a
definite constant less than $1$.

\subparagraph{Conclusion.} 
To see that $\alpha_{n}\gamma_{n}$ goes down to $0$ at least
exponentially fast, first wait $N$ steps for bounds on $\Omega$ to be
achieved ($N$ is bounded in terms of $\tau$). 
Then pick a $k>2N$, and construct the blocks starting from $\phi_{N}$.
Cut off the last block at $\phi_{k}$. The uniform exponential estimate
follows at once from our considerations of the rate of decay within
blocks. So, Proposition~\ref{p:2} follows.   

\subsection{Decay of box geometry}
\paragraph{General picture.}
In this subsection we will estimate the rate of the decay of box geometry
proving eventually that for all S-unimodal rotation-like maps 
the rate is always at least exponential. This will prove Theorems A
and C. 
In the course of inducing a subtle interaction between $\alpha_{n}$ and
$\gamma_{n}$ takes place. Namely, after a long series of non-close
returns, $\gamma_{n}$ is approximately equal to the second power
of $\alpha_{n}$. The first close return will violate
this simple relation between $\gamma_{n}$ and $\alpha_{n}$ by
decreasing  $\alpha_{n}$ stronger than
$\gamma_{n}$. If the close return is deep enough (i.e. the critical
value needs a lot of iterates to escape from the central domain)
then $\gamma_{n}$ and $\alpha_{n}$
can even become comparable. At the moment when we are leaving box maps
with close returns, $\gamma_{n}$ and $\alpha{n}$ will quickly regain 
(exponentially fast with the number of steps of non-close inducing) their
square-law relation. However, the product $\alpha_{n}\gamma_{n}$
for rotation-like maps of bounded type decreases 
asymptotically with each step of inducing by a constant uniformly
separated from $0$ and $1$. It means
that while switching between patterns of inducing
``oscillations'' between $\alpha_{n}$ and $\gamma_{n}$ destroy
monotone (known from the Fibonacci example) fashion of the decay of boxes. 
In particular, $\alpha_{n+1}/\alpha_{n}$ can become arbitrarily large
in some cases. 

Theorem C expresses what we mean by an exponential decay of 
box geometry. From Theorem C it follows immediately
that there is  a whole class of S-unimodal maps with at most
exponential decay of box geometry. The dynamics of maps from this class,
purely characterized in terms of rotation number, is 
certainly different from the ``Fibonacci pattern''.

\begin{coro}
For any S-unimodal rotation-like map with 
a rotation number of the constant type\footnote{Let us recall that
 a number $\rho$ is of the constant type if and only if all coefficients 
in its continuous fraction representation are bounded.} there exist constants
$K>0$ and $0<\kappa <1$ so that
\[\tau_{n}> K \kappa^{n}.\]
The constant $\kappa$ depends only on the upper bound of the coefficients
of the continued fraction representation of the rotation number while
the constant $K$ depends solely on the initial geometry of $\phi_{0}$,
in particular in uniformly controlled by the parameter $\epsilon$ of
Theorem A.
\end{coro} 
>From now on we will denote positive constants dependent only on the
initial geometry by $K$ and call them uniform.
 Whenever confusion of the type $K<K$ can arise we will distinguish
constants $K$ by adding appropriate  subscripts. 

\paragraph{Proof of Theorem C:}
We will start the proof with two Lemmas.
\begin{lem}\label{gamma}
The is a uniform constant $K$ so that
\[\gamma_{n}\leq K\alpha_{n}.\]
\end{lem}

\begin{proof}
The proof follows from the observation that
the image of $\gamma_{n}$ by $f_{n}$ is comparable with $\gamma_{n}$
and $f_{n*}(\gamma_{n})$ can exceed $\alpha_{n}$ by no more than a
uniform constant. 
\end{proof}

\begin{lem}\label{tau}
The ratio $\tau_{n}$ of two consecutive boxes  goes
to zero at least exponentially fast.
\end{lem}
\begin{proof}
Suppose that $n$ is so large that $\tau_{n} < \tau=0.37$
Let us recall that $\alpha_{n}\gamma_{n}$ goes to zero
at least exponentially fast. We will actually prove that $\alpha_{n}$
decreases exponentially which is easily equivalent. 
 Consider two cases:
\begin{itemize}
\item $f_{n}$ shows a close return.\newline
We push forward $\alpha_{n}$ by $f_{n*}^{a_{n}}$. Using (\ref{eq6})
 we obtain
\begin{equation}\label{trzeci}
\alpha_{n+1}^{2} \leq K\lambda^{a_{n}}\gamma_{n}\; ,
\end{equation}
where $\lambda <1$ depends only on $\tau$. Finally, by Lemma \ref{gamma}
\[\alpha_{n+1}^{4} \leq K \lambda^{2a_{n}}\alpha_{n-1}\gamma_{n-1}\]
Proposition~\ref{p:2} concludes the proof. 

\item $f_{n}$ shows a close return.\newline
Similarly as before we get that
\[\alpha_{n+1}^{2} \leq K \gamma_{n}.\]
Lemma~\ref{tau} likewise follows.
\end{itemize}
\end{proof}

Lemma \ref{tau} together with the inequality $\ref{trzeci}$
give the upper estimate of Theorem C.

To prove the opposite estimate make
$\tau$ go to $0$ arbitrary small in all distortion estimates. 
 So, we can reverse the directions
of the inequalities estimating $\alpha_{n}\gamma_{n}$ from
below. Next, observe that  $\lambda_{n}$ and $\Lambda_{n}$
appearing in the recursive scheme for 
for $\alpha_{n}\gamma_{n}$ are after a finite number
of inducing step greater than $(\frac{1}{2}
-\epsilon)^{a_{n}}$. This completes the proof of Theorem C.

Theorem A follows directly from Theorem C and Lemma~\ref{lem:8mp,1}.

We finish this section with an important technical statement
concerning the starting condition. 
    
\begin{prop}\label{sta}

For any type I real box map $\phi$ which satisfies
\[ \frac{|B_{n}|}{|B_{n'}|} \leq \tau\]
 and every $\delta>0$ there is a
number $k(\tau,\delta)$ only depending on its stated parameters so
that if only
a close return occurs and the 
critical value remains in the domain of the central branch for more
than $k$ iterations, then the starting condition
is satisfied with norm $\delta$. 
\end{prop}

\begin{proof}
Define $B^{i} = f^{-i}(B_{n})$ where $f$ is the central branch of
$\phi$. Also, suppose that $\phi$ is of rank $n$. 
  By the 
non-positive Schwarzian property,
 in each component of $B_{n}\setminus B_{n+1/2}$ there is at most one point at
 which the derivative of $f_{n}$ is equal to $1$. This point is
between fixed points of $f$, or there would be a restrictive interval.
Hence, there exists a uniform constant $K_{1}(\tau)$ so that for all 
$1<i < a_{n}$ 
\[\frac{|B^{i}\setminus B^{i-1}|}{|B_{n}|} \leq \frac{K_{1}(\tau)}{i}.\]
Push forward $\alpha_{n+1}$ by $f_{n}$. The bounded distortion yields
\[\alpha_{n+1} \leq K_{2}(\tau) \sqrt{\zeta_{n}}.\]
On the other hand (see $(\ref{eq-1})$)
\[\alpha_{n+1}\leq K_{3}(\tau) \al_{n+1}/\zeta_{n}.\]
By the real K\"{o}be Lemma and the definition of $\al_{n+1}$
\[\al_{n+1} \leq K_{4}(\tau) \sqrt{\zeta_{n}/a_{n}}.\]
Combining the above inequalities we get finally that
\[\alpha_{n+1} \leq K_{5}(\tau)/\sqrt{a_{n}},\]
which completes the proof of the Proposition.
\end{proof}

\section{Initial bounds}
Our goal is to prove that $0.37$ works
as an eventual bound of  $\tau_{n}=|B_{n}|/|B_{n-1}|$. The estimates will
be uniform in the sense of Theorem A. This will prove
Lemma~\ref{lem:8mp,1}. The reasoning  naturally splits into two parts.
In the first we will find very initial estimates by solving a certain extremal
problem for so called 'auxiliary inducing'. 
Then using these bounds to control
distortion we will refine previous estimates referring
to the particular features of both types of inducing.

\subsection{Geometrical setup}
\paragraph{Notations.}
The new domain $B_{n+1}$ is formed as the preimage by $f_{n}$ of an
immediate branch of $\phi_{n}$ which is always filled-in so that it
maps onto $B_{(n+1)'}$. Let us call this domain $(u_{n}, v_{n})$ and
let us say that $u_{n}$ is closer to the critical point. The branch
defined on $(u_{n},v_{n})$ extends at least onto $B_{n'}$ as the
image. The domain of this extension will be called $(s_{n},t_{n})$,
and again say that $s_{n}$ is closer to the critical point. Then, call
$B_{(n+1)'}$ $(-w_{n}, w_{n})$ and say that $w_{n}$ is on the side of
the critical value. Lastly, let $(-w_{n-1}, w_{n-1})$ be
$f_{n}^{-1}((-w_{n}, w_{n}))$. The reader may try to get familiar with
this notation by trying to see that the ordering of points is
\[ -w_{n-1}, -w_{n}, 0, w_{n}, s_{n}, u_{n}, v_{n}, t_{n}, w_{n-1}\]
or perhaps the reverse. These notations are applicable for both close
and non-close returns.    

Observe that 
\[\Cr_{i}(s_{n},u_{n},v_{n},t_{n-1})\leq 
\Cr_{i}(-w_{n-1},-w_{n},w_{n},w_{n-1})\; .\]
Let us mention now that the size
of the interval $(s_{n-1},t_{n-1})$ with respect to $(w_{n},w_{n-1})$
will play an important role later on when we refine the first estimates
obtained by assuming in estimates that $(s_{n-1},t_{n-1})$ is equal to
$(w_{n},w_{n-1})$.

\paragraph{Formulation of extremal problem.} Let us recall that $f_{n}$ can
be represented as a composition of the quadratic map $g$ and
the diffeomorphism $h_{n}$. The  extendability properties of $f_{n}$
are formulated in Lemma \ref{ext}. We will denote preimages of points
by $h_{n}$ by adding primes to the notation. 

>From the previous two paragraphs,
\begin{equation}\label{s0}
\tau_{n}^{2} \leq \frac{|u'_{n}|- |v'_{n}|}{|w'_{n}|+|v'_{n}|}.
\end{equation}
 Complete the ratio $|w'_{n-1}|-|w'_{n}|/|w'_{n-1}|+|w`_{n}|$
to the cross-ratio 
\[\Cr_{e}(-w'_{n-1},w'_{n},w'_{n-1},(z^{-}_{(n-2)'})')\; .\]
By the expanding property of cross-ratios, 
\[\frac{|(z^{-}_{(n-2)'})'-w'_{n}|}{|(z^{-}_{(n-2)'})'+|w`_{n}|}\leq 
\frac{1-\alpha}{1+\alpha}\frac{1+\nu_{n}^{2}}{1-\nu_{n}^{2}},\]
where $\alpha= |w'_{n}|/|w'_{n-1}|$. Denote the right-hand side of the
above inequality by $1/L$.
To find preliminary bounds we solve the following
extremal problem:
\begin{prob}
Suppose that the inclusive cross-ratio of the interval
$(u'_{n},v'_{n})$ with respect to $s'_{n},t'_{n-1}$ is equal to
to $C$ and the cross-ratio 
\[ \Cr_{e}(-s'_{n-1},u'_{n},v'_{n-1},w'_{n-3})\; .\]
is equal to $1/L\beta$, $\beta\geq 1$. The interpretation of
$\beta$ is that it 
accounts for the nesting of $(s_{n}, t_{n})$ in  $(w_{n}, w_{n-1})$ as
well as for possible additional extendability of $f_{n}$.

Find a maximum of
\[\frac{|u'_{n}|- |v'_{n}|}{|w'_{n}|+|v'_{n}|}.\]
depending on the location of points $u'_{n}$ and $v'_{n}$.
\end{prob}

\begin{sol}
We assume that the interval $s'_{n-1},t'_{n-1}$ has the unit length.
Then $|w'_{n-1}|+|w'_{n}| > L$. Denote $|s'_{n-1}|-|u'_{n}|$
by $\rho$ and $|v'_{n}-w'_{n}|$ by $\lambda$. Then
\begin{equation}\label{s1}
C=\frac{\rho+\lambda-1}{\rho\lambda}.
\end{equation}
We want
to maximize
\[\frac{\rho +\lambda-1}{L+\lambda}= C \frac{\rho\lambda}{L+\lambda},\]
 which is equivalent to minimizing
\[M=\frac{1}{\rho}(1+\frac{L}{\lambda}).\]
In the extremal point the gradients of $M$ and $C$ are parallel. Hence
\begin{equation}\label{s2}
\frac{1-\lambda}{1-\rho} = 1+\frac{\lambda}{L}.
\end{equation}
Calculate $\rho$ from~(\ref{s1}) and  substitute 
into~(\ref{s1}). As we solve the resulting quadratic equation for
$\lambda^{-1}$ we get
the first coordinate $\lambda_{0}$ of the extremal point.
\[\lambda_{0}^{-1}=1+\sqrt{(1-C)(1+L^{-1})}.\]
By algebra we calculate the second coordinate and then
 the minimum of $M$ 
\[L(\sqrt{1-C}+\sqrt{1+L^{-1}})^{2}.\]
Finally, we get   
\[\frac{C}{(\sqrt{1-C}+\sqrt{1+L^{-1}})^{2}}.\]
as a solution of our problem.
Let $T=1-\nu_{n}^{2}/1+\nu_{n}^{2}$ and $\gamma=\sqrt{1-C}$.
Then by inequality $(\ref{s0})$
\[\tau_{n+1} \leq 
\frac{\sqrt{1-\gamma^{2}}}{\sqrt{T\beta\gamma}+\sqrt{1+\frac{\beta T}{\gamma}}}.\]
\end{sol}

Set $T'= 1-\tau_{n+1}^{2}/1+\tau^{2}_{n+1}$. If $T'>T$               
then $\tau_{n+1}<\mu_{n}=\sqrt{\tau_{n-1}\tau_{n-2}}$. Observe that
$T'$ is a growing function of $T$. We will determine when
the difference $T'-T$ is positive. To this aim we bring $T'-T$ to the
common denominator and examine the sign of the numerator.
\[T((1-T)\beta(\gamma+\frac{1}{\gamma})-2) +\gamma_{2}(1+T)+2T\beta(1-T)
\sqrt(1-T)\sqrt{1+\frac{\gamma}{T}}.\]

\subsection{Computation of initial bounds}
\paragraph{Basic procedures.}
\subparagraph{General description.}
We begin with the condition $T' - T > 0$ which can be analytically
rewritten as
\begin{equation}\label{equ:1ma,1} 
T((1-T)\beta(\gamma+\frac{1}{\gamma)} - 2) + \gamma^{2}(1+T) +
2\beta T(1-T)\sqrt{1+\frac{\gamma}{\beta T}} > 0 \; .
\end{equation}

The first procedure aims to find a possibly large $T$ independent of
$\gamma$ in some range, but depending on $\beta$, which will force
this estimate to be fulfilled. 

The second condition will find $T'$ from the formula
\begin{equation}\label{equ:1ma,2}
T' = \frac{t(\gamma + \frac{1}{\gamma} +
2\sqrt{1+\frac{\gamma}{t}}) + \gamma^{2}}
   {2-\gamma^{2} + t(\frac{1}{\gamma} + \gamma + 2\sqrt{1 +
\frac{\gamma}{t}})}  
\end{equation}
where $t = \beta * T$. More precisely, an upper bound will be found
depending on $t$, but not depending on $\gamma$ varying in some
specified range. 

We also obtain as a corollary:
\begin{fact}\label{fa:9mp,1}
Let $\phi$ be a type I box mapping of rank $n$, let 
$|B_{n}|/|B_{n'}|\leq 1-\epsilon$ and assume that the central branch
of $\phi$ is $\epsilon$-extendable. If another type I box mapping of
rank $m$ is obtained from $\phi$ in a number of standard inducing
steps, then 
\[ \frac{|B_{m}|}{|B_{m'}|} \leq 1 - K(\epsilon)\]
where $K$ is a continuous function of $\epsilon$ only, positive when 
$\epsilon>0$. 
\end{fact}
\begin{proof}
>From Fact~\ref{fa:9mp,1} we see that the ratio will remain bounded
away from $1$ for two first standard inducing steps. Then,
formula~\ref{equ:1ma,1} implies that it will not deteriorate as long
as it is close to $1$ (the condition~\ref{equ:1ma,1} is clearly
satisfied for $T$ close to $0$.) 
\end{proof}

\subparagraph{Implementation of the first procedure.}
The first procedure will take five parameters. $\gamma_{u}$ and
$\gamma_{l}$ will give an upper and lower bound of the allowed range
of $\gamma$. For practical reasons, we assume $\gamma_{l} \leq 0.1$, but one
easily checks that the condition of~[\ref{equ:1ma,1}] is satisfied for
any $\gamma < 0.1$ and $0\leq T\leq 1$. Another parameter $\sigma$
gives the step size. We will cover the range $[\gamma_{l},\gamma_{u}]$
with finitely many closed intervals of length $\sigma$. We will
estimate from below the range of positive values of $T$ which
satisfy~[\ref{equ:1ma,1}] on each subinterval, and finally take the
maximum of all estimates with respect to the subintervals. Another
parameter $\beta$ whose meaning is clear. Lastly, we have $\nu$ 
which must no less than the answer (thus $\nu=1$ will always work,
but the point is to sharpen the estimate by picking $\nu$ just about
as small as possible.)

To finish the description, we have to explain how the lower estimate
is found on a subinterval $[\gamma_{1},\gamma_{2}]$.
The left-hand side of~[\ref{equ:1ma,1}] is clearly bounded from below
by
\[ T((1-T)\beta(\gamma_{2}+\frac{1}{\gamma_{2}}) - 2) + \gamma_{1}^{2}(1+T) +
2\beta T(1-T)\sqrt{1+\frac{\gamma_{1}}{\beta \nu}} > 0 \; .\]
This gives us a quadratic inequality on $T$ which is solved
algebraically to give us the answer.

\subparagraph{Implementation of the second procedure.}
This procedure will take four parameters. As previously, $\gamma_{u}$
and $\gamma_{l}$ will give the range of $\gamma$ with respect to which
the $T'$ must be minimized. Again, $\gamma < 0.01$ will give an answer
greater than $0.9$ which is better than we will ever use, so we assume
$\gamma_{l} \geq 0.01$. We use the same procedure of dividing
$[\gamma_{l},\gamma_{u}]$ into subintervals, taking a lower bound for
$T'$ on each interval, and taking the minimum with respect to all
subintervals for a final answer. The parameter $\sigma$ gives the
length of subintervals. The parameter $t$ is $\beta T$. 

On each subinterval $[\gamma_{1}, \gamma_{2}]$,
formula~[\ref{equ:1ma,2}] bounds $T'$ from below by 
\[   T' = \frac{t(\gamma_{2} + \frac{1}{\gamma_{2}} +
2\sqrt{1+\frac{\gamma_{1}}{t}}) + \gamma_{1}^{2}}
   {2-\gamma_{1}^{2} + t(\frac{1}{\gamma_{1}} + \gamma_{1} + 2\sqrt{1 +
\frac{\gamma_{2}}{t}})} \; .\] 
This is evaluated on each subinterval.   
  
\paragraph{The first estimate on $\tau$.}
The first estimate on $\tau$ is obtained by calling the first
procedure with parameters $\gamma_{l} = 0.1$, $\gamma_{u} = 1$,
$\beta=1$, $\mu = 0.7$, $\sigma = 10^{-5}$. The result is more than 
$T>0.69901$. 
Since 
\[ T = \frac{1 - \tau^{2}}{1 + \tau^{2}}\; , \]
$\tau(1) = 0.4209$ will guarantee that $T'-T>0$. We claim that from
any initial box mapping the value of $\tau$ will eventually drop below
$\tau(1)$ in a uniformly bounded number of inducing steps. Indeed,
from~(\ref{equ:1ma,1}) by a compactness argument it is clear that if 
$0 < \epsilon < T < 0.69901$, the increment $T'-T$ is bounded away
from $0$. However, $\epsilon$ is uniformly bounded away from $0$. 
On the other hand, the formula~(\ref{equ:1ma,2}) gives $T'$ as an
increasing function of $T$, thus a decreasing function of $\tau$.
Hence, once $\tau$ gets below $\tau(1)$ it will stay there. We
conclude that eventually, after a number of inducing steps bounded in
terms of $\epsilon$ only, the box ratios $\tau$ become smaller that
$\tau(1)$.  

\paragraph{Better estimates on the box ratio.}
Here we concentrate on a complex box mapping $\phi$. We call $b_{n}$ the
central domain of $\phi$. We assume that an estimate $\tau$ for the
box ratio
is already satisfied by $\phi$ as well as two preceding mappings in
the inducing sequence. We will assume that $\tau\leq tau^{1}$.
We will try to get better estimates for the
next box ratio $|b_{n+1}|/|b_{(n+1)'}$. The procedure will depend on
what $\phi$ does. 

\subparagraph{A close return for $\phi$.}
In this case, we see that opportunity to significantly improve $\beta$
in formula~(\ref{equ:1ma,2}). Indeed, $\beta=1$ corresponds to the
assumption that either of intervals between $B^{k}$ and $B^{k-1}$ maps
onto $B_{k-1}$. In fact, however, we know that each of them contains
an extended branch mapping of rank $n'$. By standard estimates, we get
that the branch mapping onto $B^{k-1}$ occupies the fraction of either
space between $B^{k}$ and $B^{k-1}$ which is no more than $\tau$.
Thus, one can take $\beta=tau^{-1}$. Next, one uses procedure two with
\[ t = \beta\frac{1-\tau^{2}}{1+\tau^{2}}\]   
$\sigma = 10^{-5}$ and the full range of $\gamma$ from $0.01$ to $1$.
For $\tau=\tau(1)$ this gives $\tau(21) = 0.29728$.

\subparagraph{A non-close return following a close return.}
In this case, we will use formula~(\ref{equ:1ma,1}) and get
an improvement from two sources. First, the lower estimate for
$\gamma$ will be quite large as a result of the box ratio being very
small in the case of a close return. Secondly, the extendability
factor normally given by 
\[ \frac{1-\tau^{2}}{1+\tau^{2}}\] will also grow, since the ratio of 
the length of $B^{k}$ to the length of $B_{(n-2)'}$ will be at most 
the ratio between lengths of $B^{1}$ and $B_{n}$ times $\tau^{2}$.
This last ratio $b$ is given by
\[ b=\sqrt{\frac{2\tau}{1+\tau}(1+\tau^{2})} \; .\]
Now $\beta$ is at least
\[ \beta \geq \frac{1-bd^{2}}{1+bd^{2}}\cdot\frac{1+d^{2}}{1-d^{2}} \]
where $d \leq \tau$. Thus, if we want to get a better estimate on the
box ratio, $d$ should be no more than this expected lower estimate. 
Using $d=0.37$ and $\tau=\tau(1)$ we get $\beta \geq 1.04694$. 
Also, $\gamma_{l} \geq \frac{1 - \tau(21)}{1 + \tau(21)}$. With
$\gamma_{u}=1$, $\sigma=10^{-5}$ and $t=0.749$, procedure one gives  
$T = 0.7483$ corresponding to $\tau(22) = 0.3795$. The meaning of this
result is that once that estimates apriori given by $\tau(1)$ and
$\tau(21)$ hold, in this case $\tau$ the box ratio will be less than
$\tau$ as long as $\tau \geq tau(22)$. This would imply that
eventually $\tau$ gets smaller than $\tau^{22}$, but we have one more
case.

\subparagraph{A non-close return followed by a non-close return.}
In this case we also use procedure one. The improvement is obtained
from a better $\gamma_{l}$ given by 
\[ \gamma_{l} = \frac{1 - \tau(1)}{1 + \tau(1)}\]
as well as better $\beta$. Here, $\beta$ can easily be estimated
\[ \beta \geq \sqrt{\frac{1 + \tau}{2\tau(1+\tau^{2})}}\]
With $\sigma=10^{-5}$, $\gamma_{u}=1$ and $t=0.749$ we get 
$T=0.7449$ which gives $\tau(23) = 0.38236$. 

Now, the combination of estimates in all cases implies that the box
ratio will eventually decrease below the maximum of $\tau^{21}$,
$\tau^{22}$ and $\tau^{23}$ which happens to be $\tau^{23}$. Again, we
argue that this will happen in a uniformly bounded number of inducing
steps. So, we can take $\tau(2) = \tau(23) = 0.38236$.

\subparagraph{A second round of estimates.}
We now repeat the same sequence of estimates in all three cases using
$\tau^{2}$ as our original bound instead of $\tau^{1}$. In the case of
a close return this gives
\[ \tau(31) = 0.27736\; .\]
In the second situation, we get $b\leq 0.79629$. With this, and
$d=0.369$, we obtain $\beta\geq 1.05793$. Also, $\gamma_{l} = 0.56572$.   
With $t=0.76$ and $\sigma=10^{-5}$ procedure one yields
$\tau(32)=0.36983$.
In the last case, we get $\gamma_{l} = 0.4468$ and $\beta=1.25582$. We
feed those into procedure one together with $\gamma_{u}=1$ and
$t=0.78$ to get $T=0.75983$ which corresponds to $\tau(33) = 0.36942$.

We see that indeed the box ratio eventually goes below $0.37$ which
proves Lemma~\ref{lem:8mp,1}.

\section{Complex induction}
\subsection{Non-decreasing moduli}
\paragraph{Statement of the result.}
Suppose now that a complex box mapping $\phi$ is given of rank $n>1$
which later undergoes $k$ consecutive steps of general inducing in the
box case. We denote by $\phi_{i}$ the complex box mappings obtained in
the process so that $\phi_{0}=\phi$ and $\phi_{i+1}$ arises from
$\phi_{i}$ in a general inducing step. 
\begin{prop}\label{prop:21kp,1}
Choose $0<j<k$. Suppose that a constant $\beta$ can chosen
independently of $B$ so that normalized critical
symbols can be chosen with norm $\beta$ for all monotone branches $B$
whose domains intersect the real line. Then, for $\phi_{j+1}$ all
normalized critical symbols can be constructed
with the same norm $\beta$. 
Assume in addition that $j\geq 3$ and 
$\phi_{j}$ does not show a rotation-like return. 
If the critical point is in the range of the real central branch, let
$l_{j}$ mean the number of
consecutive images of the critical point by the central branch of
$\phi_{j}$ which
remain in the central hole. Otherwise, set $l_{j}=0$. 
Put $l$ equal to the maximum of $l_{j}$,
$0\leq j<k$. 
Then, there is a positive function $K(l)$
such that normalized 
symbols can be constructed for all univalent branches of $\phi_{j+1}$
with norm $(1+K(l))\beta$.
\end{prop}

The first part of Proposition~\ref{prop:21kp,1} which says that the symbols
exist with norm $\beta$ follows from the proof of Lemma 3.4 of~\cite{hyde}.    
We will repeat the proof here and refine the argument to show the
second part of the statement.

\paragraph{An outline.}
The proof of Proposition~\ref{prop:21kp,1} has to be split into a number of
cases. The major dichotomy is between close returns and others.
We remind the reader that the situation is classified as a close
return if the critical value is in the central branch.
As analytic tools, we will use the behavior of moduli of annuli under
complex analytic mappings. Univalent maps transport the annuli
without a change of modulus, analytic branched covers of degree $2$
will at worst halve them, and for a sequence of nesting annuli their
moduli are superadditive (see~\cite{lehvi}, Ch. I, for proofs, or
~\cite{brahu} for an application to complex dynamics.) 

We assume that a mapping $\phi_{j}$ of rank $n$ is given as in the hypothesis of
Proposition~\ref{prop:21kp,1}. We will construct $\phi_{j+1}$ and show
that necessary separation estimates. On the level of notation,
quantities related to $\phi_{j+1}$ will be distinguished by writing a
bar above them.  
\subsection{Non-close returns}
\paragraph{Dynamical classification of branches.} First, we classify branches
of $\phi_{j+1}$ according to their parent branches (compare the description of
the inducing construction for the definition of parent branches.) The main
split is between immediate parent branches and non-immediate parent branches.
Among domains with immediate parent branches we distinguish maximal branches,
or immediate preimages of $B_{n+1}$, and others. Otherwise, a non-immediate
parent branch gets mapped forward by $\psi_{j}$. This gives a univalent branch
of $\phi_{j}$, denoted with $B'$. We also have another univalent domain of
$\phi_{j}$ which contains the critical value. This will be denoted with $B$.
Then, we distinguish three subcases according to whether $B$ and $B'$ are
independent, or one is subordinate to the other.  The results of our
computation can be summarized as follows.

\begin{lem}\label{lem:28np,1}
Suppose that a sequence of type I complex box mappings $\phi_{j}$ is given which
 satisfies the assumptions of Proposition~\ref{prop:21kp,1}. Let $1\leq j \leq k$ and assume that $\phi_{j}$ shows a non-close return. If $\beta$ is the separation norm of $\phi_{0}$, then normalized critical symbols can be constructed for  branches of 
\begin{itemize}
\item
for  immediate preimages of $B_{n+1}$ the symbol has norm $\beta$ with corrections
\[ \overline{\lambda}_{1} = \frac{\lambda_{2}(B)}{2}\: ,\: 
   \overline{\lambda}_{2} = \frac{\lambda_{1}(B)}{2}\; . \]
\item
for domains whose parent branches are not immediate, and whose domains $B'$ described above are independent with the postcritical domain $B$, the norm is $\beta$ with corrections
\[ \overline{\lambda}_{1} = \frac{\lambda_{2}(B)}{2}\; ,\]
\[ \overline{\lambda}_{2} = \frac{\alpha - \delta}{2}\]
where $\delta = \max( - \lambda_{2}(B), - \lambda_{2}(B'))$. 
\item
in all other situations the norm is at least $(1+K(l))\beta$ with $K$ a positive function of $l$ only ($l$ is defined in the statement of Proposition~\ref{prop:21kp,1}.)
\end{itemize}
\end{lem}

We proceed to prove Lemma~\ref{lem:28np,1}.
\paragraph{Immediate preimages.}
Let $B$ denote the hole which contains the critical branch. In all cases
the new central hole $B_{n+1}$ is separated from the boundary of
$B_{n}$ by an annulus of modulus at least $(\beta+\lambda_{2}(B))/2$.
We will first construct the symbols for immediate preimages of
$B_{n+1}$, meaning the preimages by the central branch. Naturally,
these preimages exist on the real line exactly if the image of the
real central branch covers the central domain. The annulus
$\overline{A}_{2}$ around $B_{n+1}$ will be the preimage by the
central branch of the encompassed by $A_{3}$ around $B$ with $B$
removed. This region consists at least of the union of $A_{3}$ and
$A'$. Then, $\overline{A}_{1}$ is the preimage of $A_{4}$. It
follows that we can take
\[ \overline{s}_{1} = \frac{\beta+\lambda_{2}(B)}{2}\; \mbox{and}\]     
\[ \overline{s}_{2} = \frac{\beta-\lambda_{1}(B)}{2} \; .\]
Of course, since components of the symbol are only lower estimates, we
are always allowed to decrease them if needed.
The annulus $\overline{A}'$ is naturally given as the preimage of the
annulus between $B_{n+1}$ and the boundary of $B_{n}$ by the central
branch, likewise $\overline{A}_{3}$ is the preimage of $A_{2}$, and
$\overline{A}_{4}$ is the preimage of $A_{1}$. Since the first two
preimages are taken in an univalent fashion, we get
\[ \overline{s}_{3} = \frac{\beta+\lambda_{2}(B)}{2} + \alpha -
\lambda_{2}(B)  \; \mbox{and}\]
\[ \overline{s}_{4} = \overline{s}_{3} +
\frac{\lambda_{1}(B)+\lambda_{2}(B)}{2} = \frac{\beta}{2} + \alpha +
\frac{\lambda_{1}(B)}{2}\; .\]
Thus, if we put 
\[ \overline{\lambda}_{1} = \frac{\lambda_{2}(B)}{2}\: ,\: 
   \overline{\lambda}_{2} = \frac{\lambda_{1}(B)}{2} \]
This gives the critical symbol for immediate preimages.

\paragraph{Reduction to maximal branches.}
A branch $\overline{B}$ can be either maximal, i.e. the  preimage of
$B_{n+1}$ by its parent branch, or it can be inside another domain of rank $n'$
nested inside the parent domain. We now argue that it is
sufficient to do the estimates for maximal branches. In the process we describe, the annulus
$\overline{A}'$ is always chosen as the preimage of the annulus
between $B_{n+1}$ and $B_{n'}$. If a maximal branch inside the
parent branch is replaced with another branch, the annuli
$\overline{A}_{i}$ with $i=1,2,3,4$ stay the same, as they are all
outside of the parent branch. So, we only need to show that
$\overline{A}'$ can be chosen larger than for the immediate preimage. 
But this is clear, since the branch defined on $\overline{B}$ has a univalent
extension onto $B_{n'}$, thus an annulus of the desired modulus will
always sit inside this extended domain. Then, there is an extra
annulus between the extended domain and the complement of the parent
branch, whose annulus we can usually bounded away from $0$.

 \paragraph{Immediate parent branches.}  
Choose a branch $\overline{B}$ of $\phi_{1}$ which is neither critical
nor immediate. Assume that the parent domain is an immediate preimage of $B_{n}$ by the
central branch. Since we assumed that $\overline{B}$ was not an
immediate preimage of $B_{n+1}$, we see that while annuli
$\overline{A}_{i}$ for $i=1,2,3,4$ can be chosen as for immediate
preimages, $\overline{A}'$ will be larger by an annulus between the 
complement of
the extended immediate domain and the rank $n$ extension of
$\overline{B}$. Since all branches of $\phi$ sit inside $B_{n}$ with
separating moduli at least $\frac{1}{2}\alpha$, all parent branches
are nested inside $B_{n}$ with moduli at least $\frac{1}{4}\alpha$. 
The extension of $\overline{B}$ will be mapped by the extension of the
immediate branch inside some parent branch. Thus, the extra
contribution to $\overline{A}'$ will be at least $\frac{1}{4}\alpha$.
That means $\overline{s}_{3}$ and $\overline{s}_{4}$ will both grow by
the this amount compared with the estimate for the immediate
preimages.   
Then, one can choose $\beta' = \frac{9}{8}\beta$ and 
\[ \overline{\lambda}_{1} = \frac{\lambda_{2}(B)}{2} -
\frac{\alpha}{8}\]
and 
\[\overline{\lambda}_{2} = \frac{\lambda_{1}(B)}{2} +
\frac{\alpha}{8}\; .\]
One checks directly that this gives a normalized symbol.

\paragraph{The case of independent $B$ and $B'$.}
Let us first consider the independent case. To pick $\overline{A}_{2}$, we
consider the annulus separating $B$ from the boundary of the domain of
its rank $(n-1)'$ extension. We claim that its modulus in all cases is
estimated from below by
$\alpha+\delta$ where $\delta$ can be chosen as the maximum of
$-\lambda_{2}(B)$ and $-\lambda_{2}(B')$.  Indeed, if $B$ is carried
onto $B_{n}$ by the extended branch, the estimate is $\alpha$ plus the
maximum of $\lambda_{1}(B)$ and $\lambda_{1}(B')$ which is at least
$\delta$ since $\lambda_{1} + \lambda_{2} \geq 0$ in any normalized
symbol. On the other hand, if $B$ is mapped by the extension onto
something different from the central hole, the estimate
$\beta-\lambda_{2}\geq \frac{3}{4}\beta$ applies which is better than 
$\alpha+\lambda_{1}\leq \frac{3}{4}\beta$. 
To pick $A_{1}$, consider the
annulus of modulus $\beta+\lambda_{2}(B)$ separating $B$ from the
boundary of $B_{n-1}$. Pull these annuli
back by the central branch to get $\overline{A}_{2}$ and
$\overline{A}_{1}$ respectively.
By the hypothesis of the induction, the estimates are
\[ \overline{s}_{1} = \frac{\beta+\lambda_{2}(B)}{2}\; \mbox{and}\]     
\[ \overline{s}_{2} = \frac{\alpha+\delta}{2} \; .\]

As always, $\overline{A}'$ is determined with modulus at least
$\overline{s}_{1}$. The annulus $\overline{A}_{3}$ will be obtained as the
preimage by the central branch of the annulus surrounding the preimage
of $B'$ in the domain of the extended branch. This has modulus at
least $\alpha+\delta$ in all cases as argued above. The annulus 
$\overline{A}_{4}$ is the preimage of the
annulus surrounding the extension in $B_{n-1}$. By induction,
\[ \overline{s}_{3} = \frac{\beta+\lambda_{2}(B)}{2} +
\alpha+\delta\; \mbox{and}\]
\[ \overline{s}_{4} = \overline{s}_{3} +
\frac{\beta + \lambda_{2}(B') - \alpha - \delta}{2} \; .\]

We put $\overline{\lambda}_{1} = \frac{\lambda_{2}(B)}{2}$ and
$\overline{\lambda}_{2}=\frac{\alpha-\delta}{2}$. We check that 
\[ \overline{s}_{3} + \overline{\lambda_{1}} = \frac{\beta}{2} + \alpha +
\lambda_{2}(B) + \delta \geq \beta - \lambda_{2}(B) +
\lambda_{2}(B) \geq \beta \; .\]

In a similar way one verifies that
\[ \overline{s}_{4} - \overline{\lambda}_{2} \geq \beta\; .\]  
Also, the required inequalities between corrections
$\overline{\lambda}$ follow directly. In this case, it is not evident
how to obtain a symbol with norm greater than $\beta$, so we will
return later to this case with more careful estimates. 

\paragraph{The case of $B'$ subordinate to $B$.}
Presently, we will consider the case of $B'$ being subordinate to $B$. This
means that the same univalent branch of rank $n-1$ transforms $B$ onto
$B_{n}$ and $B'$ onto some $B''$. Consider the annulus separating
$B_{n}$ from $B''$, and a larger annulus surrounding the previous one
in $B_{n-1}$. Their preimages first by the extended branch and then by
the central branch give us $\overline{A}_{2}$ and 
$\overline{A}_{1}$ respectively. Notice that if an annulus surrounds
the domain of the extended branch, its preimage can be used to get
another layer of $\overline{A}_{1}$.    
The estimates are
\[ \overline{s}_{1} = \frac{\alpha+\delta'}{2}\]     
where $\delta'$ is the
maximum of $\lambda_{1}(B'')$ and $0$. This is allowed, since if
$\lambda_{1}(B'')$ is negative, we can always use the fact that $B'$
is nested in $B_{n-1}$ with modulus at least
$\beta+\lambda_{2}(B')>\alpha$. 
\[ \overline{s}_{2} = \frac{\alpha-\lambda_{2}(B'')}{2} \; .\]
The annulus $\overline{A}'$ is uniquely determined with modulus
$\overline{s}_{1}$, and $\overline{A}_{3}$ will be the preimage of the
annulus separating $B''$ from $B_{n}$. Finally, $\overline{A}_{4}$
will separate the image of $\overline{A}_{3}$ from $B_{n-1}$. Again,
the separation between the extended domain and $B_{n-1}$ will give us
another layer of $\overline{A}_{4}$. The estimates are
\[ \overline{s}_{3} = \frac{\alpha+\lambda_{1}(B'')}{2} + \beta -
\lambda_{1}(B'') = \beta + \frac{\alpha-\lambda_{1}(B'')}{2} \; \mbox{and}\]
\[ \overline{s}_{4} = \overline{s}_{3} +
\frac{\lambda_{1}(B'')+\lambda_{2}(B'')}{2} = \beta +
\frac{\alpha+\lambda_{2}(B'')}{2} \; .\]

Set 
\[ \overline{\lambda}_{1} = \frac{-\alpha+\delta'}{2}\; \mbox{ and}\] 
\[ \overline{\lambda}_{2} = \frac{\alpha+\lambda_{2}(B'')}{2}\; .\]
The requirements of a normalized symbol are clearly
satisfied. Moreover, as we noted, an estimate of the separation
between the extended domain and $B_{n-1}$ will give better
$\overline{s}_{1}$, $\overline{s}_{3}$ and $\overline{s}_{4}$. If,
say, this extra modulus is 
at least $K\beta$, $K<1/2$, then $\overline{s}_{1}$,
$\overline{s}_{3}$  and $\overline{s}_{4}$ will grow by at
least $\frac{K}{2}\beta$. 

Since this situation will recur, we emphasize the following reasoning
as a lemma.
\begin{lem}\label{lem:31ka,1}
Suppose that a separation symbol $(s_{1},s_{2},s_{3},s_{4})$ 
can be represented by a normalized symbol with norm $\beta$ and
corrections $\lambda_{1}$ and $\lambda_{2}$. If a number $0<K<1/2$
exists such that 
\[ s_{1} \geq (1 + K) \alpha + \lambda_{1}\; ,\]
\[ s_{3} \geq (1 + \frac{K}{2})\beta - \lambda_{1}\; \mbox{ and}\] 
\[ s_{4} \geq (1 + K) \beta + \lambda_{2}\; ,\]
then another normalized
symbol can be built for $(s_{1},s_{2},s_{3},s_{4})$ with norm 
$(1 + \frac{K}{4})\beta$. 
\end{lem}
\begin{proof}
We will distinguish the new normalized symbol by writing primes. 
Set
\[ \lambda'_{1} = \lambda_{1} + K\frac{\alpha}{4}\; \mbox{ and}\]
\[ \lambda'_{2} = \lambda_{2} + K\frac{\alpha}{4}\; .\]
Since the corrections were increased the bounds from below on
$\lambda'_{i}$ and $\lambda'_{1} + \lambda'_{2}$ will remain in force. 
Also, the bounds from above on  will remain since $\alpha'$
will grow by the same amount as $\lambda'_{i}$.  
Now, the corrections $\lambda'_{i}$ were set so that 
\[ \alpha - \lambda_{2} = \alpha' - \lambda'_{2} \; ,\]
hence $s_{2}$ will be correctly estimated.  Finally, we
check directly that 
\[ \beta' - \lambda'_{1} = \beta - \lambda_{1} + K\frac{\beta}{8} <
s_{3}\; ,\]  
\[ \alpha' + \lambda'_{1} = \alpha + \lambda_{1}+ \frac{3}{4}K\alpha < s_{1}\; \mbox{
and}\]
\[ \beta + \lambda'_{2} = \beta + \lambda_{2} + \frac{3}{8}K\beta < s_{4}\]
which concludes the proof.
\end{proof}

Coming back to our situation, we see that a uniform
bound on the separation between the extended domain and $B_{n-1}$ will
give us a uniform increase of the norm in the case of $B'$ subordinate
to $B$ by Lemma~\ref{lem:31ka,1}.  This is provided by the following lemma: 
\begin{lem}\label{lem:30ka,1}
Under the hypotheses of Proposition~\ref{prop:21kp,1}, let $0<j\leq
k$. Let $B_{n}$ be the central hole of $\phi_{j}$. 
Choose $D$ to be the domain of an extended univalent
branch of $\phi_{j}$ of rank $n'$. Then there is a positive bound
$K(l)$ such that $D$ is nested inside $B_{n'}$ with modulus at least 
$K(l)\beta$. 
\end{lem}
\begin{proof}
Consider the previous inducing step on $\phi_{j-1}$. All extended
branches of $\phi_{j}$ are nested inside parent branches created after
the first filling of $\phi_{j-1}$. So, it is enough to prove the bound
for the parent branches. If $B'$ is a parent non-immediate branch,
consider its push-forward image. By separation bounds, the
push-forward image is surrounded inside $B_{(n-2)'}$ by an annulus of
modulus at least $\alpha$ which also separates it from the central
hole. The bound by $\alpha/2$ for the parent branch follows. Extended
immediate preimages are nested with annuli 
\[ 2^{-k-1} (\alpha + \lambda_{1}(B)) \geq 2^{-k-2}\alpha\]
where $k$ is the number of consecutive images of the critical point
inside the central hole. This gives a uniform bound in terms of $l$ of
Proposition~\ref{prop:21kp,1}.
\end{proof}

\paragraph{The case of $B$ subordinate to $B'$.}
This situation is
analogous to the situation of immediate preimages considered at the
beginning. Indeed, by mapping $B$ to $B_{n}$ and composing with the
central branch one can get a folding branch of rank $n-1$ defined on
$B$. We now see that the situation inside the domain of the rank $n-1$
extension of $B$ and $B'$ is analogous to the case of immediate
preimages, except that the folding branch maps onto a larger set
$B_{n-1}$. Like in the previous case, the estimates do not use the
separation between the extended domain and $B_{n-1}$. By Lemma~\ref{lem:30ka,1} this gives a definite improvement in terms of $\beta$. Hence 
$\overline{s}_{1}$, $\overline{s}_{3}$ and $\overline{s}_{4}$ all
improve, and we can
increase the norm of the symbol by Lemma~\ref{lem:31ka,1}.  

This concludes the proof of Lemma~\ref{lem:28np,1}.
\subsection{Close returns}
\paragraph{Topological description.}
\subparagraph{Notation.}
We will use the notation $\tilde{\varphi}$ for the type I mapping obtained from $\phi_{j}$ in the first stage of the box inducing step.
Also, the $i$-th preimage of
$B_{n}$ by the central branch will be called $B^{i}$, thus $B_{n} =
B^{0}$. We take $k$ to be the smallest $i$ such that the critical
value does not belong to $B^{i}$. Let $B^{*}$ mean the parent  domain  of
 the branch of $\tilde{\varphi}$ which contains the critical value of
 $\psi_{j}$.$B^{*}$ is the $k$-th preimage of some $B$ by the central branch,
 where $B$ is a univalent domain of $\phi_{j}$.  

\subparagraph{Classification of branches.}
According to the definition of the box inducing step, $\phi_{j+1}$ is obtained from $\tilde{\phi}$ by a box  inducing step defined for non-close returns. 
The first filling of $B^{k}$
(which becomes $B_{n+1/2}$) gives parent branches.\footnote{One must avoid confusion between parent branches of $\tilde{\phi}$ and parent branches of $\phi_{j+1}$.}
 Those will be of
two types: immediate preimages of $B^{k}$ by the central branch and
preimages of branches of rank $n+1/2$. We will split preimages of
$B_{n+1}$ into two classes. First, we consider $D$-type preimages.
Those can be mapped by a composition of  immediate branches only until they hit $B_{n+1}$. All other domains of $\phi_{j+1}$ will be called
$\delta$-type preimages. It follows that any $\delta$-type preimage
falls into a non-immediate parent branch before hitting $B_{n+1}$. A
reasoning similar to one conducted in the case of non-close returns shows that separation estimates for $\delta$-type
preimages will be the worst for  domains which are mapped by immediate
 branches only until they hit a {\em maximal} domain
  inside a non-immediate parent branch.
 Thus, we only consider this kind of $\delta$-preimages. 

The outcome of our computations will be as follows:
\begin{lem}\label{lem:28np,2}
Suppose that a sequence of type I complex box mappings $\phi_{j}$ is given which satisfies the assumptions of
 Proposition~\ref{prop:21kp,1}. Let $0\leq j \leq k$ and assume that $\phi_{j}$ shows a close return whereby the critical value stays in the central domain for $k$ steps. If $\beta$ is the separation norm of $\phi_{0}$, then normalized critical symbols 
\begin{itemize}
\item
for  immediate preimages of $B_{n+1}$ the symbol has norm $\beta$ with corrections
\[ \overline{\lambda}_{1} = \frac{\alpha}{2} - \frac{\alpha-\lambda_{2}(B)}{2^{k+1}}\: ,\: 
   \overline{\lambda}_{2} = -\frac{\alpha}{2} + \frac{\alpha+\lambda_{1}(B)}{2^{k+1}}\; . \]
\item
for other $D$-type preimages the same symbol will work as for immediate preimages, or a better symbol with norm $(1+K(k))\beta$ can be used, where $K$ is a positive function of $\beta$, but tends to $0$ with $k$ growing to infinity. 
\item
in all other situations the norm is at least $(1+K(l))\beta$ with $K$ a positive function of $l$ only ($l$ is defined in the statement of Proposition~\ref{prop:21kp,1}.)
\end{itemize}
\end{lem}

\paragraph{Immediate preimages.}
We first consider the annulus $C$ surrounding $B_{k}$ inside $B_{n}$.
Let us denote its modulus by $\gamma$. Clearly, 
\begin{equation}\label{equ:28ka,1}
\gamma \geq (\alpha+\lambda_{1}(B))(1-2^{-k})\; .
\end{equation}
 To construct $\overline{A}_{2}$, we
first consider the annulus surrounding $B^{*}$ inside $B^{k-1}$ and
separating $B^{*}$ from $B^{k}$. This annulus is obtained as the
preimage of $A_{3}\cup A'$ which existed for $b$ by a univalent map.
Then, inside $B^{*}$ the preimage of $B^{k}$ nests with modulus
$\gamma$. Thus,
\[\overline{s}_{2} = \frac{\beta - \lambda_{1}(B) + \gamma}{2}\; .\]   

Then, $\overline{A}_{1}$ can be
obtained as the preimage by the central branch of the annulus which 
separates $B^{*}$ from the
boundary of $B^{k-1}$ and is the $k$-th preimage of $A_{4}$. This gives
\[ \overline{s}_{1} = \frac{\beta-\lambda_{1}(B)+\gamma}{2} + 
\frac{\lambda_{1}(B)+\lambda_{2}(B)}{2^{k+1}}\; .\]

The annulus $\overline{A}'$ is defined in the natural way as the
preimage of the annulus surrounding $B_{n+1}$ inside $B^{k}$. We have
already estimated its modulus by $\overline{s}_{1}$.

The modulus of $\overline{A}_{3}$ is easily estimated by 
\[ \frac{\alpha-\lambda_{2}(B)}{2^{k}}\; .\]

The annulus $\overline{A}_{4}$ will be the preimage of the the skinny
annulus encircling
$B^{*}$ in $B^{k-1}$. Its modulus is bounded from below by
\[ \frac{\lambda_{1}(B)+\lambda_{2}(B)}{2^{k+1}}\; .\]

This gives
\[ \overline{s}_{3} = \frac{\alpha-\lambda_{2}(B)}{2^{k}} +
\frac{\beta-\lambda_{1}(B)+\gamma}{2} + \frac{\lambda_{1}(B)+
\lambda_{2}(B)}{2^{k+1}} 
\; ,\]
\[ \overline{s}_{4} = \overline{s}_{3} +  \frac{\lambda_{1}(B)+
\lambda_{2}(B)}{2^{k+1}}  =\]
\[ \frac{\beta-\lambda_{1}(B)+\gamma}{2} + \frac{\alpha +
\lambda_{1}(B)}{2^{k}}  \; .\] 

We now observe that $\gamma$ occurs in all estimates with the positive
sign, so we can replace it with the lower estimate~[\ref{equ:28ka,1}].

This means the choice of 
\[ \overline{\lambda}_{1} = \frac{\alpha}{2} - \frac{\alpha -
\lambda_{2}(B)}{2^{k+1}} \;\mbox{ and }\]
\[ \overline{\lambda}_{2} = -\frac{\alpha}{2} +
\frac{\alpha+\lambda_{1}(B)}{2^{k+1}}\; .  \]

>From this one can see directly that the requirements of a normalized
symbol are satisfied. 

\paragraph{$\delta$-type preimages.}
Choose a $\delta$-type preimage $\overline{B}$. It is mapped by a composition of immediate parent branches into a non-immediate parent domain which we call $B_{p}$. We argued that it suffices to consider the case when $\overline{B}$ gets mapped onto a m
Now, consider the parent branch of  $\psi_{j}(B_{p})$ (as a branch of $\tilde{\phi}$) and
call it $\tilde{B}$. We distinguish
between two cases depending on whether $\tilde{B}$ and $B^{*}$ are the same or not.   

\subparagraph{$\tilde{B}$ is distinct from $B^{*}$.}
Let $B'$ be $\psi^{k}(\tilde{B})$ and $B$ mean the $\psi^{k}(B^{*})$.
Let $C$ be the modulus on an annulus separating 
$B$ from $B_{n}$, $B'$ and the complement of $B_{n-1}$. Analogously,
let $C'$ separate $B'$ from $B$, $B_{n}$ and the complement of
$B_{n-1}$. The triple alternative between one independent and two
subordinate cases tells us that we can always choose $C+C' = \alpha$. 

To get the symbol, we use 
\[ \overline{s}_{1} = \frac{\beta - \lambda_{1}(B) + \gamma}{2}\]
obtained for immediate preimages (we skipped a positive term.) The
reader is reminded that $\gamma$
is a lower estimate on the modulus between $B^{k}$ and $B_{n}$, thus
\[ \gamma \geq (\alpha + \lambda_{1})(1-2^{-k}) \; \]
where $\lambda_{1}$ is the supremum of $\lambda_{1}(D)$ over all
monotone holes of $\phi$. 

We put 
\[ \overline{s}_{2} = \frac{\gamma+C}{2}\; .\] 
The annulus $\overline{A}'$ has modulus 
$\overline{s}_{1}$. $\overline{A}_{3}$ has two layers: one is the
preimage of the annulus between $B^{k}$ and $B_{n}$ inside the parent
branch, another is the preimage by the central branch of an annulus
separating the parent branch from $B^{k}$ inside $B_{k-1}$. Those give
\[ \overline{s}_{3} = \frac{\beta - \lambda_{1}(B)}{2} +
\frac{3}{2}\gamma + \frac{\beta-\lambda_{1}(B')}{2}+\frac{C'}{2}\]
Finally, we
put $\overline{s}_{4} = \overline{s}_{3}$.  

The only way $k$ enters the estimates is through $\gamma$. Thus, we
only consider $k=1$ which gives the smallest value of $\gamma$, namely
$\gamma=\frac{\alpha+\lambda_{1}}{2}$. 
Then, pick $\beta' = \frac{9}{8}\beta$. Take 
\[\overline{\lambda}_{1} = - \frac{\alpha}{8}\; .\]
Then 
\[ \overline{s}_{1} \geq \alpha +
\frac{\alpha+\lambda_{1}-2\lambda_{1}(B)}{4}  \geq 
\alpha = \frac{9}{8}\alpha - \overline{\lambda}_{1}\; \]
while
\[ \overline{s}_{3} \geq \beta + 2\gamma - \frac{\lambda_{1}(B) +
\lambda_{1}(B')}{2} - \frac{\gamma}{2} + \geq \beta + \frac{\alpha}{2} =
\frac{5}{4}\beta \geq \beta' - \overline{\lambda}_{1}\; .\]

Next, put 
\[ \overline{\lambda}_{2} = \frac{-C - \gamma+ 2\alpha}{2} +
\frac{1}{8}\alpha\; \mbox{ if $\overline{s}_{2}<\alpha$ and}\]
\[ \overline{\lambda}_{2} = \frac{\alpha}{8}\; \mbox{ otherwise.}\]

Thus guarantees that $\overline{s}_{2}$ is suitably bounded by 
$\beta'/2 - \overline{\lambda}_{2}$. 
We have already seen that 
\[ \overline{s}_{4} \geq \frac{5}{4}\beta \geq \beta' +
\frac{\alpha}{8}\; .\]

So we only need to check the case when $\overline{\lambda}_{2}$ is
given by the first more complicated formula.  
\[ \overline{s}_{4} - \overline{\lambda}_{2} = \beta -
\frac{\lambda_{1}(B)+\lambda_{1}(B')}{2} + 2\gamma + \frac{C+C'}{2} -
\frac{9}{8}\alpha =\]
\[ \beta - \frac{\alpha}{8} + \frac{C+C'}{2} +
\frac{2\lambda_{1} - \lambda_{1}(B) - \lambda_{1}(B')}{2} \geq 
\beta + \frac{\alpha}{2} - \frac{\alpha}{8} \geq \beta' \; .\]

Also, other inequalities between $\overline{\lambda}_{i}$ are easily
verified.           

\subparagraph{$\tilde{B}$ is equal to $B^{*}$.}
In this case, we try mapping forward by parent branches of $\tilde{\phi}$ until
images of $\psi(B_{p})$ and the postcritical branch of $\tilde{\phi}$
land in different domains of $\phi_{k}$. Then, we distinguish between
the situation in which both these domains are monotone, or one of them
is central. If both are monotone, the same estimates apply as in the
previous case. Indeed, the two distinct branches obtained here can be
used as $\tilde{B}$ and $B^{*}$ in the previous case to obtain the
same separation bounds.   

We are then left with two subordinate cases. First, assume that the
postcritical branch gets mapped onto $B^{k}$, and $\psi(B_{p})$ lands
in another domain of $\phi_{k}$, call it $B'_{p}$. Then, we estimate 
\[ \overline{s}_{1} = \frac{\beta}{2}\]
\[ \overline{s}_{2} = 0\; .\]
These are both weak estimates that hold in all cases. The annulus
$\overline{A}'$ as always will have modulus at least
$\overline{s}_{1}$. The annulus $\overline{A}_{3}$ will have two
layers. The inner one will be of modulus at least $\gamma$ and is the
preimage of the annulus which
surrounds the image of $\psi(B_{p})$ inside $B'_{p}$. The outer one
has modulus $\beta-\lambda_{1}$ and is the preimage of the annulus
that separates $B'_{p}$ from $B^{k}$ inside $B_{n}$. The worst case
is when $k=1$ which gives
\[ s_{3} \geq \frac{3}{2}\beta \; .\]
As $\overline{s}_{4}$ we take the preimage of the annulus surrounding
$\tilde{B}$ inside $B^{k-1}$ by the central branch. That will add
modulus at least $\alpha/2$. So, we can take 
\[ s_{4} = \frac{7}{4}\beta \; .\]

Clearly, a normalized symbol with 
$\beta' = \frac{9}{8}\beta$, $\lambda_{1} = -\frac{\beta}{16}$ and 
$\lambda_{2} = \frac{\beta'}{2}$ will work. 

Finally, let $\psi(B_{p})$ be mapped on $B^{k}$ with the postcritical
domain going into some $B^{p*}$. 
Let $l$ be the
least positive integer such that $\psi^{l}(B^{p*})$ is not in $B_{n}$. 
Certainly, $l\leq k$. Then, $\overline{A}_{2}$ has two
layers. One is the preimage of the annulus surrounding the image of
the postcritical branch inside $B^{p*}$, and that has modulus at least
$\gamma$. Another is obtained as the preimage of the annulus
separating $B^{p*}$ from $B^{k}$ inside $B^{k-1}$, that gives modulus
at least $\frac{\beta-\lambda_{1}(B^{t*})}{2}$ where 
$B^{t*} = \psi^{l}(B^{p*})$.  Thus,
\[ \overline{s}_{2} = \frac{\beta - \lambda_{1}(B^{t*}) + \gamma}{2}
\; .\]
 The annulus $\overline{A}_{1}$ will have three layers.
The inner one
will be the preimage of the annulus surrounding $A_{3}(B^{p*})$ inside
$B_{l-1}$. Its modulus is
\[ \frac{\lambda_{1}(B^{t*})+\lambda_{2}(B^{t*})}{2^{l+1}}\]
Next, we have the preimage of an
annulus between
$B^{l-1}$ and $B_{n}$. This is at least
$(\alpha+\lambda_{1})(1-2^{1-l})$. Finally, we have the preimage of an
annulus separating $B^{k}$ from $B^{*}$ inside $B^{k-1}$, that is at
least $\alpha/2$. All this gives
\[ \overline{s}_{1} = \frac{1}{4}\beta +
(\alpha+\lambda_{1})(1-2^{1-l}) + \frac{\lambda_{1}(B^{t*}) +
\lambda_{2}(B^{t*})}{2^{l+1}} + \overline{s}_{2}\; .\]
One easily sees that this has the least value for $l=1$, namely
\[ \overline{s}_{1} \geq \frac{\beta - \lambda_{1}(B^{t*})}{2} +
\frac{\alpha+\lambda_{1}}{4} + \frac{\lambda_{1}(B^{t*}) +
\lambda_{2}(B^{t*})}{4} + \frac{1}{4}\beta \geq\]
\[\geq \frac{3}{4}\beta + 
\frac{\alpha+\lambda_{2}(B^{t*})}{4} \geq \frac{13}{16}\beta\; .\]    
Thus, we take 
\[ \overline{s}_{1} = \frac{13}{16}\beta\]
\[ \overline{s}_{2} = \frac{\beta}{2}\; .\]
Then, we trivially have 
\[ \overline{s}_{3} = \overline{s}_{1} = \frac{13}{16}\beta\]
and 
\[ \overline{s}_{4} = \overline{s}_{3} + \frac{\alpha}{2} =
\frac{17}{16}\beta\; .\]
The additional term $\alpha/2$ in $\overline{s}_{4}$ comes from
$\overline{A}_{4}$ which is taken at least as the preimage by the
central branch of the annulus separating $B^{*}$ from $B^{k}$ inside
$B^{k-1}$. If $\overline{s}_{i}$, $i$ not $3$ were all decreased by a
sixteenth, a normalized symbol with norm $\beta$ could obviously be
set up. Thus, we get a symbol with greater norm by
Lemma~\ref{lem:31ka,1}. 

This ends the analysis of $\delta$-type preimages.       
\paragraph{$D$-type preimages.}
 Choose a non-immediate $D$-type preimage and call it $B'$.
We go
back to the estimates done previously for immediate preimages in the
case of a close return. The bound $\overline{s}_{1}$ certainly
remains, moreover, we observe that the annulus $\overline{A}_{2}$ for
immediate preimages was chosen in such a way that it also separates
their extensions rank $B_{n+1/2}$ from $B_{n+1}$, Thus,
$\overline{s}_{2}$ also remains in force for all $D$-type preimages. 
For non-immediate $D$-type preimages, the annulus $\overline{A}'$ can be chosen
with modulus equal to $\overline{s}_{4}$ in the immediate case.
Indeed, $B'$ will be mapped onto an immediate preimage by a mapping which
extends in a univalent way onto $B^{k}$. Then $\overline{A}_{3}$ is       
the annulus surrounding the domain of the rank $n+1/2$ extension of
$B'$ inside the parent branch (which is one of the extended immediate
preimages.) Its modulus is at least as large as the modulus between an
extended immediate preimage of $B^{k}$ inside $B_{k}$, and that is
half of the modulus of $B^{k}$ inside $B^{k-1}$, or at least 
\[ (\alpha+\lambda_{1}(B))\cdot 2^{-k} \geq \beta\cdot 2^{-k-2}\; .\]
So, $\overline{s}_{3}$ will grow by a half of this amount compared to the
estimate for immediate preimages. Since $\overline{A}_{4}$ remains as
in the immediate case, $\overline{s}_{4}$ will also grow by the same
correction. For simplicity, denote this correction by $K\beta$ with
$K$ uniformly bounded away from $0$ in terms of $k$ and no more than
$1/2$. Let 
\[\beta'=(1+\frac{K}{4})\beta\; .\]
Now let us mark estimates from the immediate case with bars, and the
estimates from the present case with tildes. Then put
\[ \tilde{\lambda}_{1} = \overline{\lambda}_{1} - \frac{K}{8}\beta\]
\[ \tilde{\lambda}_{2} = \overline{\lambda}_{2} + \frac{K}{8}\beta\;
.\]

A comparison with the bounds for immediate preimages shows that this
gives a legitimate normalized symbol. Thus, for $D$-type preimages the
norm of the symbol can be increased by an amount bounded away from $0$
in terms of $k$.  

This concludes the proof of Lemma~\ref{lem:28np,2}. 
\subsection{Proof of Theorem B.}
\paragraph{Statement of results.}
The estimates we have done so far show that in all cases the symbols
after inducing can be chosen with at least the same norm. This proves
the first part of Proposition~\ref{prop:21kp,1}, namely that the norms
of critical symbols stay at least $\beta$. We will now work to get the
improved claim.

\begin{lem}\label{lem:30ka,2}
In the situation of Proposition~\ref{prop:21kp,1}, let $2 < j \leq k$.
Consider a univalent hole $\overline{B}$ of $\phi_{j}$ other than an immediate
preimage of the central hole. Then, the critical symbol of
$\overline{B}$ can be chosen
with norm $\beta'\geq (1+K(l))\beta$ with $K(l)$ always positive and
depending only on $l$.
\end{lem}
\begin{proof}
Until we say otherwise, our reasoning will be valid for $j=2$ as well.
A large portion of the proof is already in Lemmas~\ref{lem:28np,1} and~\ref{lem:28np,2}.
We must must only return to the "independent" case when $\phi_{j-1}$ shows no close return while
the $B'$ (defined as the image by $\psi$ of the parent branch of $\overline{B}$) and the
postcritical branch $B$ are independent. 

Let $B_{n}$ be the central hole of $\phi_{j-1}$, $B$ the postcritical
branch, and $B'$ the push-forward of the
parent branch of $\overline{B}$. The extended branch of $B$ maps the domain
of $B$ onto some $D$, and similarly $B'$ is mapped onto $D'$ by its
extension. By Lemma~\ref{lem:30ka,1}, we know that the domains of
the extended branches are surrounded inside $B_{(n-1)'}$ by     
an annulus with modulus at least $K_{1}(l)\beta$. Now we consider a
number of cases depending on whether $D$ and $D'$ are central or not.

First, assume that $D$ is not central. Then, the same estimates we
got in our previous analysis of the independent case will be valid with
$B$ replaced with $D$. However, $\overline{A}_{1}$ will contain an
extra layer of modulus at least $\frac{\beta}{2}K_{1}(l)$. Thus,
$\overline{s}_{1}$ will be improved. This will also cause
$\overline{s}_{3}$ and $\overline{s}_{4}$ to improve. By
Lemma~\ref{lem:31ka,1}, the norm $\beta$ can be increased by 
$\frac{\beta}{8}K_{1}(l)$. 

Knowing that $D$ is central, we have an option of estimating 
\[ \overline{s}_{1} \geq K_{1}(l)\beta + \frac{\alpha+\lambda_{1}}{2}\; 
\mbox{ and}\]
\[ \overline{s}_{2} \geq \frac{\alpha+\lambda_{1}}{2}\;\]
where $\lambda_{1}$ is the supremum of $\lambda_{1}(b)$ over all
univalent holes of $\phi_{j-1}$. 

Now suppose that $D'$ is not central. In this case, use the new
estimate for $\overline{s}_{1}$ and the old estimate
$\overline{s}_{2}$. 
Also, we take
\[ \overline{s}_{3} = \overline{s}_{1} + \beta + \lambda_{2}(D')\; ,\]
\[ \overline{s}_{4} = \overline{s}_{3} + K_{1}(l)\frac{\beta}{2}\; .\]

Then put $\beta' = \beta + \frac{\beta}{4}K_{1}(l)$
where $K_{1}(l)$ does not exceed $1/2$. 
Also, choose
\[ \overline{\lambda}_{1} = \frac{3K_{1}(l)\alpha}{4} -
\frac{\alpha-\lambda_{1}}{2}\; \mbox{ and}\]
\[ \overline{\lambda}_{2} =
\frac{\alpha-\delta}{2}+\frac{K_{1}(l)}{4}\alpha
\; .\]
     
We easily check that $\overline{\lambda}_{1} +
\overline{\lambda}_{2}>0$,
also it is clear that $\overline{s}_{1}$ and $\overline{s}_{2}$ are
bounded as desired. Next, check
\[ \overline{s}_{3} + \overline{\lambda}_{1} =
\frac{7}{4}K_{1}(l)\alpha + \lambda_{1} + \beta + \lambda_{2}(D') >
\beta'\]
 and
\[ \overline{s}_{4} - \overline{\lambda}_{2} =
\frac{\alpha+\lambda_{1}}{2} + \beta + \lambda_{2}(D')
\frac{\alpha-\delta}{2} + \frac{7}{4}K_{1}(l)\alpha >\]
\[ \beta' + \frac{\lambda_{1} + \lambda_{2}(D')}{2} +
\frac{\delta+\lambda_{2}(D')}{2} \; .\]
Since both fractions are non-negative, we are done.   

Thus, we are left with the case when both $D$ and $D'$ are central. 
We first consider the situation of $\delta$ very large, i.e.
\begin{equation}\label{equ:30kp,1}
 \delta \geq \frac{\alpha}{2} - \frac{1}{4}K_{1}(l)\alpha\; .
\end{equation}
As $\lambda_{1}\\delta$, we use 
\[ \overline{s}_{1} = \frac{\alpha+\delta}{2} + K_{1}(l)\alpha\; ,\]
\[ \overline{s}_{2} = \frac{\alpha+\delta}{2}\; ,\]
\[ \overline{s}_{3} = \overline{s}_{1} + \alpha + \delta\; ,\]
\[ \overline{s}_{4} = \overline{s}_{3}\; .\]

Then, pick 
\[ \beta' = (1 + \frac{K_{1}(l)}{8})\beta\; ,\]
\[\overline{\lambda}_{1} = \frac{\delta-\alpha}{2} +
\frac{7}{8}K_{1}(l)\alpha\; ,\]
\[\overline{\lambda}_{2} = \frac{\alpha-\delta}{2} +
\frac{1}{8}K_{1}(l)\alpha\; .\]
Provided that $K_{1}(l)\leq 1/2$ we have no trouble seeing that the
inequalities between corrections are satisfied. We need to check
\[ \overline{s}_{3} + \overline{\lambda}_{1} = 2\delta +
\alpha + \frac{15}{8}K_{1}(l)\alpha  \geq 2\alpha + \frac{11}{16}\beta
> \beta'\] and
\[ \overline{s}_{4} - \overline{\lambda}_{2} = 2\delta + \alpha + 
\frac{7}{8}K_{1}(l)\alpha = 2\alpha + \frac{3}{16}\beta > \beta'\; .\]

So, we are done in this case and can assume the converse of
inequality~[\ref{equ:30kp,1}].  We proceed to note that for $j>1$
there is a positive bound $K_{2}(l)$ such that 
\begin{equation}\label{equ:30kp,2}
 \lambda_{1}(B) + \lambda_{2}(B) \geq K_{2}(l)\beta
\end{equation}
for all univalent holes $B$ of $\phi_{j}$.   
In the cases already dealt with in the proof of
Lemma~\ref{lem:30ka,2} this estimate follows. Indeed, once the norm of
the symbol is larger than $\beta$, one can always decrease the norm and at the
same time increase the corrections $\lambda_{i}$ by the same amount
getting weaker estimates. If a correction exceeds $\alpha$ as a
result, it means that it was close to $\alpha$ to begin with, and thus
the sum of the corrections was already positive. So we only need
to prove estimate~[\ref{equ:30kp,2}] in the remaining case. However,
then our original choice of 
\[ \overline{\lambda}_{1} = \frac{\lambda_{2}(B)}{2} \; \mbox{ and}\]
\[ \overline{\lambda}_{2} = \frac{\alpha - \delta}{2} \]
proves correct in view of the converse of
inequality~[\ref{equ:30kp,1}]. Thus, for $j\geq 3$ we can proceed
assuming that estimate~[\ref{equ:30kp,2}] holds for all holes of
$\phi_{j-1}$.  

In this situation, we set
\[ \overline{s}_{1} = \frac{\beta + \lambda_{2}(B)}{2}\; ,\]
\[ \overline{s}_{2} = \frac{\alpha + \delta}{2}\; ,\]
\[ \overline{s}_{3} = \overline{s}_{1} + \alpha + \lambda_{1}\; ,\]
\[ \overline{s}_{4} = \overline{s}_{3} +
\frac{\beta+\lambda_{2}(B')-\alpha-\lambda_{1}}{2}\; .\]

Let $\beta' = (1 + \frac{1}{4})K_{2}(l)\beta$.  
Put
\[ \overline{\lambda}_{1} = \frac{\lambda_{2}(B)}{2} - 
\frac{1}{4}K_{2}(l)\alpha\; \mbox{ and}\]
\[ \overline{\lambda}_{2} = \frac{\alpha-\delta}{2} +
\frac{1}{4}K_{2}(l)\alpha\; .\]

We see that $\overline{\lambda}_{1} + \overline{\lambda}_{2} > 0$ and
the bounds on $\overline{\lambda}_{i}$ also hold provided that $K_{2}<1/2$
which we can be assumed. Also, $\overline{s}_{1}$ and
$\overline{s}_{2}$ are correctly estimated. We check
\[ \overline{s}_{3} + \overline{\lambda}_{1} = \beta + \lambda_{1} +
\lambda_{2}(B) - \frac{1}{8}K_{2}(l)\beta\; .\]

By~[\ref{equ:30kp,2}], $\lambda_{1} + \lambda_{2}(B) \geq
K_{2}(l)\beta$ for any $B$, thus we get 
\[ \overline{s}_{3} + \overline{\lambda}_{1} \geq \beta +
\frac{7}{8}K_{2}(l)\beta > \beta'\; .\]

Also,
\[ \overline{s}_{4} - \overline{\lambda}_{2} = \beta +
\frac{\lambda_{2}(B)+\lambda_{2}(B')} + \frac{\lambda_{1} + \delta}{2}
- \frac{1}{8}K_{2}(l)\beta \geq \beta +\frac{3}{8}K_{2}(l)\beta >
\beta' \; .\]

This concludes the proof of Lemma~\ref{lem:30ka,2}. 
\end{proof} 

\paragraph{Proof of Proposition~\ref{prop:21kp,1}.}
The part of Proposition~\ref{prop:21kp,1} which says that norms of
critical symbols do not deteriorate in the inducing process follows
from the first round of our estimates. To prove that the norm
improves, assume that $j\geq 3$ and the critical value falls into a 
non-immediate preimage $B$. By Lemma~\ref{lem:30ka,2} we see two
things. First, there is a normalized critical symbol of $B$ with norm  
at least $(1+K(l))\beta$. Secondly, for all holes of
$\phi_{j+1}$ other than the immediate ones normalized critical symbols
can be constructed with increased norms. Thus, we are left to show
that the norm will also increase for critical symbols of the immediate
holes of $\phi_{j+1}$. This norm is at least equal to the norm of the
critical symbol of $B$ for immediate preimages, so we are done.    

\paragraph{Theorem B.}
Suppose that a mapping $\phi$ satisfies the assumptions of Theorem B.
Our first remark is this:

\begin{lem}\label{lem:8ma,1}
Under the hypotheses of theorem B, there is a function $K_{1}(\beta)$
bounded on any closed set in $(0,\infty)$ such that if in the sequence 
$\phi_{i}$ there is a rotation-like sequence 
\[\phi_{j},\phi_{j+1}, \cdots, \phi_{j+K_{1}(\beta)}\; ,\] 
then $\phi_{j+K_{1}}(\beta)$ as a real mapping satisfies the starting
condition.  
\end{lem}
\begin{proof}
This follows directly from Theorem A. We only need to check that if
$B_{n}$ is the central domain of $\phi_{j}$, then $|B_{n}|/|B_{n'}|$
is bounded away from $1$ uniformly in terms of $\beta$. However,
Proposition~\ref{prop:21kp,1} implies that as complex box mappings,
the central hole $B_{n}$ is nested inside $B_{n'}$ with modulus at
least $\beta/2$. The bound we need follows from classical analysis,
see~\cite{lehvi}. 
\end{proof}

By Corollary 5.2, we can assume that numbers $l_{j}$ are uniformly
bounded, or the starting condition is immediately satisfied. Thus, we
can regard the parameter $l$ in Proposition~\ref{prop:21kp,1} as an
absolute constant. 
Now let us fix some $j$. We observe that either $j+3+K_{1}(\beta)$
satisfies the starting condition, or there is a $j'\geq j+3$ such that   
$\phi_{j'}$ does not make a rotation-like return. But then it means
that for any $j$ either $\phi_{j+3+K_{1}}(\beta)$ already satisfies the
starting condition, or its separation norm is at least
$(1+K_{2}(l))\beta$ from Proposition~\ref{prop:21kp,1}. 

Next, we infer that for any $k$, $\phi_{k(j+3+K_{1}(\beta)}$ either
already satisfies the starting condition, or its separation norm is
at least $\beta(1+K_{2}(l))^{k}$. But if the separation norm becomes
sufficiently large, the starting condition follows. Thus, Theorem B
has been demonstrated.    

\subsection{Proof of Theorem D}
\paragraph{Superadditivity of conformal moduli.}
In the course of our complex induction we have repeatedly encountered nesting
annuli.  The modulus of their union was always estimated from below by the sum
of moduli. However, we will now show that unless the curve separating the two
nesting annuli is quite smooth, a definite increment can be added to the
modulus of the union. The following is a classical result:
\begin{fact}\label{fa:29np,1}
Let $A_{1}$ and $A_{2}$ be two disjoint open annuli situated so that $A_{1}$ 
separates $0$ from $A_{2}$ while $A_{2}$ separates $A_{1}$ from 
$\infty$. Assume further that both are contained in the ring 
$A=\{z: r <|z| < R\}$ for some $0<r<R$. 
By $C$ denote the set (annulus) of all points from $A\setminus(A_{1}\cup
A_{2})$ separated from $0$ and $\infty$ by $A_{1}\cup A_{2}$. Then, for 
every $\delta>0$ there is a number $\epsilon$ with the following property:
if  
\[ \mod A_{1} + \mod A_{2} \geq \mod A - \delta\; ,\]
then a $\rho$ exists for which the ring 
\[ \{z: \rho < |z| < (1+\epsilon)\rho\}\]
contains $C$. 
\end{fact} 

Fact~\ref{fa:29np,1} follows directly from a ``Modulsatz" of ~\cite{szkop}.

\subparagraph{Corollary.}
Our next lemma is a simple corollary to Fact~\ref{fa:29np,1}:

\begin{lem}\label{lem:26na,1} 
Let $\alpha$, $w$ and $\beta$ be non-intersecting Jordan curves, all separating
$0$ from $\infty$. Suppose that $w$ also separates $\alpha$ from $\beta$ and
that $w$ passes through $1$ and $-1$ in the complex plane. Let be the
annulus $A$ between $\alpha$ and $w$ and $B$  the annulus bounded
by $w$ and $\beta$. Suppose further that the moduli of $A$ and $B$ are both 
at least $\zeta$. Then, there is a
number $K\geq 1$ with the property that for each $\delta>0$ at least one of
these possibilities occurs:
\begin{itemize}
\item
the curve $w$ is in the Hausdorff distance less than $\delta$ from a
 $K(\zeta)$-quasicircle passing through $-1$ and $1$, or
\item
\[ \mod A + \mod B < \mod (A\cup B) - \epsilon(\delta)\]
where $\epsilon(\delta) > 0$ depends on $\delta$ only.
\end{itemize}
\end{lem}
\begin{proof}
Apply the uniformizing map which carries $\alpha$ onto a circle with radius
$r$, while $\beta$ goes to a circle with radius $R>r$ and $w$ is mapped onto a
Jordan curve $v$ in between so that $v$ passes through $1$. 
By Fact~\ref{fa:29np,1} unless the second part 
of the alternative holds, $v$ is contained in a neighborhood of the unit circle
of width $\delta$. Since the annuli of $A$ and $B$ are at least $\zeta$, for
$\delta<1/2$ the set of all inverses to uniformizing maps is normal on the ring
$\{z: 1 < |z| < 1+\delta\}$. From this, it first follows that the preimage of
the unit circle is a uniform quasicircle, and secondly that the preimage of
this ring is a narrow neighborhood of this quasicircle with width going to $0$ 
uniformly with $\delta$. So, the first part of the alternative holds.
\end{proof}

\paragraph{Proof of Theorem D.}
\begin{lem}\label{lem:26np,1}
Consider a type I complex box mapping $\phi$ of rank $n$. Suppose that $\phi$
shows a close return which results in a mapping $\phi_{1}$.  Let $v$ be the
ratio $|B_{n}|/|B_{n'}|$. Let $\beta$ be the separation norm of
$\phi$. There exists a function $V(\beta)>0$ so that
if $v < V(\beta)$ then at
least one possibility occurs:
\begin{itemize}
\item
the separation norm of $\phi$  is at least  $-K(\beta)\log v$ where $K$ is a
positive function of $\beta$ alone, or
\item
the separation norm of $\phi_{1}$ at least $\beta + \epsilon$ where 
$\epsilon>0$
 is an absolute constant. 
\end{itemize}
\end{lem}
\begin{proof}
By Lemma~\ref{lem:28np,2} we only need to consider $D$-type preimages,
including immediate ones. Also, by the same Lemma the symbol for immediate
preimages is also good for all $D$-type preimages, so it is enough to show that
unless the first part of the alternative holds, the second possibility must 
occur for immediate preimages.   
Let us track again the proof of Lemma~\ref{lem:28np,2} in the case immediate
preimages using the same notations. The annulus $\overline{A}_{2}$ is 
obtained as the preimage by the
central branch $\psi$ of a certain annulus $\Gamma$ composed of two nesting
annuli. The outer one surrounds the branch $B^{*}$ inside $B^{k}$ 
while the inner one is isomorphic to $C$ and surrounds the postcritical domain 
of $\tilde{\varphi}$ inside its parent domain $B^{*}$. Both annuli can be
mapped by $\psi^{k}$ which is univalent on $\Gamma$. Then, $B^{*}$ gets
mapped  onto a domain $B$ of the original mapping $\phi$. The outer layer 
of $\Gamma$ gets mapped onto $A_{3}(B)\cup A'(B)$. The estimate used in
Lemma~\ref{lem:28np,2} gives the modulus of $\Gamma$ as the sum of $s_{3}(B)$
and the modulus of $C$. However, we can get a positive correction here due to
the nesting of the image of $C$ inside $A'(B)$. First, map both by the
extension of the branch defined on $B$, which is univalent. Then $A'$ goes onto
the annulus surrounding $B_{n}$ inside $B_{n'}$. As the moduli of $C$ and $A'$ 
can both be bounded away from $0$ in terms of $\beta$, for a suitably
chosen $\beta$, by Lemma~\ref{lem:26na,1} 
the
nesting of the image of $C$ inside this annulus gives a definite correction
$\epsilon$ unless $B_{n}$ is $\delta(\epsilon)$-close to a quasidisc
uniform in terms of $\beta$, where $\delta(\epsilon)$ goes to $0$ with
$\epsilon$.
If this correction
occurs, it will be factored into $\overline{s}_{2}$ and consequently into all 
other components of the critical symbol. So in this case we indeed get the
second possibility allowed by Lemma~\ref{lem:26np,1}.

Thus, we can assume that $B_{n}$ normalized so that its boundary passes through
$1$ and $-1$ is $\delta$-close to a $K(\beta)$-quasidisc where
$\delta>0$ can be
specified.  It
remains to show that this last possibility implies the first part of the
alternative claimed in Lemma~\ref{lem:26np,1}. The region $B_{n}$ is the
preimage of $B_{n'}$ by a quadratic mapping composed with a univalent
transformation with bounded distortion. If $B_{n}$ is a uniform quasidisc, it
contains a round disc centered at $0$ of radius comparable to the length of
$B_{n}\cap {\bf R}$. This property carries over to $B_{n'}$. On the other hand,
again by virtue of $B_{n}$ being a uniform quasidisc, it is contained in a
round disc centered at $0$ with radius comparable to the length of the real
section of $B_{n}$. So, the modulus between the complex boxes $B_{n}$
and $B_{n'}$ is $\log v + C(\beta)$. Thus, the component $s_{1}$ of
any separation symbol can be as large as the lemma claims. But then
$A'$ always has at least the same modulus, so $s_{3}$ and $s_{4}$ can
be chosen with the same size. The lemma follows.
\end{proof}

 \begin{lem}\label{lem:26np,2}
Consider a type I complex box mapping $\phi$ of rank $n-1$. Assume that $\phi$ 
arises from some $\phi_{-1}$ in a box inducing step.
Suppose that $\phi$
shows a non-close return which results in a mapping $\phi_{1}$, and $\phi_{1}$
and also shows a non-close return whereby which gives $\phi_{2}$ after a box 
inducing step. Suppose also that in both cases the critical value falls into an
immediate preimage.   Let $v$ be the
ratio $|B_{n}|/|B_{n'}|$.   Let $\beta$ denote the separation norm of
$\phi$. There exists a function $V(\beta)>0$ so that
if $v < V(\beta)$, then at
least one possibility occurs:

\begin{itemize}
\item
the separation norm of $\phi_{1}$  is at least  $-K(beta)\log v$ where
$K>0$ is an
absolute constant, or
\item
the separation symbol for immediate preimages of $\phi_{2}$ at least $\beta + 
\epsilon$ where
$\epsilon>0$ is an absolute constant.
\end{itemize}
\end{lem}
\begin{proof}
We ask the reader to return to the proof of Lemma~\ref{lem:28np,1} in the case
of immediate preimages.   For $\phi_{2}$, the annulus $\overline{A}_{2}$ is 
constructed as the
preimage by the central branch $\psi_{1}$ of $A'(B)\cup A_{3}(B)$ where $B$ is
the postcritical domain of $\phi_{1}$. In the estimates, the modulus of 
$A'(B)\cup A_{3}(B)$ is taken to be equal to the sum of moduli of its
components. These nesting components are separated by the boundary of $B$ which
is the preimage of the boundary of $B_{n}$ by the extension of the branch
defined of $B$ which is univalent and its distortion on the boundary of $B$ is 
bounded in terms of $\beta_{0}$. Assume in addition (as will be
verified later) 
that moduli of both $A'(B)$ and $A_{3}(B)$ are bounded from below
uniformly in terms of $\beta$. Then,
Lemma~\ref{lem:26na,1} implies that either the nesting involves a positive
correction $\epsilon$, or the boundary of $B_{n}$ after normalization 
is $\delta$-close to a uniform quasidisc. Then we proceed as in the proof of 
Lemma~\ref{lem:26np,1} to prove that the first part of the alternative holds.
On the other hand, if the nesting gives a correction, the same correction will 
appear in all components of the symbol, and the second possibility occur. 
So we must show that 
the moduli of $A'(B)$ and $A_{3}(B)$ are both bounded away from $0$.
This is clear for 
$A'(B)$ which is the preimage of the annulus between $B_{n}$ and $B_{n'}$ of
modulus at least $\beta/2$. Next, we have to see how $A_{3}(B)$ arises from
the previous inducing step. Denote the postcritical domain of $\phi_{0}$ by
$b$. Then, $A_{3}(B)$ is the preimage by $\psi_{n-1}$ of $A_{2}(b)$. So it is
sufficient to show that $s_{2}(b)$ is a bounded away from $0$ proportion of 
$\beta$. Recall that $b$ is an immediate preimage, and the inspection of 
separation bounds given by Lemmas~\ref{lem:26np,1} and~\ref{lem:26np,2} for 
immediate preimages shows that $\lambda_{1}(b) \leq \beta/4$.      
\end{proof}

We will now conclude the proof of Theorem D. Consider a sequence of type I
complex box mappings $\phi_{i}$ as defined by the hypothesis of Theorem D.
Observe that if the first possibility occurs in Lemma~\ref{lem:26np,1} or
Lemma~\ref{lem:26np,2} for $\phi = \phi_{j}$, the estimate claimed by Theorem D
follows for $\phi_{j}$ with a uniform $C$.  Take an index $j$. First, look for
$j_{0}$ defined as the largest index $i$ not exceeding $j$ for which the
estimate of Theorem D follows with this $C$. Clearly,   we are done if we show
that from $j_{0}$ to $j$ the moduli grow at a definite linear rate. By
Lemma~\ref{lem:26np,1} we know that each close return causes the modulus to
grow by a constant (that is because now only the second possibility can occur.)
Thus, we are done if we show that a sequence of, say, five consecutive
non-close returns  also increases the modulus. Unless the fourth or fifth of 
those
mappings has the critical value fall into an immediate domain, we are done by
Lemma~\ref{lem:26np,1}. Otherwise, we also get an increase by
Lemma~\ref{lem:26np,2}.

\section{From unimodal to box mappings}
We will need to delve a little more deeply in the inducing
construction of~\cite{yours}. So we begin by recalling certain crucial
constructions.  

\subsection{The basics of inducing}
\paragraph{Extendability.}
The whole fragment on extendability is a repetition of arguments
of~\cite{yours} put in a somewhat different language.
The idea of extendability is to provide a condition which will imply
bounded distortion of branches, at least a fair part of all branches,
and will automatically reproduce itself by the inducing construction. 
What is required is the ``metric extendability'' condition which we
state following~\cite{kus}:
\begin{defi}\label{defi:14fp,2}
A diffeomorphism $g$ with a non-positive Schwarzian derivative defined
on an interval $(a,b)$ is said to be $\epsilon$-extendable if there is
a larger interval $(c,d)\supset (a,b)$ and an extension $\tilde{g}
\supset g$ such that $\tilde{g}$ is still a diffeomorphism with a
non-positive Schwarzian, and
\[ \frac{\tilde{g}(c) - \tilde{g}(a)}{\tilde{g}(d) - \tilde{g}(a)}
\frac{\tilde{g}(d) - \tilde{g}(b)}{\tilde{g}(c) - \tilde{g}(b)} >
\epsilon \; .\]
\end{defi}

We will call the interval $(c,d)$ from Definition~\ref{defi:14fp,2} the
{\em collar} of extendability, and its image by $g$ the {\em margin}
of extendability.

For branches of our box mappings Schwarzian derivative is non-positive
by definition. We
can therefore use the powerful ``real K\"{o}be lemma''
(see~\cite{guke}). This will say that
\begin{fact}\label{fa:fp16,1}
If a monotone branch or its restriction satisfies
$\epsilon$-extendability, the distortion of that branch measured as
the maximum logarithm of the ratio of derivatives taken at two points
is bounded as uniform function of $\epsilon$.
\end{fact}

A folding branch can also be said to be
$\epsilon$-extendable provided that it is the composition of
$x\rightarrow (x-1/2)^{2}$ with an $\epsilon$-extendable
diffeomorphism. The margin of extendability is equal to the margin
determined for this diffeomorphism. The collar is the preimage of the
margin by the complete branch.

\paragraph{Induced box maps.}
Given a mapping from class $\cal F$, we will consider its induced box
maps, that is real box mappings in sense of Definition~\ref{defi:1,1}
whose branches are iterations of $f$. 

A real box mapping will be called {\em suitable} if there is a
symmetric neighborhood of the critical point which is mapped by the
central branch inside itself. If a suitable box mapping is induced,
then this neighborhood must be a restrictive interval of the
underlying $f$. 

\paragraph{General inducing process.}
We proceed to describe the general inducing process introduced
in~\cite{yours}. This will work for any S-unimodal $f$
representable as a power law composed with a bounded distortion
diffeomorphism which is either renormalizable or whose critical orbit
is recurrent.  The main features of this process are as follows:
\begin{enumerate}
\item
Any map in the sequence is defined except on a set of points whose
forward orbits avoid neighborhoods of the critical point.
\item
All rank $0$ branches show a uniform margin of extendability which
is independent of the original map $f$
as well as of the place in the construction. This naturally implies 
$\epsilon$-extendability with a uniform $\epsilon$.
\item
For any rank $0$ branch, its collar of extendability is contained in
the smallest box containing the domain of this branch, with the
obvious exception of the branches adjacent to the boundary of a box
whose collars stick out to one side. Also, except
for the central branch the collar of extendability does not contain
the critical point. 
\item
If a map constructed in the process is not full, it is still of type
I, and any branch of
positive rank $k$ extends with the margin $B_{k'}$.
\end{enumerate}

We will now describe the process. This part is essentially a
summary of~\cite{yours}. The reader can refer there for more details.

\paragraph{The beginning.}
A full induced map exists which satisfies the extendability
properties. In most cases, one just takes the first return map of $f$
into its fundamental inducing domain. A problem may occur since the
central branch can be arbitrarily short which prevents us from
asserting any uniform margin of extendability for the branch next to
the critical one. This problem is taken care of by more inducing as
described in the section 2.2 of~\cite{kus}.  %caution 
The margin of extendability established for this induced map will
remain in force throughout the construction.

We now suppose that an induced box mapping mapping $\phi$, not
necessarily of type I, is given which is not
suitable, whose critical value is in the domain of definition of $\phi$, 
and which satisfies our postulates. We will show how to get the next
induced mapping. From the procedure, the postulates will be also
satisfied. We distinguish a few cases. 
\paragraph{The basic case.}
The basic case occurs if the critical value falls into the domain 
of a monotone rank $0$ branch. In the basic case, the construction
of the new induced mapping proceeds as follows. Define $\phi'$ to be
$\phi$ with the central branch replaced with the identity. Then
compose $\phi$ with $\phi'$. This, first of all, gives you the new
central branch which is of rank $0$.  Overall, we now have a box
mapping $\phi'_{1}$ of complicated structure.
Now consider
the extendability of rank $0$ branches of $\phi'_{1}$. The only ones of 
questionable extendability are preimages of rank $0$ branches of
$\phi'$ by the central branch of $\phi$. This is because the critical
value of $\phi$ may have entered their collars of extendability.
However, the new central branch is uniformly extendable (compare the
definition of extendability for folding branches.) To regain the
extendability of all rank $0$ branches of $\phi_{1}$, we apply the
process of {\em boundary refinement} to branches of $\phi'$ that are
images of non-extendable branches of $\phi_{1}$. The process is
described in~\cite{yours}. For us. it suffices to say that the
boundary refinement involves composing monotone branches of $\phi'$
with $\phi'$ and allows us shrink the collar of extendability of the
branch adjacent to the postcritical domain so that it no longer
contains the critical value. The easy proof is provided in~\cite{yours}.
This phenomenon is based on the fact that the boundary refinement
increases the number of iterations on branches. The boundary of the
maximal possible margin of extendability then gets into a vicinity of
the repelling fixed point $q$ of $f$, and as the number of iterations
increases, this point will be repelled from $q$, thus increasing the
margin. After the adjustment of $\phi'$ by boundary refinement 
followed again by replacing
the central branch of $\phi$ with its composition with $\phi'$, we get
a mapping $\phi_{1}$ with all branches of rank $0$ uniformly
extendable. 

This is followed by a filling-in process. The objective is to obtain a
type I mapping. To this end, monotone branches of ranks between $0$ and
$n+1$ have to be refined. Prepare $\phi'_{1}$ by replacing its central
branch with the identity map. Then take every branch of $\phi_{1}$ of
rank strictly between $0$ and $1$ and replace it with the composition
of this branch with $\phi'_{1}$. Again, some non-extendable branches
of rank $0$ may appear which are preimages of rank $0$ monotone
branches of $\phi'_{1}$. Like in the previous step, one then returns
to $\phi'_{1}$ to boundary-refine the preimages of those troublesome
branches and that takes care of the problem. This gives us $\phi_{2}$.
The mapping $\phi_{2}$ still has branches of ranks bigger than $0$,
though less than $n+1$, but the set occupied by their domains has
shrunk. Then, exactly the same filling-in step is performed with
$\phi_{2}$ instead of $\phi_{1}$. In the limit of filling-in a type I
induced map is regained. Checking the properties of this mapping, we
notice that compared with the domain of definition of $\phi$,
this limit map is not defined on the Cantor set points which forever
stay in domains of ranks between $0$ and $n+1$. But the orbits of
these points forever avoid $B_{n+1}$. Other postulates are easily
satisfied from the construction.  

\paragraph{The box case.}
The box case occurs when the critical value of $\phi$ is
found in a monotone branch of positive rank. We then follow the {\em
standard inducing step} defined earlier on branches of positive rank.
The only difference is that monotone branches of rank $0$ are adjusted
by boundary refinement to ensure their uniform extendability. 

\paragraph{Close returns.}
The remaining case is when the critical value falls into the central
domain. This again follows the process described in the {\em standard
inducing step} with the exception that some monotone branches of
$\phi'$. One first constructs the mapping $\tilde{\phi}$ in the same
way it was described in the standard inducing step. The difference may
occur is some branches of this map of rank $0$ do not have standard
extendability. That is helped by going back to $\phi'$ and
boundary-refining some of its monotone branches of rank $0$.  Then a
dichotomy occurs for $\tilde{\phi}$ since it can show either a box or
a basic return. We simply follow the appropriate step as described above.  

\paragraph{The refined inducing process.}
In addition to the general inducing process just described we will need
yet another inducing procedure, which we call {\em refined inducing
process}.   
The refined inducing also allows type II mappings in the sense of
Definition~\ref{defi:9mp,2} and
switches back forth between the two cases. Let us assume that an
induced map $\phi$ of type I or II and rank $n$ is given. We first
consider the case when $\phi$
is like a map coming from the general inducing process, i.e. it is of
type I and satisfies the properties true for maps obtained in the
general inducing. 

\subparagraph{Inducing on a type I map.}
Let us first assume no close return. The first filling is
done in the usual way. A map $\phi'$ is built by replacing the central
branch with the identity, and then the central branch of $\phi$ is
replaced with its composition with $\phi'$. If that leads to
non-extendable rank $0$ branches, we return to $\phi'$ to apply the
boundary-refinement, and then compose the central branch with this
modified $\phi'$. In what we get all monotone branches of positive
rank have rank $n$, while the new central branch has rank $n$ or $0$
depending on whether the basic or box case occurred. In the basic
case, we exactly follow the general inducing step (there is no
difference between general and refined inducing in the box case.)
Otherwise, we stop at
this stage.  The only thing we do is to adjoin the new central domain to
the box structure as $B_{n+1}$. This gives a type II mapping which
differs from the outcome of the general step by a lack of filling-in.

In a close return, the mapping $\tilde{\phi}$ is obtained as in the
general step, and then the process just described for non-close
returns is used. 
We observe that this modified step  preserves all properties
inductively claimed for the general step.

\subparagraph{Inducing on a type II mapping.}
The procedure depends on whether or not the range of the real central branch
covers the critical point. If so,  we start with a filling-in of
all monotone branches of rank $n'$. This is a general observation that there
always is a passage from type II to type I by filling-in of monotone
branches of rank $n'$ to obtain only monotone branches of rank $n$.
After that, we follow the refined inducing step for type I maps just
described. 

If the image of the central branch does not cover the critical point,
we follow the general step as in the refined step for type I mappings,
and again skip
the filling-in stage. In this case, if the rank of  $\phi$
was $B_{n}$ and, then after this step all branches will be
of rank $n'$. We will thus remove the old $B_{n}$ from the box structure and
replace it with the new central domain. This means that this
particular step does not increase the rank.       

This completes the description of the refined inducing.

\paragraph{Boxes shrink uniformly.}
\subparagraph{Decay in terms of return time.}
We start with a general lemma.
\begin{lem}\label{lem:27gp,1}
For a map $f\in{\cal F}$, let $k$ denote the maximum depth of almost
parabolic points with periods less than the return time of the
restrictive interval if $f$ is renormalizable, or $\infty$ if not. 
For any $\delta > 0$ there exists an integer $m(\delta,k)$ with the
property that if the stopping time on the central branch of a type I induced
map in the sequence constructed from $f$ the general inducing process
is at least $m(\delta,k)$,
then the size of any domain of positive rank, as well as of the
central domain, are less than $\delta$. The function $m(\delta,k)$ is 
independent of $f$. Conversely, for
any integer $m$ there is a number $\delta(m)>0$ so that if the return
time of the central branch is no more than $m$, the length of the central
domain is greater than $\delta(m)$. 
\end{lem}
\begin{proof}
We first prove the bound from below on the rate of decay. 
We begin by noting that $\cal F$ is a normal family in the $C^{2,1}$
topology. Indeed, all members of this family are in the form
$h_{f}(z-1/2)^{2}$. Diffeomorphisms $h_{f}$ are of negative Schwarzian
derivative and uniformly $\epsilon$-extendable by the principle
of~\cite{miszczu}. It is a well known fact the Schwarzian derivative
of an $\epsilon$-extendable iterate of a one-dimensional map with
finitely many polynomial-type singularities is bounded from below
uniformly in terms of $\epsilon$ (see a proof of a very similar
estimate in~\cite{harm}.) Thus the normality follows.

Now, proceed by contradiction and consider a limit $g$ of maps from
$\cal F$ which have increasing return times on the central branch
while the sizes of positive rank or folding domains remain bounded
away from $0$. One easily sees that $g$ has a homterval, i.e., an
interval on which all iterations of $g$ are monotone. By the general
result of~\cite{blyu}, $g$ must have a non-repelling, thus neutral
cycle. We also notice that $g$ continues not to expand cross-ratios,
thus by~\cite{preston} this neutral orbit is unique and the critical
point is in the immediate basin of one point, say $p$. Now carry out
the inducing process for $g$. The critical point and $p$ will always
stay together in the central branch, since branches in the inducing
construction are separated either by preimages of the fixed point.
Next, it is a property of the construction that for any branch,
no intermediate images enter the central domain. This statement is
verified by induction. The easiest way to see it is by reviewing the
description we provided asking which branches can change when we
change $f$ on its central domain only. By induction, for $\phi$ the
central branch is the only one. Now, observe that when
constructing $\phi'$ to be precomposed with other branches, we always
replace the central branch by the identity, which eliminates the
dependence except for the new central branch. From this observation it
follows that return times on the central branch in the inducing
process for $g$ cannot jump the period of $p$. Thus, after finitely
many steps an induced map is obtained which exhibits a close return
(which must be low, i.e. the image of the real central branch does not
cover the critical point). Now, if we take a map $f$ from the sequence which
allegedly contradicts the claim of the lemma which is very
close to $g$ in the $C^{2}$ topology, the construction is
conducted in the same way for $f$, since the course of the
construction only depends on where the critical value falls. The map
$f$ will show a low return, but will recover from it after a large
number of steps.  By Lemma~\ref{lem:10ma,1} the central branch of this
map is uniformly $\epsilon$-extendable.
Since  it takes a long time for the critical value to escape the
central domain, and this time can be made arbitrarily large by
choosing $f$ close enough to $g$, we can obtain a map $f$ with an
almost parabolic point of arbitrary depth, contradiction.  

Finally, we prove the bound from above on the rate of decay of boxes. 
This follows immediately by induction. By construction, each central
branch defined on $B_{n}$ is a composition of the quadratic polynomial with a
diffeomorphism mapping on $B_{n-1}$. In the case of a close return,
this should be applied to mappings $\phi_{i}$ described in the
standard inducing step. But the derivative of this diffeomorphism is
bounded from above in terms of the return time. So, by each step the
box shrinks only by a bounded factor in terms of the return time, and
the number of inducing steps is certainly no more than the return
time.  
\end{proof}

\subparagraph{Decay in terms of the rank.}
\begin{lem}\label{lem:10ma,1}
If $\phi$ is a type I or II induced box mapping of rank $n$ derived from some
$f\in {\cal F}$, then 
\[ \frac{|B_{n}|}{|B_{n'}|} \leq 1- \epsilon\]
where $\epsilon$ is a n absolute constant independent of $f$. 
\end{lem}
\begin{proof} 
If $\phi$  is full, the ratio is indeed bounded away from $1$ since a
fixed proportion of $B_{0}$ is occupied by the domains of two branches
with return time $2$. In a sequence of box mappings this ratio remains
bounded away from $1$ by Fact~\ref{fa:9mp,1}.   
\end{proof}

For type I mappings, we know that monotone branches of rank $r > 0$ are
extendable to $B_{r'}$, and the central branch is extendable to
$B_{(r-1)'}$, also in type II mappings.    
Note that  it implies $\epsilon$-extendability with $\epsilon$ an
absolute constant. 
\subsection{Main proposition}
\begin{prop}\label{prop:28gp,1}
For every $\delta>0$, there are a fixed integer $N$ and $\beta_{0}>0$,
both independent of
the dynamics, for which
the following holds. Given a mapping from $\cal F$ a mixture of
general and refined inducing gives an induced map $\phi$ which
satisfies at least one of these conditions:
\begin{itemize}
\item
$\phi$ is of type I, suitable and of rank less than $N$,
\item
$\phi$ is of type I, satisfies the starting condition with norm
$\delta$, its central branch is $\tau$-extendable, and $\delta\leq
\delta(\tau)$ in the sense of the hypothesis of Fact~\ref{fa:8mp,1},
\item
$\phi$ has a structure of a complex box mapping with separation norm
greater than $\beta_{0}$ and the rank less than $N$.  
\end{itemize}

If the depths of almost parabolic points with periods less than the
return time of the restrictive interval (taken as $\infty$ in the
non-renormalizable case) are bounded by $k$, the bound $N$ on the rank
can be replaced by a bound the return time on the central branch by 
a function $N(k)$.
\end{prop}

\paragraph{Taking care of the basic case.}
Let $\delta>0$ be given small enough so that if for a full map induced
from some $f\in {\cal F}$ the central domain and all domains of
positive rank are shorter that $\delta$, then this map satisfies the
assumption of Fact~\ref{fa:8mp,1}. This is certainly possible as all domains of
positive rank are separated from the boundary of $B_{0}$ by the
branches of return time $2$ adjacent to the boundary of $B_{0}$ and
the central branch is uniformly extendable (the $\tau$ is uniformly
bounded away from $0$ in Fact~\ref{fa:8mp,1}.) 

Then, for any $f$ we construct an induced
map $\varphi(\delta)$ according to the following procedure. We follow
the general inducing process until {\em one step before} the central
domain and all domains of positive rank become smaller than $\delta$.
By Lemma~\ref{lem:27gp,1}, the stopping time on the central branch of
$\varphi(\delta)$ is uniformly bounded in terms of $\delta$, hence by
an absolute constant.  We can be
prevented from being able to construct $\phi(\delta)$ by earlier
hitting a suitable map. That, however, means that the return time of
the restrictive interval is uniformly bounded in terms of $\delta$, so
are done as far the proof of Proposition~\ref{prop:28gp,1} is concerned.

We take $\varphi(\delta)$ as the beginning point for further inducing. 
Note that if in the general or refined inducing a full map is ever
derived from $\varphi(\delta)$, we are already done with the proof of
Proposition~\ref{prop:28gp,1}, since the starting condition holds.
Thus, we can assume that general and refined inducing on
$\varphi(\delta)$ encounters exclusively box cases.  

\paragraph{A way to control stopping times.}
As another step in the proof of Proposition~\ref{prop:28gp,1}, we
offer a convenient way of controlling stopping times on the central
branch in terms of the rank. This will let us derive the estimates
claimed in terms of return times in the absence of almost parabolic
points of great depth.   Let us label the induced mappings which follow
$\varphi(\delta)$ in the general inducing procedure
as $\varphi_{0} := \varphi(\delta)$, $\varphi_{1}$ is the next
one, and so on. 
\begin{lem}\label{lem:28gp,2} 
Under the hypothesis and in notations of
Proposition~\ref{prop:28gp,1}, there exists a
sequence of integer constants $C(j,k)$, independent of $f$, such that if
the stopping time of the central branch of $\varphi_{j}$ exceeds
$C(j,k)$, then the starting condition holds.
\end{lem}
\begin{proof}
We proceed by induction with respect to $j$, starting from $j=0$ where
the lemma is obvious. Since the stopping time on the central branch is
bounded by $C(j,k)$, the size of the central domain is bounded away from
$0$. Otherwise, the derivative on the central branch would follow to
be uniformly small. But, by Lemma~\ref{lem:27gp,1} if the stopping time
on the central branch of $\varphi_{j+1}$ becomes large,
the size of its central domain will become small as a uniform
function depending also on $k$. Observe, however, that the smallness
of the central domain of
$\varphi_{j+1}$ relative to the central domain of $\varphi_{j}$
determines box ratio in the starting condition for
$\varphi_{j+1}$. Also, the extendability of the central
branch  is uniformly 
bounded in terms of $j$ by Lemma~\ref{lem:10ma,1}. The lemma follows.
\end{proof}  

The meaning of Lemma~\ref{lem:28gp,2} is that in the absence of almost
parabolic point of great depth, the requirement regarding the stopping
time of the
central branch being bounded in a certain way can be
replaced by a simpler condition that the number of box steps leading
to the suitable map from $\varphi(\delta)$ is appropriately bounded. 
\subsection{Finding a hole structure}
We consider the sequence $(\varphi_{k})$ of consecutive box mappings
obtained from $\varphi_{0} := \varphi(\delta)$ in the course of the
refined inducing. The objective of this section is to show that for
some $k$ which is bounded independently of everything else in the
construction, a uniformly bounded hole structure exists which extends
$\varphi_{k}$ as a complex box mapping in the sense of
Definition~\ref{defi:5mp,1}.

\paragraph{The case of multiple type II maps.}
We will prove the following lemma:
\begin{lem}\label{lem:28gp,5}
Consider some $\varphi_{m}$ of rank $n$. There is a function $k(n)$
such that if the mappings $\varphi_{m+1},\ldots,\varphi_{m+k(n)}$ are
all of type II, then $\varphi_{m+k(n)}$ has a hole structure
which makes it a complex box mapping. 
The separation norm of the hole structure is bounded away from $0$
depending solely on $n$.  
\end{lem}
\begin{proof}
If a sequence of type II mappings occurs, that means that the image of
the central branch consistently fails to cover the critical point. The
rank of all branches is fixed and equal to $n$. The central domain
shrinks at least exponentially fast with steps of the construction at
a uniform rate by Lemma~\ref{lem:10ma,1}. 
Thus, the ratio of the length of the central domain of $\varphi_{m+k}$
to the length of $B_{n}$ is bounded by a function of $k$ which also
depends on $n$, and for a fixed $n$ goes to $0$ as $k$ goes to
infinity.

This means that we will be done if we show that a small enough value
of this ratio ensures the existence of a bounded hole structure.   
We choose two symmetrical circular arcs which intersect the line in
the endpoints of $B_{n}$ at angles $\pi/4$ to be the boundary of the
box.($\alpha$
will be chosen in a moment, right now assume $\alpha<\pi/2$).
We take the
preimages of the box by the monotone branches of rank $n$. They are
contained in similar circular sectors circumscribed on their domains
by Poincar\'{e} metric considerations given in~\cite{miszczu}.    
Now, the range of the central branch is not too short compared to the
length of $B_{n}$ since it at
least covers one branch adjacent to the boundary of the box (otherwise
we would be in the basic case). We claim that the
domain of that external branch constitutes a proportion of the
box bounded depending on the return time of the central branch.
Indeed, one first notices that the return time of the external branch
is less than the return time of the central branch. This follows
inductively from the construction. But the range of the external
domain is always the whole fundamental inducing domain, so the domain
of this branch cannot be too short. On the other hand, the size of the
box is uniformly bounded away from $0$ by Lemma~\ref{lem:27gp,1}.
To obtain
the preimage of the box by the complex continuation of the central
branch, we write the central branch as $h(z-1/2)^{2}$. The preimage by
$h$ is easy to handle, since it will be contained in a similar
circular sector circumscribed on the real preimage. Since the
distortion of $h$ is again bounded in terms of $n$, the real range of
$(x-1/2)^{2}$ on the central branch will cover a proportion of the
entire $h^{-1}(B_{n})$ which is bounded away from $0$ uniformly in
terms of $n$. Thus, the preimage of the complex box by the central
branch will be contained in a star-shaped region which meets the real
line at the angle of $\pi/4$ and is contained in a rectangle built on
the central domain of modulus bounded in terms of $n$.
    
It follows that this preimage will be contained below in the complex box
if the middle domain is sufficiently small. As far as the separation norm is
concerned, elementary geometrical considerations show that it is
bounded away from $0$ in terms of $n$.
\end{proof}
\paragraph{Formation of type I mappings.}
We consider then same sequence $\varphi_{k}$ and we now analyze the
cases when the image of the central branch covers the critical point.
Those are exactly the situations which lead to type I maps in the
refined inducing.  

\subparagraph{A tool for constructing complex box mappings.}  
The reader is warned that the notations used in this technical
fragment are ``local'' and should not be confused with symbols having
fixed meaning in the rest of the paper.
The construction of complex box mappings (choice of a bounded hole
structure ) in the remaining cases will be based on the technical work
of~\cite{limo}. We begin with a lemma which appears there without
proof.
\begin{lem}\label{lem:2ha,1}
Consider a quadratic polynomial $\psi$ normalized so that
$\psi(1)=\psi(-1)=-1$, $\psi'(0)=0$ and $\psi(0)=a\in (-1,1)$.
The claim is that if $a<\frac{1}{2}$, then $\psi^{-1}(D(0,1))$ is
strictly convex.
\end{lem}
\begin{proof}
This is an elementary, but somewhat complicated computation. We will
use an analytic approach by proving that the image of the tangent line
to $\partial\psi^{-1}(D(0,1))$ at any point is locally strictly
outside of $D(0,1)$ except for the point of tangency.   
We represent points in $D(0,1)$ in polar coordinates $(r,\phi)$
centered at $a$ so that $\phi(1)=0$, while for points in preimage we
will use similar polar coordinates $r',\phi'$. By school geometry we
find that the boundary of $D(0,1)$ is given by
\[ (r+a\cos\phi)^{2} + a^{2}\sin^{2}\phi = 1\; .\]
By symmetry, we restrict our considerations to $\phi\in[0,\pi]$. We
can then change the parameter to $t=-a\cos\phi$, which allows us to
express $r$ as a function of $t$ for boundary points, namely
\begin{equation}\label{equ:2ha,1}
r(t) = \sqrt{1-a^{2}+t^{2}} + t\; .
\end{equation}

Now consider the tangent line at the preimage of $(r(t),t)$ by $\psi$.
By conformality, it is perpendicular to the radius joining to $0$, so
it can be represented as the set of points $(r',\phi')$
\[ r'=\frac{\sqrt{r(t)}}{\cos^{2}(\theta/2)}\:,\; \phi' = \frac{\pi}{2}
- \frac{\phi}{2} + \frac{\theta}{2} \]
where $\theta$ ranges from $-\pi$ to $\pi$. The image of this line is
given by $(\hat{r}(\theta),\phi-\theta)$ where 
\[  \hat{r}(\theta) = \frac{2r(t)}{1 + \cos\theta} \; .\]
Let us introduce a new variable $t(\theta) := -a\cos(\phi-\theta)$ so
that $t=t(0)$. Our task is to prove that
\begin{equation}\label{equ:2ha,2}
r(t(\theta)) < \hat{r}(\theta)
\end{equation} for values of $\theta$ is some
punctured neighborhood of $0$. This will be achieved by comparing the
second derivatives with respect to $\theta$ at $\theta=0$.  
By the formula~[\ref{equ:2ha,1}]
\[ r(t(\theta)) = \sqrt{1-a^{2}+t^(\theta)^{2}} + t(\theta)\; .\]
The second derivative at $\theta=0$ is
\[ \frac{1-a^{2}}{(\sqrt{1-a^{2}+t^(\theta)^{2}})^{3}}(a^{2} - t^{2}) - t -
\frac{t^{2}}{\sqrt{1-a^{2}+t^(\theta)^{2}}}\; = \]
\[ =-\sqrt{1-a^{2}+t^(\theta)^{2}}-t +
\frac{1-a^{2}}{(\sqrt{1-a^{2}+t^(\theta)^{2}})^{3}}\; .\] 
The second derivative of the right-hand side of the desirable
inequality~[\ref{equ:2ha,2}] is more easily computed as
\[ \frac{\sqrt{1-a^{2}+t^(\theta)^{2}}+t}{2} \; .\]
Thus, the proof of the estimate~[\ref{equ:2ha,2}], as well as the
entire lemma, requires showing that
\[ \frac{3}{2}(\sqrt{1-a^{2}+t^(\theta)^{2}}+t) -
\frac{1-a^{2}}{(\sqrt{1-a^{2}+t^(\theta)^{2}})^{3}}  > 0\]
for $|a|<1/2$ and $|t|\leq |a|$.  For a fixed $a$, the value of this
expression increases with $t$. So, we only check $t=-a$ which reduces
to 
\[ \frac{3}{2}(1-a) - (1 - a^{2}) > 0\]
which indeed is positive except when $a\in [1/2,1]$.
\end{proof}
\subparagraph{The main lemma.}
Now we make preparations to prove another lemma, which is
essentially Lemma 8.2 of~\cite{limo}.
Consider three nested intervals $I_{1}\subset I_{0}\subset I_{-1}$
with the common midpoint at $1/2$. Suppose that a map $\psi$ is
defined on $I_{1}$ which has the form $h(x-1/2)^{2}$ where $h$ is
a polynomial diffeomorphism onto $I_{-1}$ with non-positive Schwarzian
derivative. We can think $\psi$ as the central branch of a generalized
box mapping.  We denote 
\[ \alpha := \frac{|I_{1}|}{|I_{0}|}\]. Next, if $0 < \theta\leq \pi/2$
we define $D(\theta)$ to be the union of two regions symmetrical with
respect to the real axis. The upper region is defined as the
intersection of the upper half plane with the disk centered in the
lower $\Re = 1/2$ axis so that its boundary crosses the real line at
the endpoints of $I_{0}$ making angles $\theta$ with the line. So,
$D(\pi/2)$ is the disk having $I_{0}$ as diameter.
\begin{lem}\label{lem:7np,1}
In notations introduced above, if the following conditions are satisfied:
\begin{itemize}
\item
$\psi$ maps the boundary of $I_{-1}$ into the boundary of $I_{0}$,
\item
the image of the central branch contains the critical point,
\item 
the critical value inside $I_{0}$, but not inside $I_{1}$,
\item
the distance from the critical value to the boundary of $I_{0}$ is no
more than the (Hausdorff) distance between $I_{-1}$ and $I_{0}$,
\end{itemize}
then $\psi^{-1}(D(\theta))$ is contained in $D(\pi/2)$ and the
vertical strip based on $I_{1}$.
Furthermore, for every $\alpha < 1$ there is a choice of
$0<\theta(\alpha)<\pi/2$ so that 
\[\psi^{-1}(D(\theta(\alpha)))\subset D(\theta(\alpha))\] with a
modulus at least $K(\alpha)$, and $\psi^{-1}(D(\theta(\alpha)))$ is
contained in the
intersection of two convex angles with vertices at the endpoints of
$I_{1}$ both with measures less than $\pi-K(\alpha)$. 
Here, $K(\alpha)$ is a continuous positive function. 
\end{lem}
\begin{proof}
By symmetry, we can assume that the critical
value, denoted here by $c$, is on the left of $1/2$. 
Then $t$ denotes the right endpoint of $I_{0}$, and $t'$ is the other
endpoint of $I_{0}$. Furthermore, $x$ means the right endpoint of
$B_{n-1}$. By assumption, $h$ extends to the
range $(t',x)$. 
To get the information about the preimages of points 
$t,t',c,x$ one considers their cross-ratio 
\[ C = \frac{(x-t)(c-t')}{(x-c)(t-t')} \geq \frac{1+\alpha}{4}\]
where we used the assumption about the position of the critical value
relative $I_{0}$ and $I_{-1}$. The cross ratio will not be decreased by
$h^{-1}$. In addition, one knows that $h^{-1}$ will map the disk of
diameter $I_{0}$
inside the disk of  diameter $h^{-1}(B_{n})$ by the Poincar\'{e}
metric argument of~\cite{miszczu}.  As a consequence of the
non-contracting property of the cross-ratio, we get
\begin{equation}\label{equ:1hp,1}
\frac{h^{-1}(c) - h^{-1}(t')}{h^{-1}(t) - h^{-1}(t')} <
\frac{1+\alpha}{4}\; .
\end{equation}
When we pull back the disk based on $h^{-1}(I_{0})$, we will get a
figure which intersects the real axis along $I_{1}$. Notice that by
the estimate~[\ref{equ:1hp,1}] and Lemma~\ref{lem:2ha,1}, the preimage
will be convex, thus necessarily contained in the vertical strip based
on $I_{1}$. Its height in the imaginary direction is
\begin{equation}\label{equ:1hp,2}
 \frac{|I_{1}|}{2}\sqrt{\frac{h^{-1}(t) - h^{-1}(c)}
{h^{-1}(c) - h^{-1}(t')}} <
 \frac{|I_{1}|}{2}\sqrt{\frac{3-\alpha}{1+\alpha}}\; ,
\end{equation}

 where we used the estimate~[\ref{equ:1hp,1}] in the last inequality.    
Clearly, 
\[ \psi^{-1}(D(\pi/2))\]
is contained in the disk of this radius
centered at $1/2$. To prove that 
\[ \psi^{-1}(D(\pi/2))\subset D(\pi/2)\; , \] in view of the 
relation~[\ref{equ:1hp,2}] we need
\begin{equation}\label{equ:2ha,3}
 \alpha\sqrt{\frac{3-\alpha}{1+\alpha}} < 1\; 
\end{equation}
By calculus one readily checks that this indeed is the case when
$\alpha<1$.
To prove the uniformity statements, we first observe that 
\[ \psi^{-1}(D(\theta)) \subset \psi^{-1}(D(\pi/2)) \; \]
for every $\theta < pi/2$. Since [\ref{equ:2ha,3}] is a sharp
inequality, for every $\alpha<1$ there is some range of values of
$\theta$ below $\pi/2$ for which $\psi^{-1}(D(\theta))\subset
D(\theta)$ with some space in between. We only need to check the
existence of the angular sectors. For the intersection of 
$\psi^{-1}(D(\theta))$ with a narrow strip around the real axis, such
sectors will exist, since the boundary intersects the real line at
angles $\theta$ and is uniformly smooth. Outside of this narrow strip,
even $\psi^{-1}(D(\pi/2))$ is contained in some angular sector by its
strict convexity.
\end{proof}

The assumption of
extendability to the next larger box is always satisfied in our
construction. 
\subparagraph{The case when there is no close return.}
We now return to our construction and usual notations. We consider a
map $\varphi_{k}$, type II and of rank $n$,
 whose central branch covers the critical
point, but without a close return. Then:
\begin{lem}\label{lem:2ha,2} 
Either the Hausdorff distance from $B_{n}$ to $B_{n-1}$ exceeds the
Hausdorff distance from $B_{n-1}$ to $B_{n-2}$, or $\varphi_{k}$ has a
hole structure uniformly bounded in terms of $n$.
\end{lem}
\begin{proof}
Suppose the condition on the Hausdorff distances fails. We choose the
box around $B_{n}$ and the hole around $B_{n+1}$ by
Lemma~\ref{lem:7np,1}. Observe that the quantity $\alpha$ which plays
a role in that Lemma is bounded away from $1$ by
Lemma~\ref{lem:10ma,1}. The box is
then pulled back by these monotone branches and its preimages are
inside similar figures built on the domains of branches by the usual 
Poincar\'{e} metric argument of~\cite{miszczu}. For those monotone
branches, the desired bounds follow immediately.
\end{proof}

\subparagraph{The case when the close return occurs.}
Finally, we have to deal with the case when the image of the central
branch of $\varphi_{k}$, which as always can be assumed of type II and
some rank $n$, covers the critical point, but also makes it a close
return.
\begin{lem}\label{lem:2ha,3}
In the situation described above at least one the three statements
holds:
\begin{itemize}
\item
the starting condition holds as claimed by
Proposition~\ref{prop:28gp,1}, 
\item
 a hole structure uniformly
bounded in terms of $k$ exists for an induced map $\varphi'$ obtained
from $\varphi_{k}$ and from which $\varphi_{k+1}$ can be extracted by
more inducing,
\item
 the Hausdorff distance from
$B_{n+1}$ to $B_{n}$ exceeds the Hausdorff distance from $B_{n}$ to
$B_{n-2}$.
\end{itemize}
\end{lem}
\begin{proof}
Suppose that the condition on the Hausdorff distances does not hold.
We look into the description of the inducing process in the close
case, and seek an opportunity to apply Lemma~\ref{lem:7np,1}. We build
a sequence of temporary boxes which are preimages of $B_{n}$ by the
central branch until for one of them, say $B_{n+m}$, the critical
value escapes. First, we remark by Lemma~\ref{lem:28gp,2} that $m$ is
uniformly bounded in terms of $k$, or the starting condition holds and
we have  nothing more to prove. 
To construct $\varphi'$ we just fill-in all branches of positive rank
so that they map onto $B_{n+m-1}$. The map $\varphi'$ is of type II
and no longer shows a close return, and indeed $\varphi_{k+1}$ can be
obtained from $\varphi'$ by a step of the refined inducing construction.   
The hole structure for $\varphi'$ is obtained by repeating the
argument of Lemma~\ref{lem:2ha,2} with the additional information that
the central branch is extendable to a margin equal to the whole
$B_{n-2}$.
\end{proof}

\paragraph{Proof of Proposition~\ref{prop:28gp,1}.}
This is just a summary of
the work done in this section. We claim that we have proved that
either a map with a box structure can be obtained from
$\varphi(\delta)$ in a uniformly bounded number of steps of the
refined inducing process, or the Proposition~\ref{prop:28gp,1} holds
anyway. Since Lemmas~\ref{lem:28gp,5},~\ref{lem:2ha,2}
and~\ref{lem:2ha,3} all provide uniform bounds for the hole
structures in terms of $k$ or the rank which is bounded in terms of
$k$, it follows that the hole structure is bounded or the starting
condition holds anyway.  It also follows that by
Lemma~\ref{lem:28gp,2} that if the
inducing fails within this bounded number of steps because of a
suitable map being reached, then the stopping time on the central
branch of the suitable map is bounded, hence
Proposition~\ref{prop:28gp,1} again follows. 

So, we need to prove that claim. If the claim fails, then by
Lemma~\ref{lem:28gp,5} the situations in which the image of the
central branch covers the critical point have to occur with definite
frequency. That is, we can pick a function $m(k)$ independent of other
elements of the construction which goes to infinity with $k$ such that
among $\varphi_{1},\ldots,\varphi_{k}$ the situation in which the
critical point is covered by the image of the central branch occurs at
least $m(k)$ times. But each time that happens, we are able to
conclude by Lemmas~\ref{lem:2ha,2} and~\ref{lem:2ha,3} that the
Hausdorff distance between more deeply nested boxes is more than
between shallower boxes. Initially, for $\varphi_{0}$ whose  rank
was $n$, the $B_{n}$ distance between and $B_{n-1}$ was a fixed
proportion of the
diameter of $B_{n}$  So
only a bounded number of boxes can be nested inside $B_{n-1}$ with
fixed space between any two of them (or at least between every other
pair in the situation of Lemma~\ref{lem:2ha,3}.) 
So we have a bound on
the value of $m(k)$, thus on $k$. This proof of the
claim is a generalization of the reasoning used in~\cite{limo}.

The claim concludes the proof of Proposition~\ref{prop:28gp,1}.    
\subsection{Proofs of main theorems}
\paragraph{Proof of Theorem~\ref{t:1}.}
By~\cite{yours} an S-unimodal non-renormalizable mapping without
attracting or indifferent periodic points is expansion-inducing
provided that box ratios shrink to $0$. Thus, Theorem~\ref{t:1}
follows directly.   

\paragraph{Proof of Theorem~\ref{th:9ma,2}.}
We will prove this theorem by showing that if the return time on the
restrictive interval is sufficiently large, it follows that for the suitable
map obtained in the inducing construction of~\cite{yours} which is of type I
and rank $n$ the ratio $|B_{n'}|/|B_{n}|$. Namely, given a $D$ in
Theorem~\ref{th:9ma,2}, this ratio should be at least $D$. Next,
specify $\delta$ in Proposition~\ref{prop:28gp,1} equal to that $D$.
Of the three outcomes of Proposition~\ref{prop:28gp,1} the integer $N$
will give us the minimum return time in Theorem~\ref{th:9ma,2}. In
other cases the starting condition either holds, in which case we are
done, or we get a hole structure with a uniform separation norm. Then
we use Theorem B. We learn that after a bounded number of box steps,
which by Lemma~\ref{lem:28gp,2} means either a starting condition again
or a uniformly bounded return time on the central branch, we get the
starting condition with norm $\delta$. So, in all cases either the
starting condition holds with norm $\delta$ and can only improve, or
the return time on the central branch is uniformly bounded.
Theorem~\ref{th:9ma,2} follows.


\begin{thebibliography}{99}
\bibitem{blyu} Blokh, A. \& Lyubich, M.: {\em Non-existence of
wandering intervals and structure of topological attractors for one
dimensional dynamical systems}, Erg. Th. \& Dyn. Sys. {\bf 9} (1989), 751-758
\bibitem{brahu} Branner, B. \& Hubbard, J.H.: {\em The iteration of
cubic polynomials, Part II: patterns and parapatterns}, Acta Math. {\bf 169}
(1992), pp. 229-325
\bibitem{dohu} Douady, A. \& Hubbard, J.H.: {\em On the dynamics of
polynomial-like mappings}, Ann. Sci. Ec. Norm. Sup. (Paris) {\bf 18} (1985),
pp. 287-343 
\bibitem{gracz} Graczyk, J.: {\em Ph.D. Thesis}, Mathematics Department of
Warsaw University (1990); also: {\em Dynamics of non-degenerate upper
maps}, preprint of Queen's University at Kingston, Canada (1991)  
\bibitem{harm} Graczyk, J. \& \'{S}wi\c{a}tek, 
G.:{\em Critical circle maps near 
bifurcation}, Stony Brook IMS preprint (1991), Proposition 2
\bibitem{guke} Guckenheimer, J.: {\em Limit sets of S-unimodal maps
with zero entropy}, Commun. Math. Phys. {\bf 110} (1987), pp. 655-659
\bibitem{gujo} Guckenheimer, J. \& Johnson, S.: {\em Distortion of
S-unimodal maps}, Annals of Math. {\bf 132}, 71-130 (1990) 
\bibitem{hokel} Hofbauer, F. \& Keller, G.: {\em Some remarks about
recent results on S-unimodal maps,}
Annales de l'Institut Henri Poincare, Physique Theorique {\bf 53},
413-425 (1990) 
\bibitem{invmes} Jakobson, M.: {\em: Absolutely continuous invariant
measures for one-parameter families of one-dimensional maps}, Commun.
Math. Phys. {\bf 81} (1981), pp. 39-88 
\bibitem{yours}
Jakobson, M. \& \'{S}wi\c{a}tek, G.: {\em Metric properties of
non-renormalizable S-unimodal maps}, preprint IHES no. {\bf
IHES/M/91/16} (1991)
\bibitem{kus}
Jakobson, M. \& \'{S}wi\c{a}tek, G.: {\em Quasisymmetric conjugacies between
unimodal maps}, Stony Brook preprint {\bf 16} (1991)
\bibitem{nowike}
Keller, G. \& Nowicki, T.: {\em Fibonacci maps revisited}, manuscript (1992)
\bibitem{lehvi} Lehto, O.\& Virtanen, K.:{\em Quasikonforme
Abbildungen}, Springer-Verlag, Berlin-Heidelberg-New York (1965)
\bibitem{opo} Lyubich. M: {\em Combinatorics, geometry and attractors
of quasi-quadratic maps}, Stony Brook preprint {\bf 18}(1992) 
\bibitem{limo} Lyubich, M,, Milnor, J.: {\em The dynamics of the
Fibonacci polynomial}, Jour. of the AMS {\bf 6} (1993), pp. 425-457
\bibitem{ma} Martens, M.: {\em Ph.D. thesis}, Math. Department
of Delf University of Technology (1990); also: IMS preprint {\bf 17} (1992)
\bibitem{miln} Milnor, J.: {The Yoccoz theorem on local connectivity
of Julia sets. A proof with pictures.}, class notes, Stony Brook, (1991-92)
\bibitem{preston} Preston, C.: {\em Iterates of maps on an interval},
Lecture Notes 
in Mathematics, Vol. 999. Berlin,Heidelberg,New York: Springer (1983)
\bibitem{miszczu} 
Sullivan, D.: {\em Bounds, quadratic differentials and renormalization
conjectures}, to appear in American Mathematical Society Centennial 
Publications, Volume {\bf 2}, American Mathematical Society,
Providence, R.I. (1991)
\bibitem{hyde} \'{S}wi\c{a}tek, G.: {\em Hyperbolicity is dense in the real
quadratic family}, manuscript, version of May 18, 1993  
\bibitem{szkop} Teichm\"{u}ller, O.: {\em Untersuchungen \"{u}ber konforme
und quasikonforme Abbildung}, Deutsche Mathematik {\bf 3}, pp. 621-678 (1938) 
\bibitem{yoc} Yoccoz, J.-C.: unpublished results
\bibitem{tanver} Veerman, J.J.P. \& Tangerman F. M.: {\em Scalings in
circle maps (I)}, Commun. in Math. Phys. {\bf 134} (1990), pp. 89-107 
\end{thebibliography}
\end{document}